\documentclass[11pt]{article}
\usepackage{amsmath,amssymb,graphicx,float,epsf}
\usepackage{csquotes}
\usepackage{appendix}
\usepackage{multirow}
\def\RR{\hbox{I\kern-.2em\hbox{R}}}
\newcommand{\qed}{\hbox to 0pt{}\hfill$\rlap{$\sqcap$}\sqcup$ \vspace{3mm}}
\numberwithin{equation}{section}
\restylefloat{figure}
\usepackage{graphicx} 
\usepackage{amsmath}
\usepackage{amsthm}
\usepackage{amsfonts}
\usepackage[caption = false]{subfig}
\usepackage{wrapfig}
\usepackage{enumerate}
\usepackage{enumitem}
\usepackage{array}
\usepackage{caption}
\usepackage[margin=1cm]{caption}
\usepackage{authblk}

\usepackage[utf8]{inputenc}
\usepackage{geometry}
\vbadness=\maxdimen

\newtheorem{Th}{\textbf{Theorem}}

\usepackage{graphicx}
\usepackage{tikz}
\usetikzlibrary{shapes.geometric, arrows}
\tikzstyle{rect} = [draw, rectangle, fill=blue!20, text width=6em, text centered, minimum height=2em]
\tikzstyle{elli} = [draw, ellipse, fill=red!20, minimum height=2em]
\tikzstyle{circ} = [draw, circle, fill=white!20, minimum width=8pt, inner sep=5pt]
\tikzstyle{diam} = [draw, diamond, fill=white!20, text width=6em, text badly centered, inner sep=0pt]
\tikzstyle{line} = [draw, -latex']
\usepackage{a4wide} 
\usepackage{authblk} 

\usepackage{varioref}
\usepackage{hyperref}
\hypersetup{colorlinks=true, linkcolor=violet, citecolor=cyan, urlcolor=teal} 

\usepackage{tikz}
\usetikzlibrary{shapes.geometric, arrows}
\usetikzlibrary{calc}
\tikzstyle{cs} = [rectangle,minimum width=1cm, minimum height=1cm, text centered,
draw=blue,fill=blue!5, thick]
\tikzstyle{ce} = [rectangle,minimum width=1cm, minimum height=1cm, text centered,draw=orange,fill=orange!5,thick]
\tikzstyle{ci} = [rectangle,minimum width=1cm, minimum height=1cm, text centered,
draw=red!70!red,fill=red!5, thick]
\tikzstyle{cr} = [rectangle,minimum width=1cm, minimum height=1cm, text centered,draw=green,fill=green!5,thick]
\tikzstyle{cv} = [rectangle,minimum width=1cm, minimum height=1cm, text centered,draw=green,fill=green!5,thick]
\tikzstyle{ct} = [rectangle,minimum width=1cm, minimum height=1cm, text centered,draw=green,fill=green!5,thick]
\tikzstyle{arrow} = [thick,->,>=stealth]
\graphicspath{{Fig3/}}

\date{}

\begin{document}

\title{Modeling H1N1 Influenza Transmission and Control: Epidemic Theory Insights Across Mexico, Italy, and South Africa}

	\author[1]{\small Md Kamrujjaman\thanks{Corresponding author  Email: kamrujjaman@du.ac.bd}}
	\author[1]{\small  Kazi Mehedi Mohammad\thanks{Email: mehedimim.me@gmail.com }}

	\affil[1]{\footnotesize Department of Mathematics, University of Dhaka, Dhaka 1000, Bangladesh}
	
		\affil[+]{\small  Both authors contributed equally to this work.}

\pagestyle{myheadings}
\thispagestyle{plain}		
	
	\maketitle
	\vspace{-1.0cm}
	\noindent\rule{6.35in}{0.02in}\\
	{\bf Abstract.}\\
This study incorporates mathematical analysis, focusing on developing theories and conducting numerical simulations of Influenza virus transmission using real-world data.   The terms in the equations introduce parameters which are determined by fitting the model for matching clinical data sets using non-linear least-square method. The purpose is to determine the wave trend, critical illness factors and forecast for Influenza in national levels of  Mexico, Italy, and South Africa and to investigate the effectiveness of control policy and making some suggestions of alternative control policies. 	Data were extracted from the following sources: published literature, surveillance, unpublished reports, and websites of Centres For Disease Control and Prevention (CDC) \cite{CDC}, Natality report of U.S. clinics and World Health Organization (WHO) Influenza Disease Dashboard \cite{WHO}. 
We included total 120 weeks data (which are calculated as per thousand) from October 01, 2020 to March 31, 2023 \cite{CDC}, throughout this study.
Numerical and sensitivity analysis are carried out to determine some prevent strategies. The objectives of local and global sensitivity analysis is to determine the dominating parameters and effective correlation with $\mathcal{R}_0$. We presented data fitting, Latin hypercube sampling, sensitivity indices, Partial Rank Correlation Coefficient, p-value, estimation of the nature of $\mathcal{R}_0$ from available data to show validation of the model with these counties. The aim is to determine optimal control strategies with drug administration schemes, treatments which represent the efficacy of drug inhabiting viral production and preventing new infections, minimizes the systematic cost based on the percentage effect of the drug. The existence and the uniqueness of the optimal pair are discussed and an optimal system is obtained that is solved numerically by a compatible method. We are interested in finding the qualitative behaviour of threshold which determines whether the disease dies out or persist in the population. Finally, we present series of numerical examples and the effect of different parameters on the compartments to verify theoretical results. 
\\[2mm]	
\noindent\rule{6.35in}{0.02in}

	\vspace{10mm}	
	
	
	\clearpage
	\noindent {\bf Highlights:}\\
	\begin{enumerate}
		\item 	Developing a 
		optimal control framework for a data-driven mathematical model of Influenza H1N1.
		\item 	Analyzing the wave trend, critical illness factors, and forecasting Influenza at the national level in Mexico, Italy,  and South Africa, while also assessing the effectiveness of control policies.
		\item Estimating parameters using clinical datasets through the non-linear least squares method.
		\item  The goal of local and global sensitivity analysis is to determine the factors that are prominent and how well they correlate with the fundamental reproduction number, $\mathcal{R}_0$.
			\item Using data fitting, sensitivity indices, LHS,  PRCC, p-values, and estimating the nature of threshold to validate the model across these countries.
			\item 	The aim is to determine optimal control strategies with drug administration schemes, treatments which represent the efficacy of drug inhabiting viral production and preventing new infections, minimizes the systematic cost based on the percentage effect of the drug.
	 
	\end{enumerate}	
	\clearpage
	
\section{Introduction}
Epidemics and pandemics have ravaged populations time and time again throughout human history, frequently radically changing the path of history and causing entire civilizations to fall. Because they shed light on the implementation of practical and efficient disease-control strategies, mathematical models that clarify the dynamics of infectious illnesses are essential to public health. Influenza-like illness, characterized by unpleasant bodily symptoms, typically lasts for 2–7 days, resolving on its own in the absence of complications from other disorders \cite{Stability Bound-5, Stability Bound-8, Stability Bound-9}. The severity of the illness is significantly influenced by the individual's immune system. The production of airborne particles and aerosols containing viruses, which occur during regular speaking and breathing activities, provides the basis for respiratory transmission of the virus. Sneezing, particularly if accompanied by increased mucus production due to infection, is a more efficient means of expelling viruses from the nasal cavity \cite{Stability Bound-14, Stability Bound-17, Stability Bound-18}.

On June 11, 2009, the World Health Organization (WHO) proclaimed the unique strain of influenza A (H1N1) to be a pandemic after it was discovered in both Mexico and the United States \cite{WHO}. During the incubation phase, individuals infected with A (H1N1) are not transmissible and do not exhibit any symptoms. The A (H1N1) virus can incubate for one to four days, and a patient's infectious phase can extend from the day before symptoms appear to seven days after they do. The influenza A (H1N1) virus causes a variety of symptoms, such as cough, nausea, diarrhea, fever, chills, headache, sore throat, muscle aches, runny nose, shortness of breath, and joint pain \cite{Stability Bound-14,Stability Bound-13}.
Although many elderly individuals had previously been exposed to a similar H1N1 virus strain and were protected by their antibodies, the pandemic resulted in the deaths of over 200,000 individuals worldwide  \cite{WHO}.

Numerous researchers have turned to the well-established SEIR model or its adaptations to elucidate how individuals transition through various compartments, representing different disease phases across the entire population over time \cite{Stability Bound-14, Stability Bound-21}. These models help in the creation and evaluation of public health initiatives by offering insightful information about the dynamics of disease transmission.
Annual influenza vaccination stands as a crucial preventive measure against influenza. However, the virus's rapid mutation poses a significant challenge. A vaccine designed for one year may not confer effective protection in the following year \cite{Stability Bound-23}. Furthermore, the phenomenon of antigenic drift can occur after the vaccine's development, diminishing its efficacy. Consequently, outbreaks of influenza, particularly among high-risk populations, become more likely.

Controlling the spread of influenza and other infectious diseases requires not just vaccination but also additional preventive measures. These include avoiding close contact with sick individuals, practicing proper respiratory hygiene by covering coughs and sneezes, and maintaining good hand hygiene through regular handwashing \cite{ CDC, WHO,Stability Bound-14}. By adopting a multi-faceted approach that combines vaccination with these preventive measures, it is feasible to lessen the effects of infectious diseases and the load they place on the populace in impacted areas.
In most pandemic scenarios, an exponential growth curve is initially observed, followed by a gradual flattening, which corresponds to reducing the epidemic peak \cite{Marcheva Book}. Understanding the transmission dynamics of a new infectious disease outbreak becomes crucial for flattening the curve, especially in the absence of established treatments or vaccinations. Mathematical models serve as essential tools for public health authorities in making decisions to optimize control measures, as these decisions heavily rely on short- and long-term predictions provided by these models \cite{Stability Bound-4}.

During seasonal Influenza virus outbreaks, epidemiological models have played a significant role in many countries \cite{WHO}. Countries like Italy and Mexico have widely adopted these models to gain insights into recent situations, assess the impact of outbreak control measures, explore alternative interventions, and provide guidance to other similar settings. Furthermore, in resource-limited settings, prediction models with multiple features would be invaluable to healthcare workers for patient monitoring \cite{CDC,Marcheva Book}.

We introduce a model for Influenza virus transmission, explicitly incorporating vaccinated and treatment subpopulations. We take it that there is an equal opportunity for interaction between the exposed and infected compartments for every healthy host.  In this study, we consider a  six compartments’ potential Susceptible-Vaccinated-Exposed
Infectious-Treatment-Removal (SVEIRT) mathematical model. Appendix \ref{appen} contains the main model as well as basic reproduction number, endemic equilibrium (EE), and disease-free equilibrium (DFE) results, which are essential for our subsequent discussions in the paper.
\subsubsection*{Research Gap in Recent Studies}
In recent studies \cite{Stability Bound-5,  Stability Bound-8, Stability Bound-9,Stability Bound-6, Stability Bound-7} identifies a gap in studying simple SEIR epidemic model analytical findings. This study extends this gap by implementing not only the vaccination and treatment strategies but also six compartmental optimal control strategies analytically for suggesting effective control scenarios. Stability analysis along with data fitting for specific region are represented in \cite{ Stability Bound-17, Stability Bound-18, Stability Bound-16}, but this study addresses this gap by numerical simulation with multinational clinical real data by non-linear least square methods. Moreover, local and global sensitivity analysis along with PRCC are performed to show effective correlation with parameters with $\mathcal{R}_0.$ Only optimal control with analytical studies are reflected in \cite{Optimal Control-1, Optimal Control-2, Optimal Control-3, Optimal Control-4, Optimal Control-5, Optimal Control-6}. Building upon this, this study delves deeper investigation with qualitative behaviour with the threshold $\mathcal{R}_0$ using several methods such as sensitivity, relative bias to determine the existence and persistence of the disease.

\begin{figure}[H]
	\centering  
{\includegraphics[width=3.5 in]{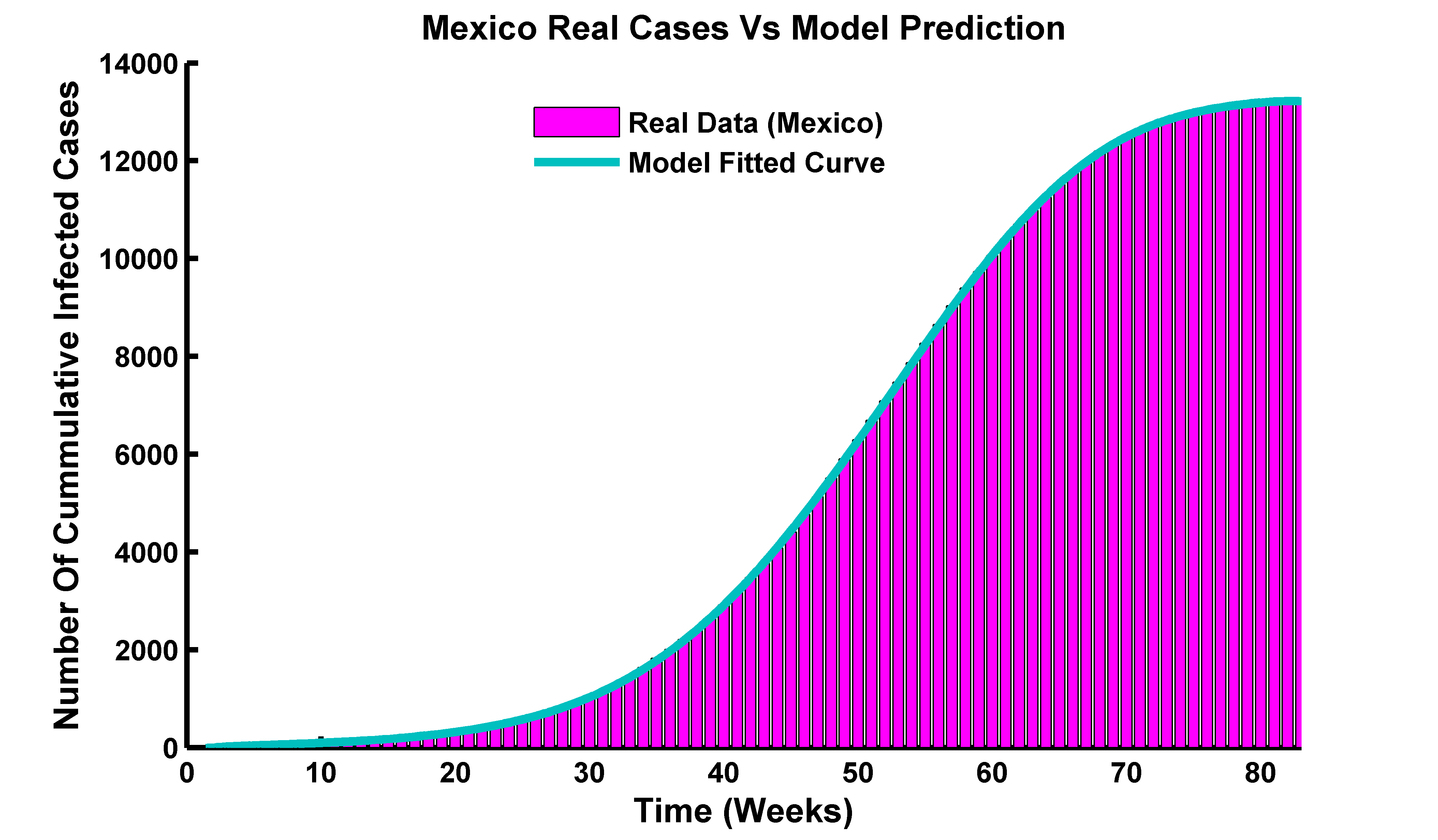}}\\
	\caption{Weekly total of Mexico data cases for 85 weeks beginning on October 15, 2021.}
	\label{FF1}
\end{figure}

Based on existing literature, the aims of this paper are:
\begin{itemize}
	\item Performing optimal control analysis for managing disease outbreaks.
	\item Implementing vaccination and treatment strategies to mitigate disease spread and outbreaks.
	\item Various statistical tools are employed in sensitivity analysis to quantify the correlation between input parameters and model outputs, while considering the influence of other variables.
	\item Analyze the optimal control model theoretically, looking at the persistence and existence of its solutions.
\end{itemize}
In relation to the objectives, we conducted a qualitative analysis of the optimal control problem, examining properties such as stability, existence, uniqueness, and sensitivity, without relying solely on numerical computations. Our goal was to identify optimal control strategies involving drug administration schemes and treatments aimed at inhibiting viral production and preventing new infections, with the aim of minimizing systemic costs based on the drug's effectiveness. We also studied data from multiple countries and compared model solutions to pave the way for preventing disease outbreaks. In order to verify theoretical results, we also presented a number of computational solutions and evaluated the effects of various factors on the compartments. Additionally, we evaluated important parameters since, when working with sickness data, accurate parameter values are essential for trustworthy quantitative forecasts across time intervals. Figure \ref{FF1} provides an immediate justification for the results that have been displayed.

The structure of this paper is as follows:  Optimal Control strategy is elaborately discussed in Section  \ref{Section-Optimal Control}. Parameter estimation, model validation, assessment of $\mathcal{R}_0$ and effective reproduction numbers are elaborated in Section \ref{Section-Parameter-Estimation-model-validation}. Then,  relative influence on $\mathcal{R}_0$ are presented in Sections  \ref{Subsection-Relative-Influence}. In Section \ref{Section-Sensitivity Analysis}, the sensitivity index is used to present the sensitivity analysis, LHS and PRCC procedure, also relative bias of $\mathcal{R}_0$ are depicted. Theoretical outcomes are presented in Section \ref{Section-Numerical-Simulation} to show behaviour of compartments of model for various parameters. After that, validity of PRCC result with phase plane is summarized. This part concludes with a comparison of the model solution with the real data analysis utilizing case studies from South Africa, Mexico, and Italy. Finally, Section \ref{Section-Concluding-Remarks} gives an overview of the results' analysis and summary. Remark that the  Mathematical Model is discussed  in  Appendix \ref{appen} establishing fixed locations  and figuring out the fundamental reproduction number, $\mathcal{R}_0$ (with control and without control).

\section{Optimal Control for Influenza}\label{Section-Optimal Control}
The ideas of control theory have proven useful in many different contexts. As a method of treatment for the recommended therapeutic agents, too high a dosage might not be preferred by patients, and too low a dosage might not be helpful. The most effective treatment regimens can lower the danger of drug side effects, viral mutation, and expensive medication load \cite{Optimal Control-1}. It is crucial to maintain appropriate drug levels in a bodily compartment in order to avoid issues resulting from the dangerous side effects of the treatment. Our primary objective is to determine the optimal medication dosage to halt the epidemic while lowering the possibility of negative side effects. In this case, the controlling factor is the pharmaceutical input, and one can partially influence a patient's response to medication by being aware of the magnitude of the drug \cite{Optimal Control-5}.\\
Individuals with influenza infection who are asymptomatic are given antiretroviral therapy (ART) to enhance their health. Numerous administration strategies are employed to enhance the quality of life for patients while also preventing the emergence of medication resistance, stopping the spread of new virus strains, minimizing major side effects, and lowering prescription costs. This section's main objective is to provide a mathematical framework that may be utilized to choose the optimal drug delivery strategy for improving patient health, especially in settings with limited resources \cite{Optimal Control-2, Optimal Control-3}.

\subsection{Mathematical Model with Vaccination and Treatment}
Two types of antiviral drugs are commonly used to reduce the viral load and limit the number of infected T cells. One subclass of influenza drugs, known as reverse transcriptase inhibitors (RTIs), stops viral RNA from becoming viral DNA, hence stopping the propagation of new infections. Protease inhibitors (PIs) comprise the other category. These agents impede the assembly of essential viral proteins that have inadvertently been generated by infected host cells. In this sense, RTIs lower the rate of infection of activated CD$4^{+}$T cells, while PIs inhibit the generation of fresh infectious virions. Thus, both medications reduce the spread of the infection. The major goal of this section is to identify the optimal strategy for administering antiviral pharmaceutical therapy to treat influenza infection, with a particular emphasis on maximizing CD$4^{+}$T cell count while minimizing drug toxicity or overall cost.Following the initiation of combined chemotherapy, which consists of a combination of RTIs and PIs, the amount of viral particles generated by an actively infected CD$4^{+}$T cell is reduced.\\
The best control countermeasures to stop the influenza virus from spreading, given that prevention is controlled by offering campaigns in the form of counseling to individuals in the susceptible subpopulation, indicated by $w_1$. Vitamin treatment control, indicated by $w_2$, is administered to members of the $I$ subpopulation in an attempt to hasten healing. The level of successful precaution (mask use, public place disinfection, physical separation to decrease the effective interactions of exposed and infected individuals with vulnerable individuals) is represented by the control function $w_1(t)$. In order to reduce the spread of influenza (seasonal flu), a novel optimal control problem is formulated in this part, which introduces time-dependent controls $w_1(t),\; w_2(t)$, and $w_3(t)$. Wearing a mask and maintaining a minimum 3-foot distance from others is known as physical distancing, or just avoiding crowds. To improve recovery, we first substitute the prior treatment recovery rate with $(1+w_2)\gamma$ and $w_2=1$. The control $(1+w_2)$ denotes the level of therapy (antiviral treatment), represents $100\%$ effective treatment. For the increased treatment of infected patients, $(1+w_3)$ is the control. Regarding $w_1(t)$, which represents the normalized RTI dosage as a function of time, the compartment $S$ will be changed to $(1-w_1(t))S$, which depicts the social distancing effectiveness. The subpopulation $I$ will change to become $(1+w_2(t))I$ and $(1+w_3(t))I$ if we additionally allow $w_2(t)$ and $w_3(t)$ to represent the normalized PI dosage. Taking into account all of them, the system \eqref{new_model}'s optimum control problem is as follows:
\begin{align}\label{op_syseqn}
	\begin{cases}
		\vspace{0.2cm}
		\displaystyle \frac{dS}{dt}=\Lambda-(1-w_1(t))(\beta_1 E+\beta_2 I)S-(\mu+\phi) S\\
		\vspace{0.2cm}
		\displaystyle \frac{dV}{dt}=\phi S-(1-\varepsilon)(\beta_1 E+\beta_2 I)V-\mu V\\
		\vspace{0.2cm}
		\displaystyle \frac{dE}{dt}=(1-w_1(t))(\beta_1 E+\beta_2 I)S-(\alpha+\mu) E\\
		\vspace{0.2cm}
		\displaystyle \frac{dI}{dt}= \alpha E+(1-\varepsilon)(\beta_1 E+\beta_2 I)V -(\mu+\delta) I-(1+w_2(t))\gamma_1 I-(1+w_3(t))\gamma I\\
		\vspace{0.2cm}
		\displaystyle \frac{dR}{dt}=(1+w_3(t))\gamma I-\mu R\\
		\displaystyle \frac{dT}{dt}=(1+w_2(t))\gamma_1 I -\mu T
	\end{cases}
\end{align}
and the non-negative initial conditions are: $S_0=S(0) \geq 0\;,\;V_0=V(0) \geq 0\;,\;E_0=E(0) \geq 0\;,\;I_0=I(0) \geq 0,\;R_0=R(0) \geq 0\;$ and $T_0=T(0) \geq 0$.\\
and $S(t),\;V(t),\;E(t),\;I(t),\;R(t),\;T(t)$ are free at final time $T_f.$\\
The ideal control parameters $0\leq w_1(t),\;w_2(t),\;w_3(t)\leq 1$ represents the proportion of therapeutic effects on CD$4^{+}$T cell interaction with virus and virions generated by infected cells. The model's \eqref{new_model} schematic representation is presented in Figure \ref{fignewmodel} (see Appendix \ref{appen}). It is also remarked that the qualitative analysis of the model \eqref{op_syseqn} is presented in Appendix \ref{qaulity}.

\subsection{Optimal Control Problem}
Our primary objective is to minimize the systemic cost related to the administered chemotherapy (RTIs and PIs) while optimizing the benefit of CD$4^{+}$T cell count. The objective functional for this study is defined as:
\begin{align}\label{opti-obj-function}
	J(w_1, w_2, w_3)=\int_0^{t_f} \left[a_1 E+ a_2 I+(a_3 w_1^2+a_4 w_2^2+a_5 w_3^2)\right]dt
\end{align}
In order to prevent the spread of influenza infection, this goal function must be reduced. The benefit based on CD$4^{+}$T cells is indicated by the variables $E(t)$ and $I(t)$ in the model \eqref{op_syseqn}, while the other terms represent the systematic costs of the treatment approaches. An increase in CD$4^{+}$T cells and a reduction in systemic medication costs underpin the benefits of therapy. Additionally, the intended weight on benefit and cost is represented by the positive constants $a_3,\;a_4,\;a_5$, and the severity of the drug's adverse effects is indicated by the constants $w_1^2,\;w_2^2,\;w_3^2$. Given the possibility of a nonlinear relationship between the effects of therapy on CD$4^{+}$T cells, the cost function is expected to be nonlinear, which is why a quadratic cost function was used. Because influenza can multiply and develop resistance to treatment after a certain amount of time, we set a requirement for treatment time $t\in[0,T_f]$, limited treatment window, that tracks the worldwide consequences of these events.Treatment may also have negative side effects, and the longer the course of treatment, the more severe the side effects become. Time $t=0$ denotes the start of the treatment, while time $t=T_f$ denotes its termination.\\
where $t_f$ denotes the end time, the weight constants $a_1$, $a_2$ are linked to the infectious population, and the weight constants $a_3$, $a_4$, and $a_5$ are linked to the cost of interventions. Here, we consider all control efforts, $w_1(t)$, $w_2(t)$, and $w_3(t)$, to be Lebesgue measurable on the interval $[0, t_f]$ and bounded. Because of the population's non-linear involvement, we take into consideration the quadratic goal functional in this case. In the case when $w_1=w_2=w_3=1$, a control effort of $100 \%$ is applied to reduce effective contact and enhance treatment, respectively. In the event that $w_1(t)=0$, $w_2(t)=0$, and $w_3(t)=0$, no control effort is applied to apply \cite{Optimal Control-6}. Finding the best control variables, $w_1^*(t)$, $w_2^*(t)$, and $w_3^*(t)$, for decreased effective contact and increased treatment, respectively, is the primary objective, such that,
\begin{align*}
	J(w_1^*,w_2^*,w_3^*)=\min_{\mathcal{D}} J(w_1,w_2,w_3)
\end{align*}
and $\mathcal{D}$, the control effort set, is determined by,
\begin{align*}
	\mathcal{D}=&\{(w_1,w_2,w_3) : [0,t_f]\rightarrow[0,1],(w_1,w_2,w_3)\;\; \text{is a Lebesgue measurable and }\\
	& 0\leq w_1(t),w_2(t),w_3(t)\leq 1\}
\end{align*}
In order to get the best control pair, we look for $w_1^*(t)$, $w_2^*(t)$, and $w_3^*(t)$, such that
\begin{align*}
	J(w_1^*,w_2^*,w_3^*)=\min\limits_{w_1,w_2,w_3\in \mathcal{D}}J(w_1,w_2,w_3),
\end{align*}
subject to the state constraints described in \eqref{op_syseqn}.\\
The main objectives of this challenge are to prove that the optimality system is unique, as well as to establish the existence and definition of the optimal control.

\subsection{The Presence and Identifiability of an Optimal Control Pair}
The model (\ref{op_syseqn}) provides us with $N'(t)\leq \Lambda-\mu N$.
Then there exists $M \in \mathbb{R}^+$ such that,
\begin{align*}
	N(t)\leq N_0 e^{-\mu t}+\frac{\Lambda}{\mu}(1-e^{-\mu t})=M, \; t\in [0, T]
\end{align*}
when, $t\rightarrow\infty$ then $\displaystyle N(t)\leq\frac{\Lambda}{\mu}$. Since, $N(t)=S(t)+V(t)+E(t)+I(t)+R(t)+T(t)$ and there is a solution to the system (\ref{op_syseqn}) if the state variables $S(t),\;V(t),\;E(t),\;I(t),\;R(t)$, and $T(t)$ are limited above. Now, utilizing this situation, we analyze the system using Fleming and Rishel's proposed theory to look into the possibility of an optimal control \cite{Optimal Control-4}.

\begin{Th}\label{op_th_exist}
	If each of the functions $\widetilde{F_i}$, for $i=1,2,\cdots,n$ and the partial derivatives $\frac{\partial\widetilde{F_i}}{\partial x_j}$ for $i=1,2,\cdots,n$. are continuous in $\mathcal{R}^{n+1}$ space, then there exists a unique solution $(x_1=\phi_1(t),\phi_2(t),\cdots,x_n=\phi_n(t))$ of the system of differential equations,
	\begin{align*}
		x_i'=\widetilde{F_i}(t,x_1,x_2,\cdots,x_n) \;\;\text{for}\;\; i=1,2,\cdots,n.
	\end{align*}
	with initial conditions $x_i(t_0)=x_i^0$ for $i=1,2,\cdots,n$ and the solution also satisfies the initial conditions. Given the objective functional:
	\begin{align*}
		J(w_1,w_2,w_3)=\int_{0}^{t_f}\left[a_1E+a_2I+\left(a_3w_1^2+a_4w_2^2+a_5w_3^2\right)\right]dt,
	\end{align*}
	where, $\displaystyle W=\{w_1(t),w_2(t),w_3(t)\}$, is piecewise continuous such that
	$0\leq w_1(t),w_2(t),w_3(t)\leq 1$ for all $t\in [0,t_f]$ subject to the system's equations \eqref{op_syseqn} with\\ $S(0)=S_0,\;V(0)=V_0,\;E(0)=E_0,\;I(0)=I_0,\;R(0)=R_0$ and \; $T(0)=T_0$, then a perfect pair of controls $w_1^*, w_2^*, w_3^*$ exists such that 
	\begin{align*}
		J(w_1^*,w_2^*,w_3^*)=\min\{J(w_1,w_2,w_3)|(w_1,w_2,w_3)\in W\}.
	\end{align*}
\end{Th}
\begin{proof}
	The requirements from the existence, uniqueness, positivity, and boundedness of model \eqref{new_model} solutions are followed in order to demonstrate the optimal control's existence and uniqueness, which is  already established, as a consequence, we have to verify them. Let, $\mathbf{f}(t,\mathbf{X},\mathbf{W})$ be on the system's right side \eqref{op_syseqn} for $0\leq t\leq t_f$, where $V\in\mathbb{R}^6$, $\mathbf{W}\in\mathbb{R}^3$ where, $\mathbf{X}=(S,V,E,I,R,T)$ and $\mathbf{W}=(w_1,w_2,w_3)$. According to \cite{Optimal Control-4,Optimal Control-6} the following prerequisites must be met in order for the optimal control to exist:
	\begin{itemize}
		\item[$(i)$] The class of all beginning circumstances that satisfy each state condition and state variable, along with an optimal control variable $w_1,\;w_2,\;w_3$ in the acceptable control set, is not empty. That suggests
		\begin{align*}
			|\mathbf{f}(t,0,0)|\leq C,\;\; |\mathbf{f_X}(t,\mathbf{X},\mathbf{W})|\leq C(1+|\mathbf{W}|),\;\;\text{and }\;\; |\mathbf{f_W}(t,\mathbf{X},\mathbf{W})|\leq C
		\end{align*}
		\item[$(ii)$] The suitable control set $\mathbf{W}$ is closed and convex,
		\item[$(iii)$] The system \eqref{op_syseqn} has continuous right-hand side equations, which are further bounded above by the sum of the bounded control and the state. The best control variables, $w_1,\;w_2,\;w_3$, can be expressed as a linear function with coefficients that depend on the state variables and the passage of time. That makes clear,
		\begin{align*}
			|\mathbf{f}(t,\mathbf{X},\mathbf{W})|=\mathbf{A}(t,\mathbf{X})+\mathbf{G}(t,\mathbf{X})\mathbf{W}
		\end{align*}
		and $\displaystyle |\mathbf{f}(t,\mathbf{X},\mathbf{W})|\leq C_1(1+|\mathbf{X}|+|\mathbf{W}|).$
		\item[$(iv)$] The objective functional integrand $J(w_1,w_2,w_3)$ and $\mathcal{H}$ is convex on the appropriate control set and this is bounded below by the factor $-C_2+C_1|\mathbf{W}|^{\beta}$, where $C_1,\;C_2$ are positive constants and $\beta>1.$
	\end{itemize}
	
	To prove the hypothesis $(i)-(iv)$ for the uniqueness of the optimal control pair, let us consider the system,
	\begin{align}\label{op_th_exist_sys_1}
		\begin{cases}
			\vspace{0.2cm}
			\displaystyle \frac{dS}{dt} &=\widetilde{F_1}(t,S,V,E,I,R,T) \\
			\vspace{0.2cm}
			\displaystyle \frac{dV}{dt} &=\widetilde{F_2}(t,S,V,E,I,R,T) \\
			\vspace{0.2cm}
			\displaystyle \frac{dE}{dt} &=\widetilde{F_3}(t,S,V,E,I,R,T) \\
			\vspace{0.2cm}
			\displaystyle \frac{dI}{dt} &=\widetilde{F_4}(t,S,V,E,I,R,T) \\
			\vspace{0.2cm}
			\displaystyle \frac{dR}{dt} &=\widetilde{F_5}(t,S,V,E,I,R,T) \\
			\displaystyle\frac{dT}{dt} &=\widetilde{F_6}(t,S,V,E,I,R,T)
		\end{cases}
	\end{align}
	where, $\widetilde{F_1},\;\widetilde{F_2},\;\widetilde{F_3},\;\widetilde{F_4},\;\widetilde{F_5}$ and $\widetilde{F_6}$ symbolize the system's right side (\ref{op_syseqn}) as well as a few constants, $c_1$, $c_2$, and $c_3$. Let, $w_1(t)=c_1$,\;$w_2(t)=c_2$ and $w_3(t)=c_3$. The functions $\widetilde{F_i}$ for $i=1,2,\cdots,6,$ must be partial derivatives of linear functions with respect to all state variables and constants. Because of this, the functions are continuous everywhere, as are their partial derivatives. According to Theorem \ref{op_th_exist}, we can thus assert that there is a singular solution $S(t)=\phi_1(t)$, $V(t)=\phi_2(t)$, $E(t)=\phi_3(t)$, $I(t)=\phi_4(t)$, $R(t)=\phi_5(t)$ and $T(t)=\phi_6(t)$, which fulfills the prerequisites. As a result, there is anything in the constituting set of controls and matching state variables.\\
	Now for any three controls $w_1,w_2,w_3 \in \mathbf{W}$ and $\theta_1,\theta_2 \in [0,1]$, $0<\theta_1 w_1+\theta_2 w_2+(1-\theta_1-\theta_3)w_3\leq 1$.
	Hence, the set $\mathbf{W}$ is closed and convex.\\ 
	Now, contrasting \eqref{op_th_exist_sys_1} with \eqref{op_syseqn} we have,
	\begin{align*}
		\begin{cases}
			\vspace{0.2cm}
			\displaystyle \widetilde{F_1} \leq& \mu N-(1-w_1)(\beta_1 E+\beta_2 I)S-\phi S\\
			\vspace{0.2cm}
			\displaystyle \widetilde{F_2} \leq& \phi S-(1-\varepsilon)(\beta_1 E+\beta_2 I)V\\
			\vspace{0.2cm}
			\displaystyle \widetilde{F_3} \leq& (\beta_1 E+\beta_2 I)(1-w_1)S-\alpha E\\
			\vspace{0.2cm}
			\displaystyle \widetilde{F_4} \leq& \alpha E+(1-\varepsilon)(\beta_1 E+\beta_2 I)V-\delta I-(1+w_2)\gamma_1 I-(1+w_3)\gamma I\\
			\vspace{0.2cm}
			\displaystyle \widetilde{F_5} \leq& \gamma(1+w_3)I \\
			\displaystyle \widetilde{F_6} \leq& \gamma_1(1+w_2) I
		\end{cases}
	\end{align*}
	as a matrix formation,
	\begin{align}\label{op_th_exist_sys_2}
		\bar{F}(t,\mathbf{\bar{X}},\mathbf{W}) \leq \bar{m}_1\left(t,\begin{pmatrix}
			S \\ V\\ E\\ I\\ R\\ T
		\end{pmatrix}\right)\mathbf{\bar{X}}(t)+\bar{m}_2\left(\begin{pmatrix}
			S \\ V\\ E\\ I\\ R\\ T
		\end{pmatrix}\right)\mathbf{W}
	\end{align}
	where, 
	\begin{align*}
		\bar{m}_1\left(t,\begin{pmatrix}
			S \\ V\\ E\\ I\\ R\\ T
		\end{pmatrix}\right)=\begin{pmatrix}
			a_{11} & 0 & -\beta_1 S & -\beta_2 S & 0 & 0 \\
			\phi & a_{22} & -\lambda\beta_1V & -\lambda\beta_2V & 0 & 0\\
			(\beta_1 E+\beta_2 I) & 0 & a_{33} & \beta_2S & 0& 0\\
			0 & a_{42} & \alpha+\lambda\beta_1V & a_{44} & 0 &0 \\
			0 & 0 & 0& \gamma & 0 &0\\
			0 & 0 & 0& \gamma_1 & 0 &0
		\end{pmatrix}
	\end{align*}
while, $a_{11}=-(\beta_1E+\beta_2I)-\phi$, $a_{22}=-\lambda(\beta_1E+\beta_2I)$, $a_{33}=-\alpha+\beta_1S$, $a_{42}=\lambda(\beta_1E+\beta_2I)$ and $a_{44}=-(\delta+\gamma+\gamma_1)+\lambda\beta_2V.$\\
	Also
	\begin{align*}
		\bar{m}_2\left(t,\begin{pmatrix}
			S \\ V\\ E\\ I\\ R\\ T
		\end{pmatrix}\right)= \begin{pmatrix}
			(\beta_1 E+\beta_2 I)S \\ 0 \\ -(\beta_1 E+\beta_2 I)S \\ -\gamma_1 I-\gamma I\\ \gamma I\\ \gamma_1 I
		\end{pmatrix}
	\end{align*}
	Here, all the parameters are constant and non-negative. Therefore from (\ref{op_th_exist_sys_2}) we have,
	\begin{align*}
		|\bar{F}(t,\mathbf{\bar{X}},\mathbf{W})| &\leq ||\bar{m}_1|||\mathbf{\bar{X}}|+|S+E+I|\left|(w_1(t),w_2(t),w_3(t))\right|\\
		&\leq q(|\mathbf{\bar{X}}|+\left|(w_1(t),w_2(t),w_3(t))\right|) 
	\end{align*}
	Consequently, a linear function in the state and control variables bounds the right side of the state system \eqref{op_syseqn} \cite{Optimal Control-4}.\\
	Additionally, the integrand 
	\begin{align*}
		\mathcal{L}=a_1 I+a_2 E+(a_3 w_1^2+a_4 w_2^2+a_5 w_3^2)
	\end{align*}
	of the objective function $J$ and $\mathcal{H}$ is convex and satisfies 
	\begin{align*}
		\mathcal{H}(w_1,w_2,w_3)= -C_2+C_1|(w_1,w_2,w_3)|^{\beta}=-C_2+C_1|\mathbf{W}|^{\beta}
	\end{align*}
	where, $w_1 >0\;,\;w_2 >0\;$ and $w_3>0$ and $\beta=2>1$, according to \cite{Optimal Control-4,Optimal Control-6}. This part proves the uniqueness of the solution.\\
	In order to confirm the existence theorem, we represent the system's right-hand side equations \eqref{op_syseqn} as 
	\begin{align*}
		\mathbf{f}(t,\mathbf{X},\mathbf{W})=(f_1,f_2,f_3,f_4,f_5,f_6)
	\end{align*}
	where, $f_1,f_2,f_3,f_4,f_5,f_6$ denotes the right hand side of optimal control model \eqref{op_syseqn} respectively. It is obvious that $\mathbf{f}(t,\mathbf{X},\mathbf{W})$ belongs to class $C^1$ and that $|\mathbf{f}(t,0,0)|=\Lambda$ and we have
	\begin{align*}
		|\mathbf{f_X}(t,\mathbf{X},\mathbf{W})|=\left|\begin{pmatrix}
			a_{11} & 0 & a_{13} & a_{14} & 0 & 0\\
			\phi & a_{22} & -\lambda\beta_1 V& -\lambda\beta_2V & 0& 0\\
			a_{31} & 0 & a_{33} & a_{34} & 0& 0\\
			0 & a_{42} & a_{43} & a_{44} & 0 & 0\\
			0 & 0 & 0 & a_{54} & -\mu & 0\\
			0 & 0& 0& a_{64} & 0& -\mu
		\end{pmatrix}\right|
	\end{align*}
	where, $a_{11}=-(1-w_1)(\beta_1E+\beta_2I)-(\mu+\phi)$, $a_{13}=-(1-w_1)\beta_1S$, $a_{14}=-(1-w_1)\beta_2S$,\\
	$a_{22}=-\lambda(\beta_1E+\beta_2I)-\mu$, $a_{31}=(\beta_1E+\beta_2I)(1-w_1)$, $a_{33}=-(\alpha+\mu)+\beta_1(1-w_1)S$,\\
	$a_{34}=\beta_2(1-w_1)S$, $a_{42}=\lambda(\beta_1E+\beta_2I)$, $a_{43}=\alpha+\lambda\beta_1V$,\\
	$a_{44}=-(\mu+\delta+(1+w_2)\gamma_1+(1+w_3)\gamma)+\lambda\beta_2V$, $a_{54}=\gamma(1+w_3)$, $a_{64}=\gamma_1(1+w_2).$\\
	Moreover,
	\begin{align*}
		|\mathbf{f_W}(t,\mathbf{X},\mathbf{W})|=
		\left|\begin{pmatrix}
			(\beta_1E+\beta_2I)S & 0&0\\
			0&0&0\\
			-(\beta_1E+\beta_2I)S & 0&0\\
			0 & -\gamma_1 I & -\gamma I\\
			0&0& \gamma I\\
			0 & \gamma_1I & 0
		\end{pmatrix}\right|
	\end{align*}
	Since, each of the compartments $S,V,E,I,R,T$ are bounded, then a constant $C$ exists such that 
	\begin{align*}
		|\mathbf{f}(t,0,0)|\leq C,\;\; |\mathbf{f_X}(t,\mathbf{X},\mathbf{W})|\leq C(1+|\mathbf{W}|)\;\;\text{and}\;\; |\mathbf{f_W}(t,\mathbf{X},\mathbf{W})|\leq C
	\end{align*}
	Definition applies, thus $\mathbf{W}$ is closed. Using any controls provided by the parameters $w_1,\;w_2,\;w_3\in \mathbf{W}$ and $\theta_1,\;\theta_2\in[0,1].$ Then
	\begin{align*}
		\theta_1w_1+\theta_2w_2+(1-\theta_1-\theta_2)w_3\geq 0
	\end{align*}
	with \begin{align*}
		\theta_1w_1\leq \theta_1,\;\; \theta_2w_2\leq \theta_2,\;\;\text{and }\;\; (1-\theta_1-\theta_2)w_3\leq (1-\theta_1-\theta_2)
	\end{align*}
	Therefore, we have 
	\begin{align*}
		\theta_1w_1+\theta_2w_2+(1-\theta_1-\theta_2)w_3\leq \theta_1+\theta_2+(1-\theta_1-\theta_2)=1
	\end{align*}
	i.e. 
	\begin{align*}
		0\leq \theta_1w_1+\theta_2w_2+(1-\theta_1-\theta_2)w_3\leq 1
	\end{align*}
	for all $w_1,\;w_2,\;w_3\in W$ and $\theta_1,\;\theta_2\in[0,1].$ This results in, $\mathbf{W}$ is convex and consequently, the requirement (ii) is satisfied. In addition, the control of system \eqref{op_syseqn}'s right hand side is bilinear and continuous. Consequently, it can be expressed as follows:
	\begin{align*}
		\mathbf{f}(t,\mathbf{X},\mathbf{W})=\mathbf{A}(t,\mathbf{X})+\mathbf{B}(t,\mathbf{X})\mathbf{W}
	\end{align*}
	where, 
	\begin{align*}
		\mathbf{A}(t,\mathbf{X})=\begin{pmatrix}
			\Lambda-(\beta_1E+\beta_2I)S-(\mu+\phi)S\\
			\phi S-(1-\varepsilon)(\beta_1E+\beta_2I)V-\mu V\\
			(\beta_1E+\beta_2I)S-(\alpha+\mu)E\\
			\alpha E+(1-\varepsilon)(\beta_1E+\beta_2I)-(\mu+\delta+\gamma+\gamma_1)I\\
			\gamma I-\mu R\\
			\gamma_1 I-\mu T
		\end{pmatrix}
	\end{align*}
	and $\mathbf{B}(t,\mathbf{X})$ is similar as $\mathbf{f_W}(t,\mathbf{X},\mathbf{W})$ as described above, they are vector-valued functions of $\mathbf{X}$ and the boundedness of solution provides
	\begin{align*}
		\mathbf{f}(t,\mathbf{X},\mathbf{W})\leq C_1(1+|\mathbf{X}|+|\mathbf{W}|),
	\end{align*}
	where, $C_1$ relies on the coefficients of the system. Hence the condition (iii) is satisfied.\\
	For the verification of the convexity of the integrand $\mathcal{L}$ of our objective functional, $J$ we can present that, 
	\begin{align*}
		(1-\epsilon)J(t,\mathbf{X},\mathbf{W})+\epsilon J(t,\mathbf{X},\mathbf{W})\leq J(t,\mathbf{X},(1-\epsilon)\mathbf{W},\epsilon V)
	\end{align*}
	for, $0<\epsilon<1$ and $J(t,\mathbf{X},\mathbf{W})=a_1E+a_2I+(a_3w_1^2+a_4w_2^2+a_5w_3^2)$\\
	Now,
	\begin{align*}
		&(1-\epsilon)J(t,\mathbf{X},\mathbf{W})+\epsilon J(t,\mathbf{X},\mathbf{V})-J(t,\mathbf{X},(1-\epsilon)\mathbf{W}+\epsilon\mathbf{V})\\
		=\;&(1-\epsilon)[a_1E+a_2I+a_3w_1^2+a_4w_2^2+a_5w_3^2]+\epsilon[a_1E+a_2I+a_3w_1^2+a_4w_2^2+a_5w_3^2]-\\
		&[a_1E+a_2I+a_3((1-\epsilon)w_1+\epsilon v_1)^2+a_4((1-\epsilon)w_2+\epsilon v_2)^2+a_5((1-\epsilon)w_3+\epsilon v_3)^2]\\
		=\;&a_3\epsilon(1-\epsilon)(w_1-v_1)^2+a_4\epsilon(1-\epsilon)(w_2-v_2)^2+a_5\epsilon(1-\epsilon)(w_3-v_3)^2\leq 0
	\end{align*}
	since, $a_3,\;a_4,\;a_5>0$,\;$J(t,\mathbf{X},\mathbf{W})$ is convex in the set $\mathbf{W}$. As a consequence, we need to show that $J(t,\mathbf{X},\mathbf{W})\leq -C_1+C_1|\mathbf{W}|^{\beta},$ where, $C_1>0$ and $\beta >1.$ For our case 
	\begin{align*}
		J(t,\mathbf{X},\mathbf{W})=a_1E+a_2I+(a_3w_1^2+a_4w_2^2+a_5w_3)^2\geq -C_1+C_1|\mathbf{W}|^2.
	\end{align*}
	where, $C_2$ depends on the upper bound on (infected) CD$4^{+}$T cells, and $C_1>0$ since $a_3,\;a_4,\;a_5>0$ and $\beta=2.$ As a result, we can say that the problem has an ideal control pair \cite{Optimal Control-1, Optimal Control-4}.
\end{proof}

\subsection{The Optimality System and the Essential Prerequisites for the Optimal Control}
The necessary conditions for an optimal control problem is provided by the Poncryagin's Maximum Principle \cite{Optimal Control-6}. By applying this technique, the issue of maximizing the needed objective function $J$ subject to the state system \eqref{op_syseqn} is transformed into the problem of maximizing the Hamiltonian $\mathcal{H}$, pointwise with regard to control parameters $w_1$, $w_2$, and $w_3$. Consequently, for describing the optimal controls $w_1^*$, $w_2^*$, and $w_3^*$, it suffices to deduce the Hamiltonian $\mathcal{H}$ rather than deriving the objective function $J$ specified in \eqref{opti-obj-function}. We now derive the required requirements that must be met by an optimal control and related states using Pontryagin's Maximum principle. According to the definition, the Hamiltonian $(\mathcal{H})$ and Lagrangian $(\mathcal{L})$ relate to the optimal control problem is defined as,
\begin{align*}
	\mathcal{L}= a_1 E+a_2 I+(a_3 w_1^2+a_4 w_2^2+a_5 w_3^2)
\end{align*}
and 
\begin{align}\label{op_hamiltonian}
	\mathcal{H}=&a_1E+a_2I+(a_3 w_1^2+a_4 w_2^2+a_5 w_3^2) \nonumber\\
	& + \Psi_S\left[\Lambda-(1-w_1(t))(\beta_1 E+\beta_2 I)S-\mu S-\phi S\right] \nonumber\\
	& + \Psi_V\left[\phi S-(1-\varepsilon)(\beta_1 E+\beta_2 I)V-\mu V\right] \nonumber\\
	& + \Psi_E\left[(1-w_1(t)(\beta_1 E+\beta_2 I)S-\alpha E-\mu E)\right] \nonumber\\
	& + \Psi_I\left[\alpha E+(1-\varepsilon)(\beta_1 E+\beta_2 I)V-\mu I-\delta I-(1+w_2(t))\gamma_1 I-(1+w_3(t))\gamma I\right] \nonumber\\
	& + \Psi_R\left[\gamma(1+w_3(t))I-\mu R\right] \nonumber\\
	&+\Psi_T\left[\gamma_1(1+w_2(t))I-\mu T\right]
\end{align}
where, $\Psi_S$,$\Psi_V$,$\Psi_E$,$\Psi_I$,$\Psi_R$ and $\Psi_T$ are adjoint variables.\\
or by simplicity this can be written as 
\begin{align*}
	\mathcal{H}=&a_1E+a_2I+(a_3 w_1^2+a_4 w_2^2+a_5 w_3^2)+\sum_{i=1}^{6}\Psi_i(t)F_i.
\end{align*}
Hence, in the model \eqref{op_syseqn}, $F_i$ is represented as the right hand side of the differential equation of the $i^{th}$ state variable.\\
We use Pontryagin's maximal principle to get the answer to the Influenza optimum control issue, \eqref{op_syseqn}. We now describe the following theorem using the equation \eqref{op_hamiltonian}.
\begin{Th}\label{op_th_main}
	An optimal control $\mathbf{W^*}=(w_1^*,\;w_2^*,\;w_3^*)$ exists there and related solutions $S^*,\;V^*,\;E^*,\;I^*,\;R^*$ and $T^*$ of the ideal system of controls \eqref{op_syseqn}, that minimizes the objective functional $J(w_1,w_2,w_3)$ over $\mathcal{D}$. Then there presents adjoint variables
	$\Psi_S(t),\;\Psi_V(t),\;\Psi_E(t),\;\Psi_I(t),\;\Psi_R(t)$\; and \;$\Psi_T(t)$ satisfying the equations:
	\begin{align}\label{op_cond_eqns}
		\frac{d\Psi_S}{dt} =& \Psi_S[(1-w_1(t))(\beta_1 E+\beta_2 I)+\mu+\phi]-\Psi_V \phi +\Psi_E[-(1-w_1(t))(\beta_1 E+\beta_2 I)] \nonumber\\
		\frac{d\Psi_V}{dt} =& \Psi_V[(1-\varepsilon)(\beta_1 E+\beta_2 I)+\mu]+\Psi_I[-(1-\varepsilon)(\beta_1 E+\beta_2 I)] \nonumber\\
		\frac{d\Psi_E}{dt} =&-a_1+\Psi_S[(1-w_1(t))\beta_1 S]+\Psi_V[(1-\varepsilon)\beta_1 V]+\Psi_E[-(1-w_1(t))\beta_1 S+\alpha+\mu]+ \nonumber\\
		&\Psi_I[-(1-\varepsilon)\beta_1 V] \nonumber\\
		\frac{d\Psi_I}{dt}=& -a_2+\Psi_S[(1-w_1(t))\beta_2S]+\Psi_V[(1-\varepsilon)\beta_2 V]+\Psi_E[-(1-w_1(t))\beta_2 S]+\Psi_I[-(1-\varepsilon)\beta_2V+\nonumber\\
		&\mu+\delta+(1+w_2(t))\gamma_1+(1+w_3(t))\gamma]+\Psi_R[-\gamma(1+w_3(t))]+\Psi_T[-\gamma_1(1+w_2(t))] \nonumber\\
		\frac{d\Psi_R}{dt}=&\Psi_R\mu \;\;,\;\;\frac{d\Psi_T}{dt}=\Psi_T \mu
	\end{align}
   under situations of transversality $\displaystyle \Psi_S(t_f)=\Psi_V(t_f)=\Psi_E(t_f)=\Psi_I(t_f)=\Psi_R(t_f)=\Psi_T(t_f)=0.$\\
	and the following provide the optimality conditions: 
	\begin{align}\label{op_cond_var}
		w_1^*(t)= &\max\left\{0,\;\min\left\{1,\;\frac{(\Psi_E-\Psi_S)(\beta_1 E^*+\beta_2 I^*)S^*}{2a_3}\right\}\right\} \nonumber\\
		w_2^*(t)= &\max\left\{0,\;\min\left\{1,\;\frac{(\Psi_I-\Psi_T)\gamma_1I^*}{2a_4}\right\}\right\} \nonumber\\
		w_3^*(t)= &\max\left\{0,\;\min\left\{1,\;\frac{(\Psi_I-\Psi_R)\gamma I^*}{2a_5}\right\}\right\}
	\end{align}
\end{Th}

\begin{proof}
	The boundedness of the solutions to the system \eqref{op_syseqn} and the convexity of the integrand of $J$ with respect to $w_1,\; w_2$ and $w_3$ are readily demonstrable. Additionally, we are able to demonstrate the system's Lipschitz property with reference to the state variables. Utilizing the corollary and these qualities in \cite{Optimal Control-6}, the optimal control's existence is established. The transversality conditions and adjoint equations can be found by applying the Pontryagin's Maximum principle in such a way that,
	\begin{align*}
		\frac{d\Psi_S}{dt}=&-\frac{\partial \mathcal{H}}{\partial S},\;\;\; \Psi_S(t_f)=0\\
		\frac{d\Psi_V}{dt}=&-\frac{\partial \mathcal{H}}{\partial V},\;\;\; \Psi_V(t_f)=0\\
		\frac{d\Psi_E}{dt}=&-\frac{\partial \mathcal{H}}{\partial E},\;\;\; \Psi_E(t_f)=0\\
		\frac{d\Psi_I}{dt}=&-\frac{\partial \mathcal{H}}{\partial I},\;\;\; \Psi_I(t_f)=0\\
		\frac{d\Psi_R}{dt}=&-\frac{\partial \mathcal{H}}{\partial R},\;\;\; \Psi_R(t_f)=0\\
		\frac{d\Psi_T}{dt}=&-\frac{\partial \mathcal{H}}{\partial T},\;\;\; \Psi_T(t_f)=0
	\end{align*}
	where we obtain the costate functions by using the Hamiltonian associated with $S,V,E,I,R$, and $T$, we obtain,
	\begin{align*}
		\Psi_S'(t)=-\frac{\partial \mathcal{H}}{\partial S}=&\Psi_S[(1-w_1(t))(\beta_1 E+\beta_2 I)+\mu+\phi]-\Psi_V\phi+\Psi_E[-(1-w_1(t))(\beta_1 E+\beta_2 I)]\\
		\Psi_V'(t)=-\frac{\partial \mathcal{H}}{\partial V}=&\Psi_V[(1-\varepsilon)(\beta_1 E+\beta_2 I)+\mu]+\Psi_I[-(1-\varepsilon)(\beta_1 E+\beta_2 I)]\\
		\Psi_E'(t)=-\frac{\partial \mathcal{H}}{\partial E}=&-a_1+\Psi_S[(1-w_1(t))\beta_1 S]+\Psi_V[(1-\varepsilon)\beta_1 V]+\\
		&\Psi_E[-(1-w_1(t))\beta_1 S+\alpha+\mu]+\Psi_I[-(1-\varepsilon)\beta_1 V]\\
		\Psi_I'(t)=-\frac{\partial \mathcal{H}}{\partial I}=&-a_2+\Psi_S[(1-w_1(t))\beta_2S]+\Psi_V[(1-\varepsilon)\beta_2V]+\Psi_E[-(1-w_1(t))\beta_2S]+\\
		&\Psi_I[-(1-\varepsilon)\beta_2V+\mu+\delta+(1+w_2(t))\gamma_1+(1+w_3(t))\gamma]+\Psi_R[-\gamma(1+w_3(t))]+\\
		&\Psi_T[-\gamma_1(1+w_2(t))]\\
		\Psi_R'(t)=-\frac{\partial \mathcal{H}}{\partial R}=&\Psi_R \mu\\
		\Psi_T'(t)=-\frac{\partial \mathcal{H}}{\partial T}=&\Psi_T \mu
	\end{align*}
	Assessing the adjoint system at the optimal control derivatives.\\
	Because at final time, $S(t),\;V(t),\;E(t),\;I(t),\;R(t)$, and $T(t)$ do not have stable values $T_f$, the corresponding adjoints values $\Psi_S(t),\;\Psi_V(t),\;\Psi_E(t),\;\Psi_I(t),\;\Psi_R(t)$ and $\Psi_T(t)$ at the final time are zero. The optimal control $w_1^*,\;w_2^*$ and $w_3^*$ on the interior of the control set can be solved by using the optimality conditions,
	\begin{align*}
		\frac{\partial \mathcal{H}}{\partial w_1}\big|_{w_1=w_1^*}=0,\;\frac{\partial \mathcal{H}}{\partial w_2}\big|_{w_2=w_2^*}=0,\;\text{and}\; \frac{\partial \mathcal{H}}{\partial w_3}\big|_{w_3=w_3^*}=0.
	\end{align*}
	The optimality requirements on the control set $\mathcal{D}$ are defined as follows:
	\begin{equation}\label{op_Deqn}
		\mathcal{D}=\left\{(w_1(t),w_2(t),w_3(t)) | 0\leq (w_1(t),w_2(t),w_3(t))\leq 1\right\}
	\end{equation}
	as 
	\begin{align}\label{op_sys_www}
		\frac{\partial \mathcal{H}}{\partial w_1}=&2a_3w_1+\Psi_S(\beta_1 E^*+\beta_2 I^*)S^*-\Psi_E(\beta_1 E^*+\beta_2 I^*)S^*=0 \nonumber\\
		\frac{\partial \mathcal{H}}{\partial w_2}=&2a_4w_2-\Psi_I\gamma_1I^*+\Psi_T\gamma_1I^*\\
		\frac{\partial \mathcal{H}}{\partial w_3}=&2a_5w_3-\Psi_I\gamma I^*+\Psi_R\gamma I^* \nonumber
	\end{align}
	Solving \eqref{op_sys_www} we get,
	\begin{align*}
		w_1^*=&\frac{(\Psi_E-\Psi_S)(\beta_1 E^*+\beta_2 I^*)S^*}{2a_3}\\
		w_2^*=&\frac{(\Psi_I-\Psi_T)\gamma_1 I^*}{2a_4}\\
		w_3^*=&\frac{(\Psi_I-\Psi_R)\gamma I^*}{2a_5}
	\end{align*}
	Using three control standard intervals, we obtain:
	\begin{equation*}  
		w_1^*(t) = 
		\begin{cases}
			\vspace{0.2cm}
			\displaystyle 0\;, &\quad\text{if}\;\;\;\;\frac{(\Psi_E-\Psi_S)(\beta_1 E^*+\beta_2 I^*)S^*}{2a_3}\leq 0\\
			\vspace{0.2cm}
			\displaystyle \frac{(\Psi_E-\Psi_S)(\beta_1 E^*+\beta_2 I^*)S^*}{2a_3}\;, &\quad\text{if}\;\;\;\; 0 <\frac{(\Psi_E-\Psi_S)(\beta_1 E^*+\beta_2 I^*)S^*}{2a_3} <1 \\
			\displaystyle 1\;, &\quad\text{if}\;\;\;\; \frac{(\Psi_E-\Psi_S)(\beta_1 E^*+\beta_2 I^*)S^*}{2a_3}\geq 1\\
		\end{cases}
	\end{equation*}
	can be rewriteable in a condensed form:
	\begin{align*}
		w_1^*(t) = \max\left\{0,\;\min\left\{1,\;\frac{(\Psi_E-\Psi_S)(\beta_1 E^*+\beta_2 I^*)S^*}{2a_3}\right\}\right\}
	\end{align*}
	Similarly, for $w_2^*$ and $w_3^*$ we obtain,
	\begin{align*}
		w_2^*=&\max\left\{0,\;\min\left\{1,\;\frac{(\Psi_I-\Psi_T)\gamma_1I^*}{2a_4}\right\}\right\}\\
		w_3^*=&\max\left\{0,\;\min\left\{1,\;\frac{(\Psi_I-\Psi_R)\gamma I^*}{2a_5}\right\}\right\}
	\end{align*}
	Additionally, the Hamiltonian's $\mathcal{H}$ second derivative with regard to the parameters $w_1(t)$, $w_2(t)$, and $w_3(t)$ is positive, indicating a maximum at $\mathbf{W^*}=(w_1,w_2,w_3).$ That reflects,
	\begin{align*}
		\frac{\partial^2\mathcal{H}}{\partial w_i^2}=2a_j,\;\; i=1,2,3\;\;\text{and}\;\;j=3,4,5\;\;\text{since}\;\;a_j\geq 0 
	\end{align*}
	Based on our findings, the optimality system is made up of the optimality condition \eqref{op_cond_var}, the corresponding adjoint system \eqref{op_cond_eqns}, and the state system \eqref{op_syseqn} with the initial conditions. As a result, we have the optimality system at $\mathbf{W^*}=(w_1^*(t),w_2^*(t),w_3^*(t)).$\\
	We achieve the desired compact form of the optimal control efforts, which are given by equations, taking into consideration the bounds on the controls (\ref{op_cond_var}). Therefore, the equations \eqref{op_syseqn}, \eqref{op_cond_eqns} with the optimality conditions (\ref{op_cond_var}) and the non-negative primary conditions and the transversality conditions, $\Psi_S(t_f)=\Psi_V(t_f)=\Psi_E(t_f)=\Psi_I(t_f)=\Psi_R(t_f)=\Psi_T(t_f)=0$ form the optimality system \cite{Optimal Control-4, Optimal Control-6}.\\
	We obtain unique solutions of the optimality system for the interval $[0,t_f]$ because the state variables and the adjoint functions are bounded and the state and adjoint system has Lipschitz structure with respect to the relevant variables. Thus, for $t \in [0,t_f]$, we may state that the boundedness and uniqueness of the solutions to the optimality system stated previously exist.
\end{proof}

\subsection{Numerical Simulation of Optimal Control}
\begin{figure}[H]
	\centering  
	\subfloat[]{\includegraphics[width=2.7 in]{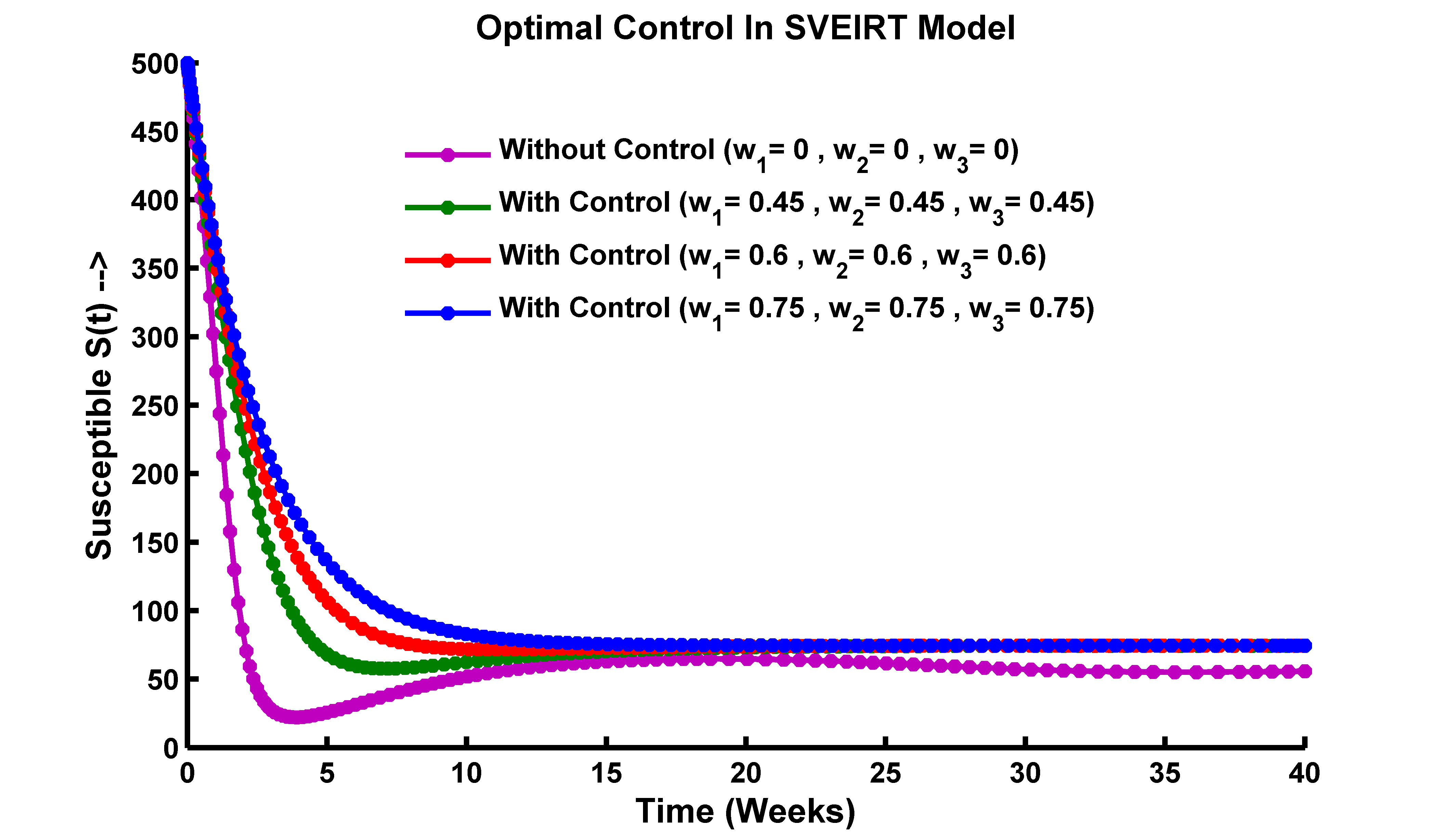}}
	\subfloat[]{\includegraphics[width=2.7 in]{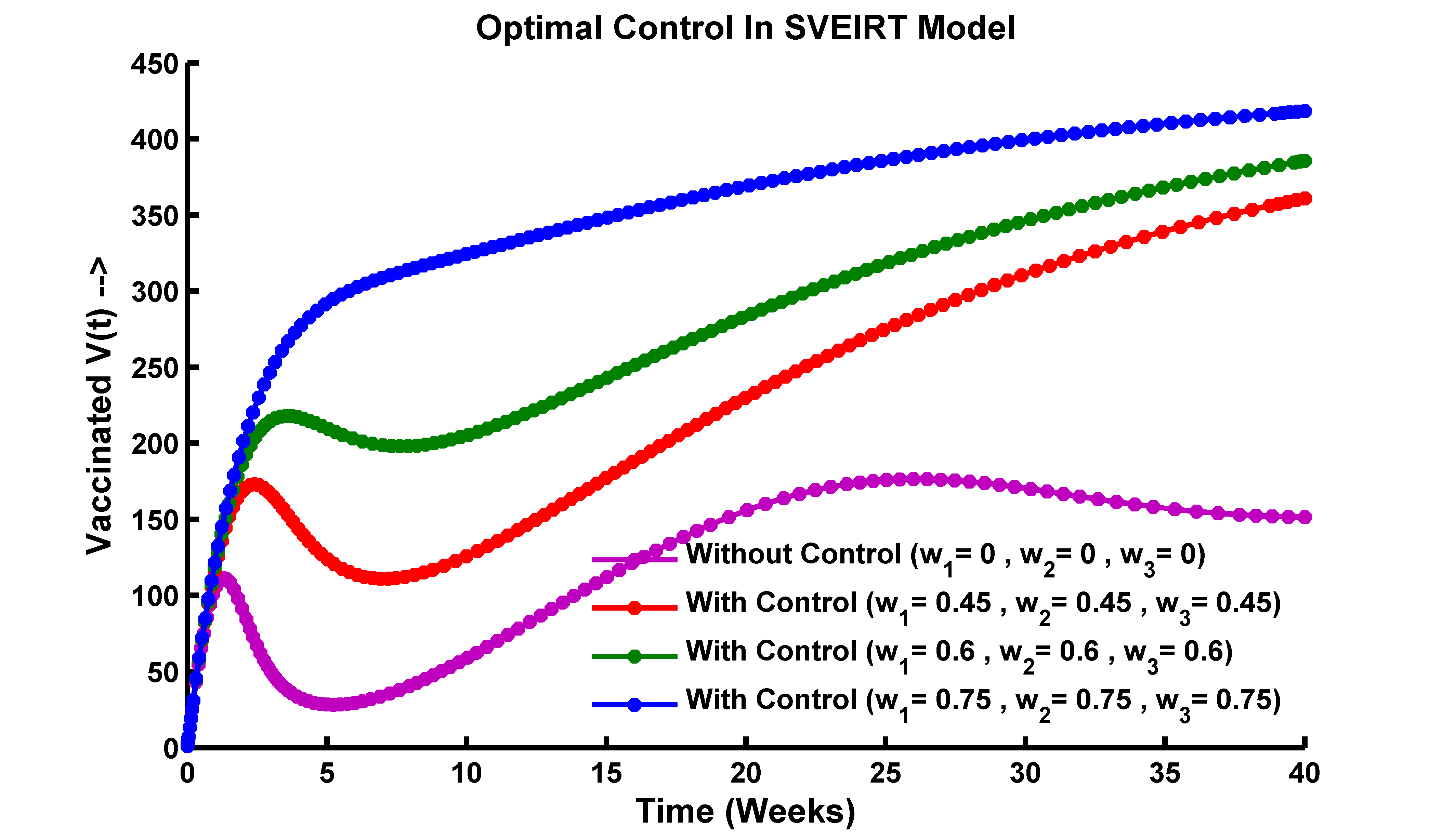}}\\
	\subfloat[]{\includegraphics[width=2.7 in]{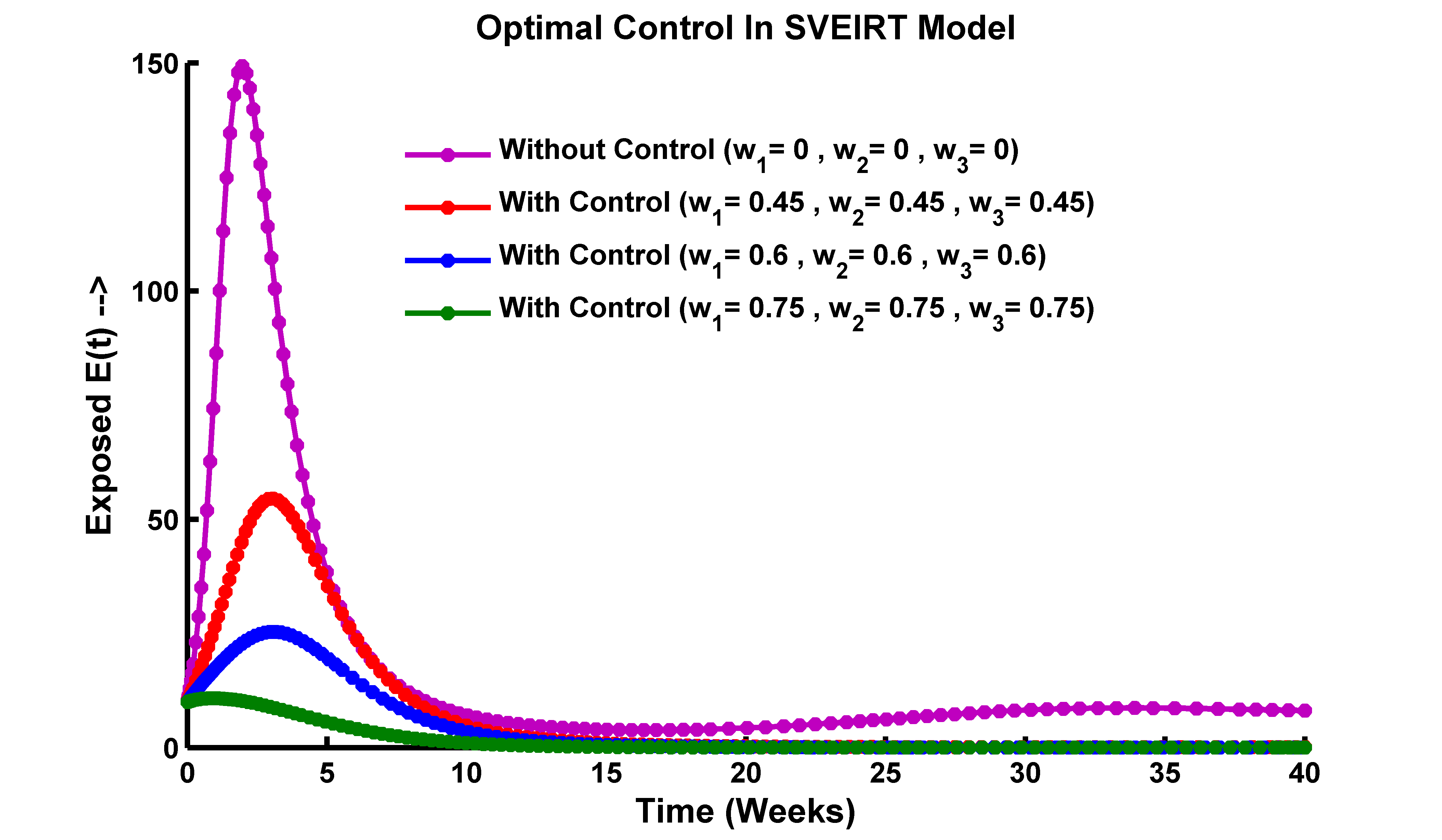}}
	\subfloat[]{\includegraphics[width=2.7 in]{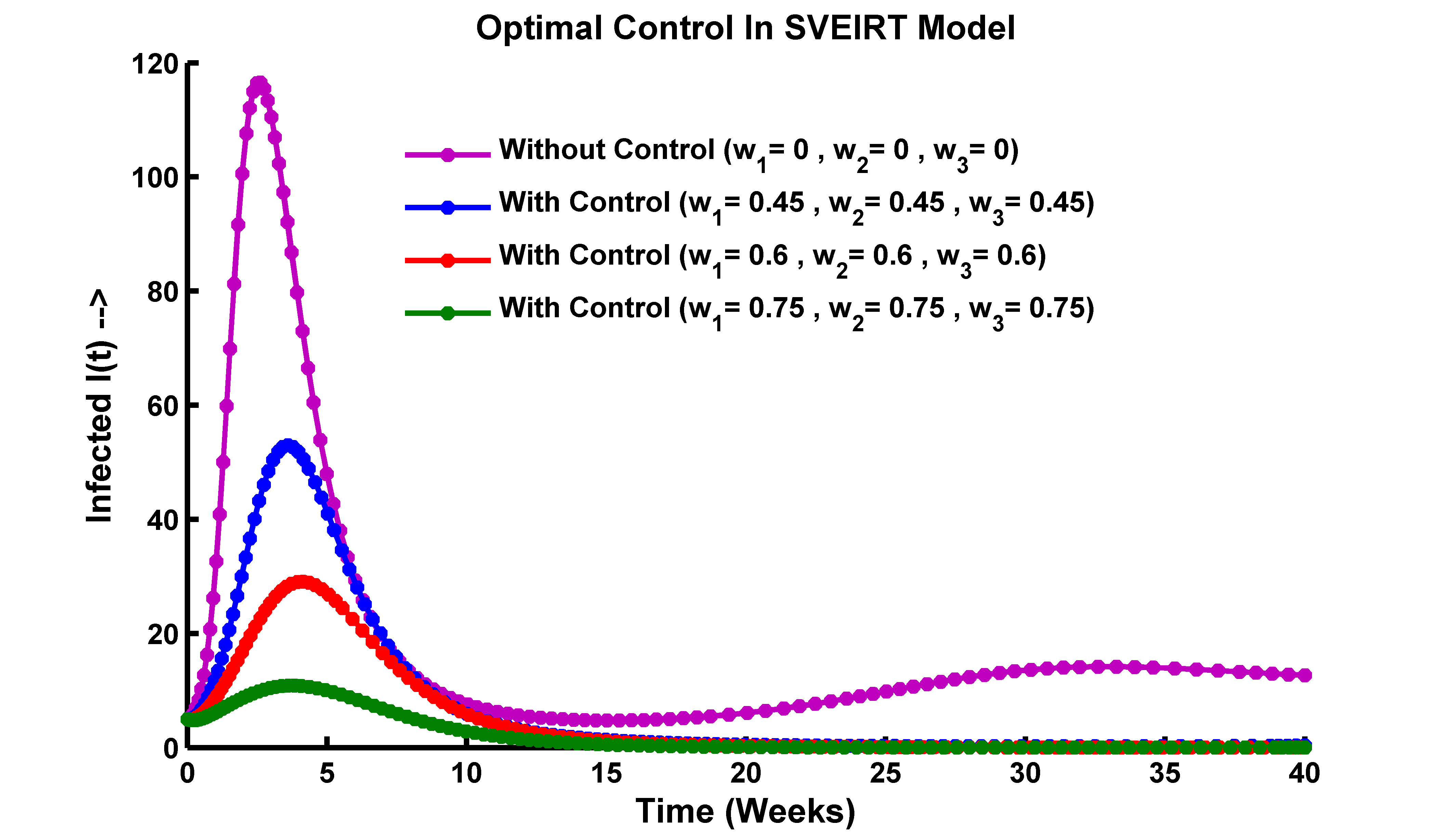}}\\
	\subfloat[]{\includegraphics[width=2.7 in]{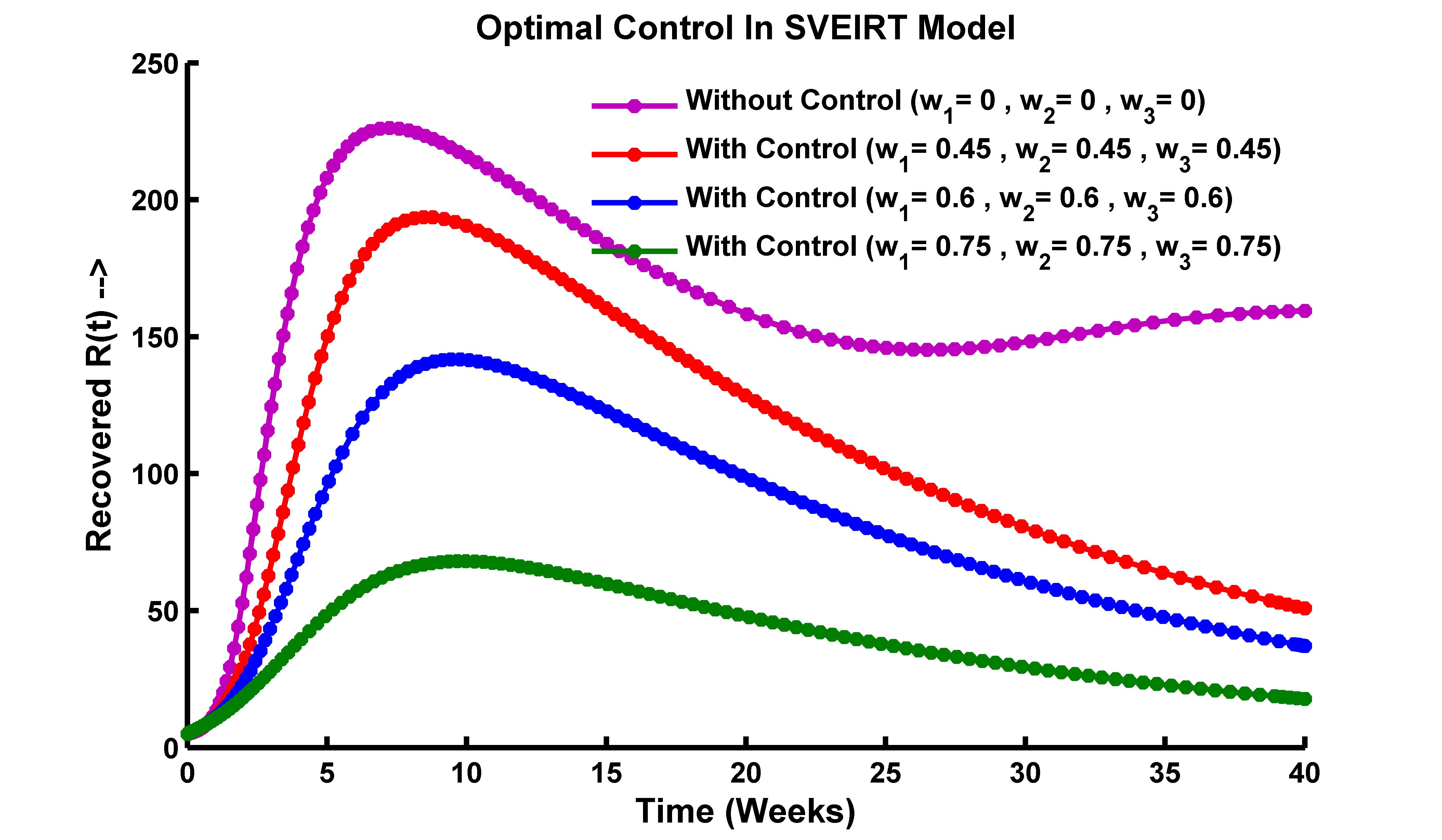}}
	\subfloat[]{\includegraphics[width=2.7 in]{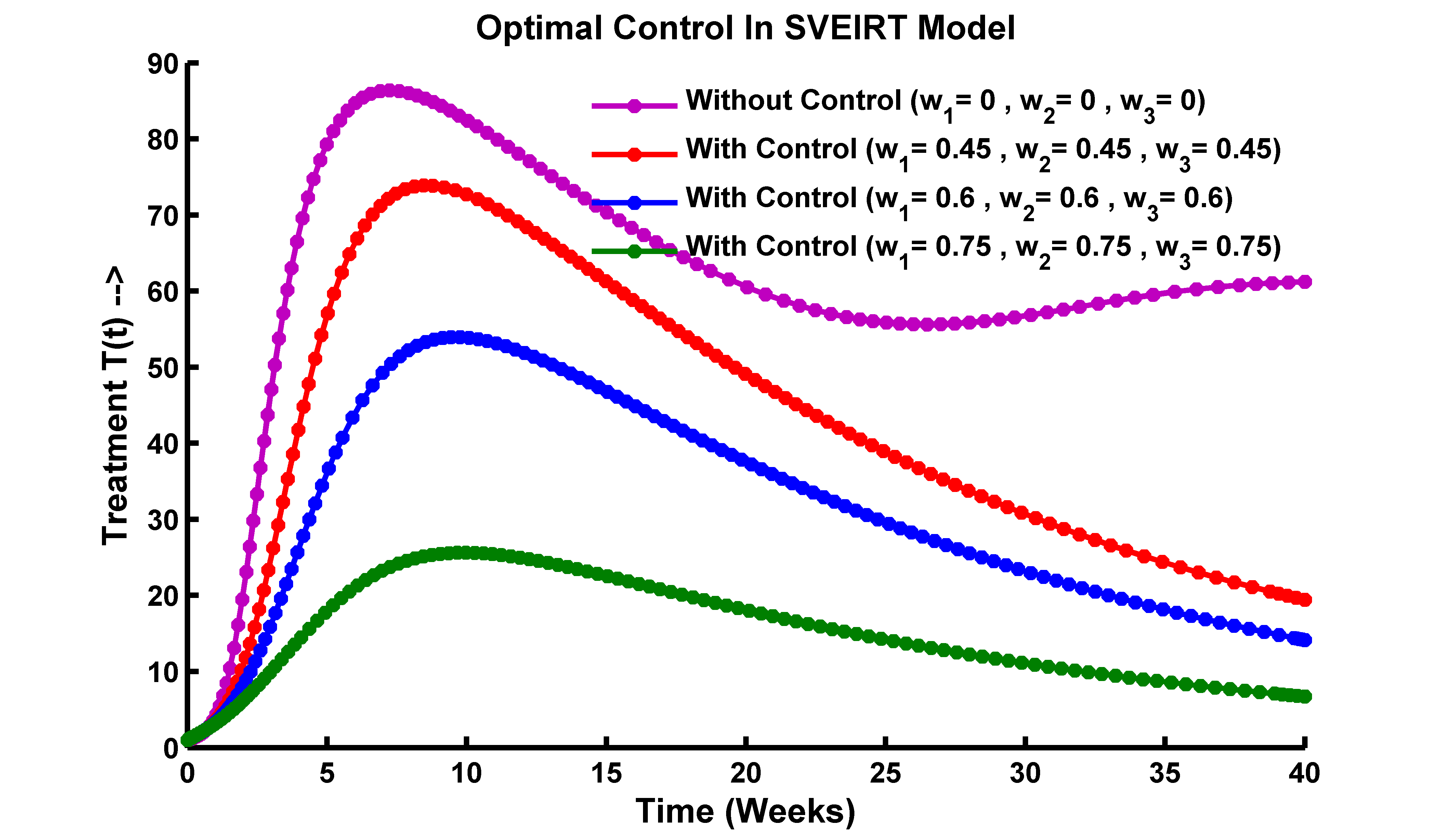}}
	\caption{Simulation of the model (a) $S(t)$ compartment, (b) $V(t)$ compartment (c) $E(t)$ compartment (d) $I(t)$ compartment (e) $R(t)$ compartment and (f) $T(t)$ compartment showing with control and without control with respect to control parameters $w_1,w_2,w_3$, where $w_1=w_2=w_3\in[0,0.45,0.6,0.75]$.}
	\label{Optimal-control vs w}
\end{figure}

\begin{figure}[H]
	\centering  
	\subfloat[]{\includegraphics[width=2.7 in]{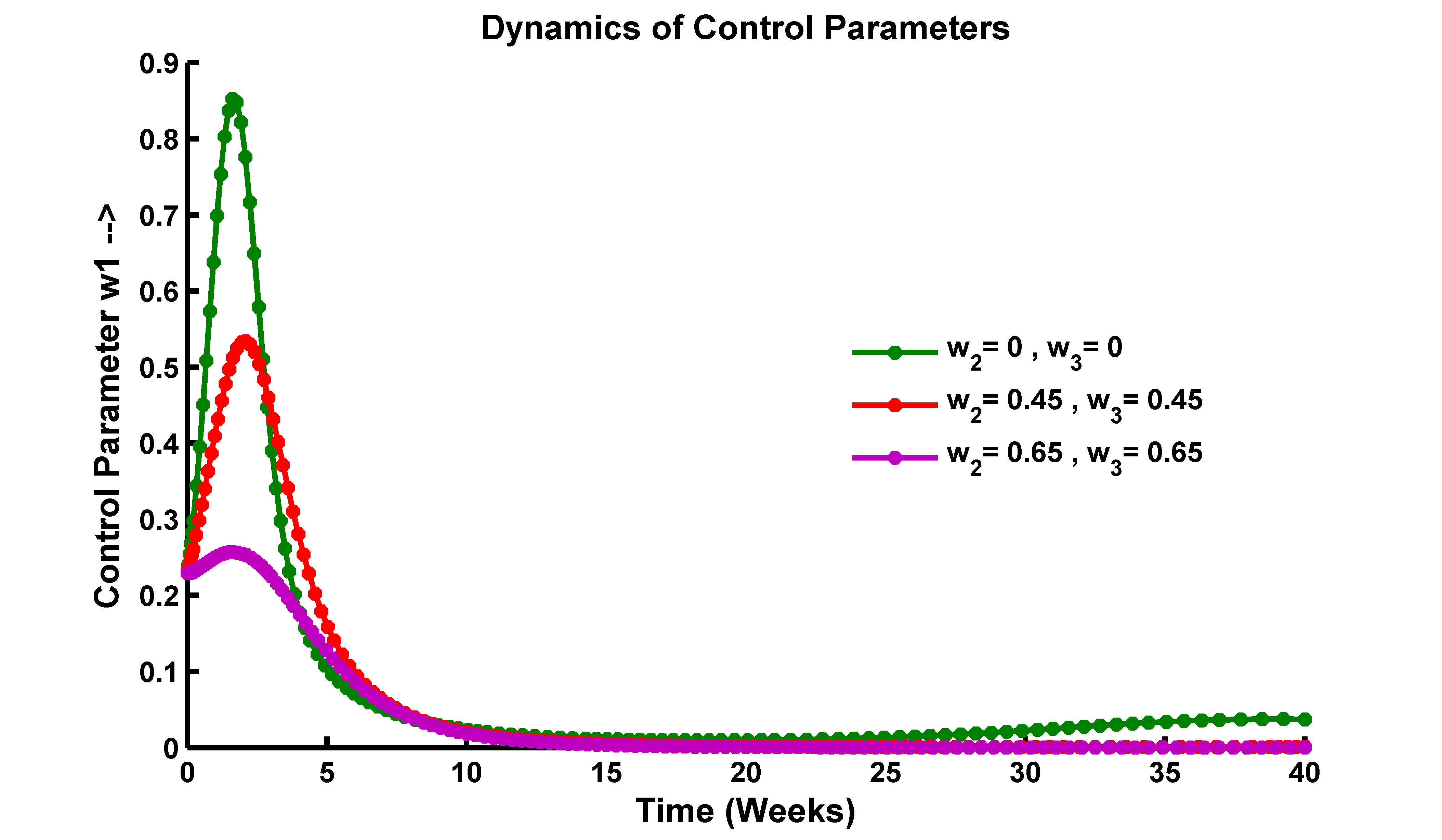}}
	\subfloat[]{\includegraphics[width=2.7 in]{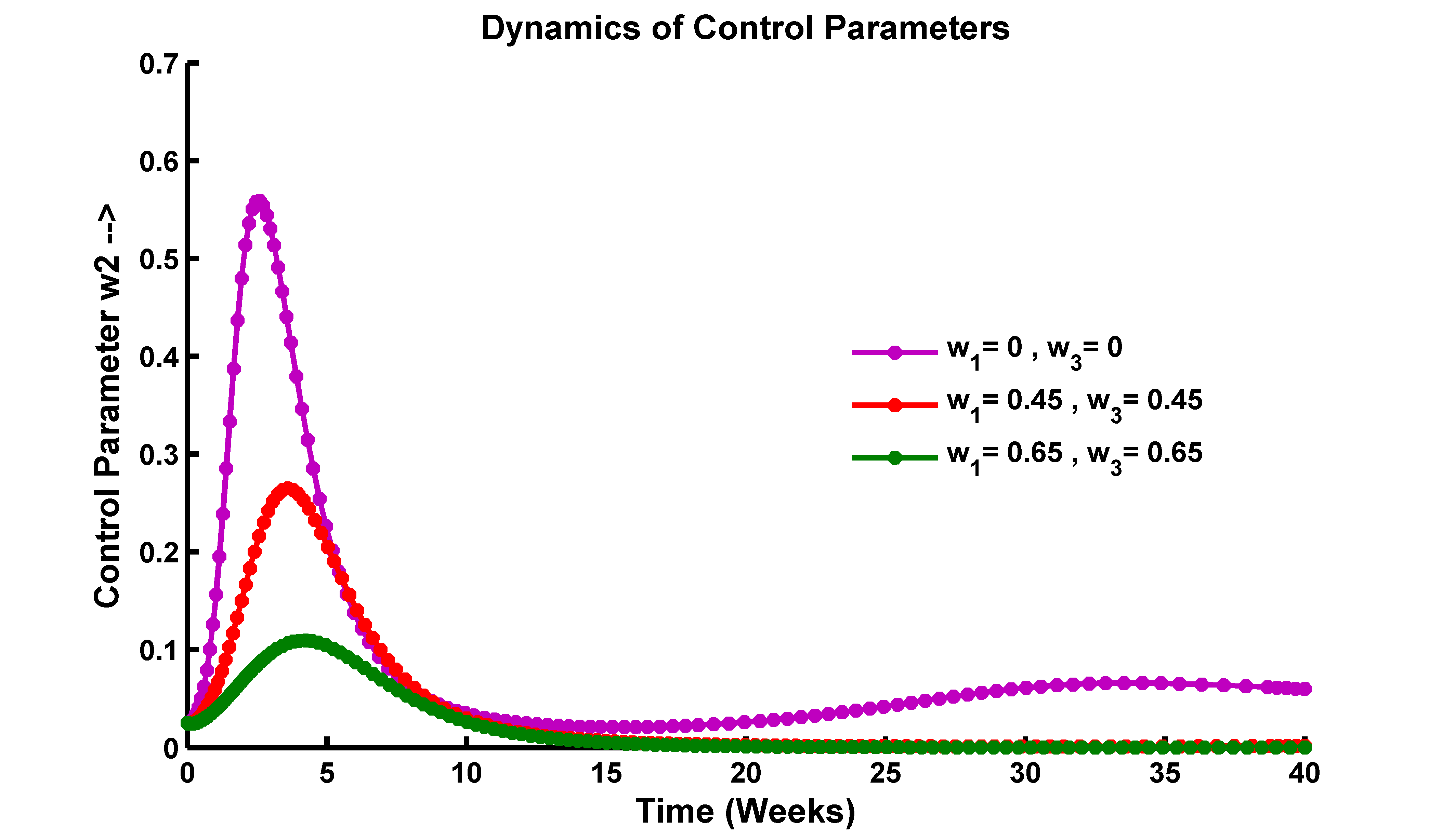}}\\
	\subfloat[]{\includegraphics[width=2.7 in]{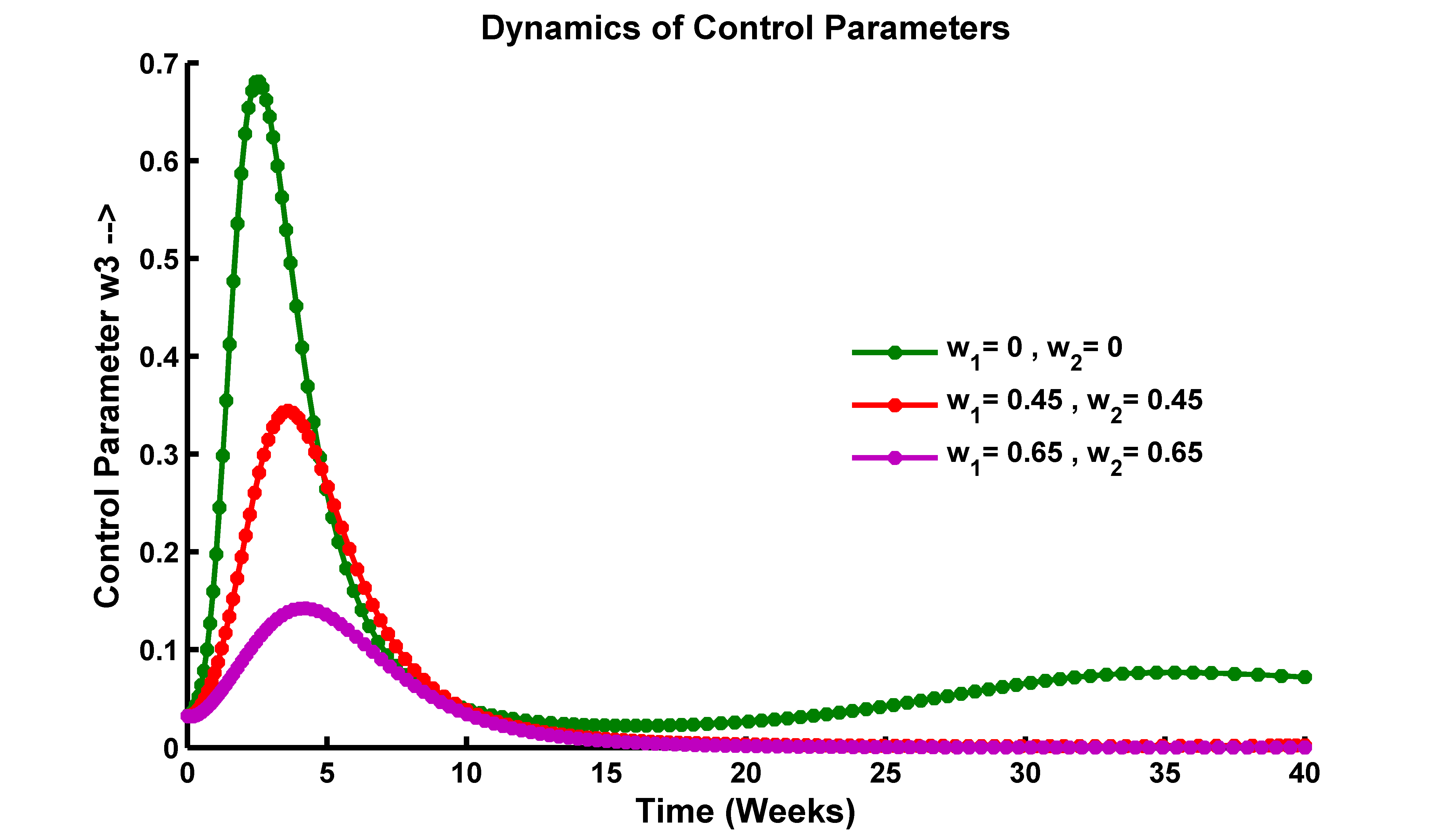}}
	\caption{Phase portrait of different control efforts (a) control strategies $w_1$ (b) control strategies $w_2$ and (c) control strategies $w_3$.}
	\label{Optimal-control parameters}
\end{figure}

This section conducts a numerical simulation of the suggested influenza model with optimum control \eqref{op_syseqn} in order to study how control techniques affect illness transmission and control costs. Table \ref{tableparameter} lists the values of the influenza model's parameters. The state variables initial conditions are $S(0)=500,\;V(0)=1,\;E(0)=1,\;I(0)=0,\;R(0)=0$ and $T(0)=0$ and the weight constants $a_1=20$, $a_2=20$, $a_3=45$, $a_4=25$ and $a_5=50$ taking consideration of four different cases bases on three controls, $w_1$ (social separation and mask-wearing to cut down on meaningful interactions), $w_2$ (treatment of hospitalized patients with the virus) and $w_3$ (increased rate of treatment), To examine each control strategy's efficacy, we simulate the control model \eqref{op_syseqn} that comprises of
\begin{enumerate}
	\item[(i)] adoption of an approach with no control based on $w_1=0,\;w_2=0$ and $w_3=0$
	\item[(ii)] implementation of three controls simultaneously, and in this scenario $w_1=0.45,\;w_2=0.45$ and $w_3=0.45$
	\item[(iii)] implementation of three controls simultaneously, and in this scenario $w_1=0.6,\;w_2=0.6$ and $w_3=0.6$
	\item[(iv)] implementation of three controls simultaneously, and in this scenario $w_1=0.75,\;w_2=0.75$ and $w_3=0.75$
\end{enumerate}
Figure \ref{Optimal-control vs w}(a)-(f) shows the variations in susceptible, immunized, exposed, infected, recovered, and treated subpopulations' densities, respectively, with and without control, taking into account four distinct scenarios and the various control profiles are displayed in Figure \ref{Optimal-control parameters}. Figure \ref{Optimal-control vs w} reflects that in the first scenario to minimize meaningful interaction by keeping a physical distance that is $w_1=0,\;w_2=0$ and $w_3=0$. In the fourth scenario, all three controls were active at the same time ($w_1\neq0,\;w_2\neq0$ and $w_3\neq0$), the number of people who were exposed drops off quickly in comparison to the first, second, and third cases, whenever $w_1=w_2=w_3=0.45$, $w_1=w_2=w_3=0.6$ and in absence of control scenario. It entails keeping a safe distance between hospitalized, exposed, and sick people in order to reduce the risk of contracting from susceptible others. The sick subpopulation has significantly decreased, as has the exposed subpopulation for the first scenario ($w_1=0,\;w_2=0$ and $w_3=0$) and fourth strategy ($w_1=w_2=w_3=0.75$) and slight decrease is observed in the second, third strategy where ($w_1=w_2=w_3=0.45$) and $w_1=0=w_2=w_3=0.6$, respectively. Figure \ref{Optimal-control vs w} depicts the details observation.Only improved inpatient patient care for both the exposed and infected populations where, $w_1=w_2=w_3=0$ and $w_1=w_2=w_3=0.45$ has no appreciable effect on the lower infection rate. It implies that the community will continue to be infected.

Figure \ref{Optimal-control vs w}(a) and \ref{Optimal-control vs w}(b) reflects that, with control strategy has effective impact on susceptible class and vaccinated class. For the third and fourth case in control strategies where $w_1=w_2=w_3=0.6$ and $w_1=w_2=w_3=0.75$, the portion of susceptible and vaccinated individuals increases in the community than without control (i) and strategy (ii). This suggests that subpopulation strategies (iii) and (iv), for both exposed and afflicted individuals, are useful and beneficial in lowering infection rates in the community.

The biological effects of each of these four approaches for both hospitalized and recovered patients are displayed in Figure \ref{Optimal-control vs w}(e)-(f). Physical distancing control and treatment control are considered in the third technique ($w_1=w_2=w_3=0.6$) and fourth strategy where treatment control increases ($w_1=w_2=w_3=0.75$) are considered. Giving treatment to hospitalized patients alone will not be sufficient to minimize the illness in the community if physical distancing measures like mask wear and isolation are not strictly enforced. Therefore, the techniques (ii), (iii), and (iv) are more effective for hospitalized and recovered patients than the situation (i) without a control.

Figure \ref{Optimal-control vs w}(e), reflects how the recovered subgroup changed both with and without various control measures. Since the exposed, infected, and hospitalized populations decline more quickly for the third and fourth strategies (in which all three controls are implemented concurrently) than for the other strategies (i) and (ii), the number of recovered individuals likewise declines noticeably for the fourth strategy. Without any control and while using the second control technique, there is an increase in the subpopulation that has been recovered from and treated. That happens as a result of the community's high exposure and infection rates.

Moreover, Figure \ref{Optimal-control parameters} shows the control profiles of implementation of all three controls ($w_1,\;w_2,\;w_3$) simultaneously, when joint control implementation occurs, roughly takes approximately 80\%-100\% physical distance as an input ($w_1$) for 3-4 weeks and 80-90\% input of treatment rates $(w_2)$ and $(w_3)$ for 6-7 weeks prior to reaching their bottom limitations, which yields superior outcomes compared to the initial two techniques (i), (ii).

\begin{figure}[H]
	\centering  
	\subfloat[]{\includegraphics[width=2.7 in]{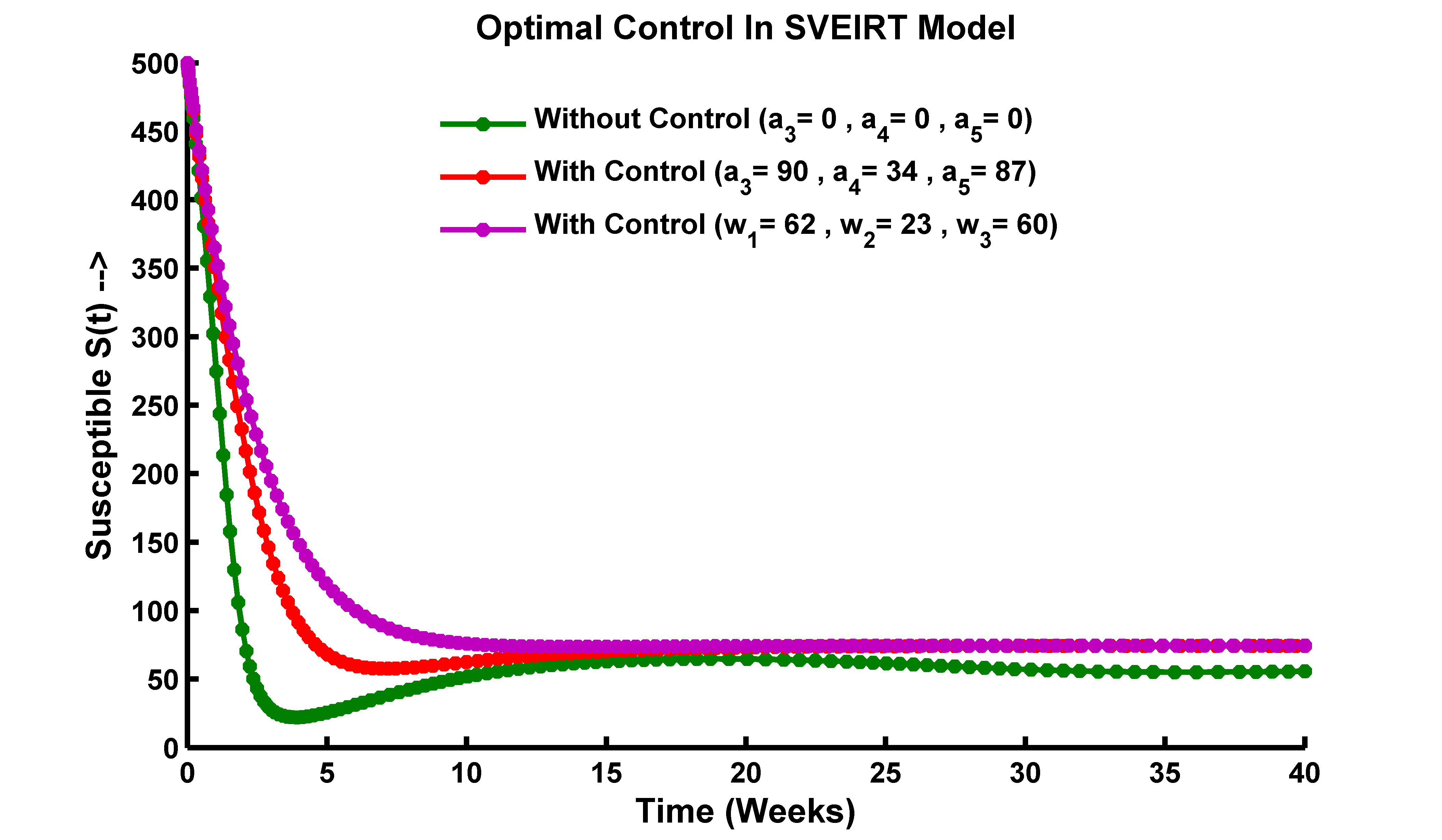}}
	\subfloat[]{\includegraphics[width=2.7 in]{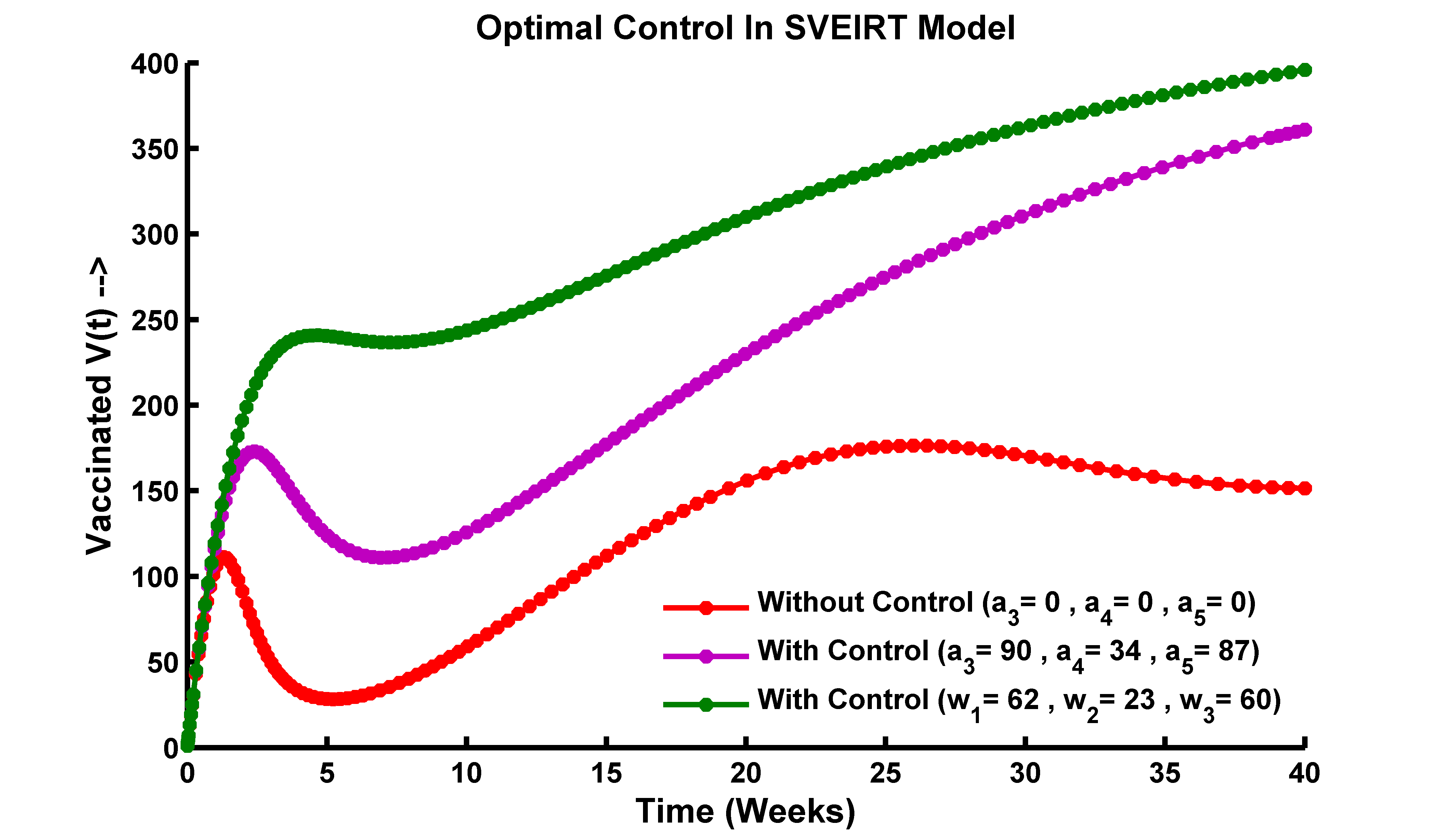}}\\
	\subfloat[]{\includegraphics[width=2.7 in]{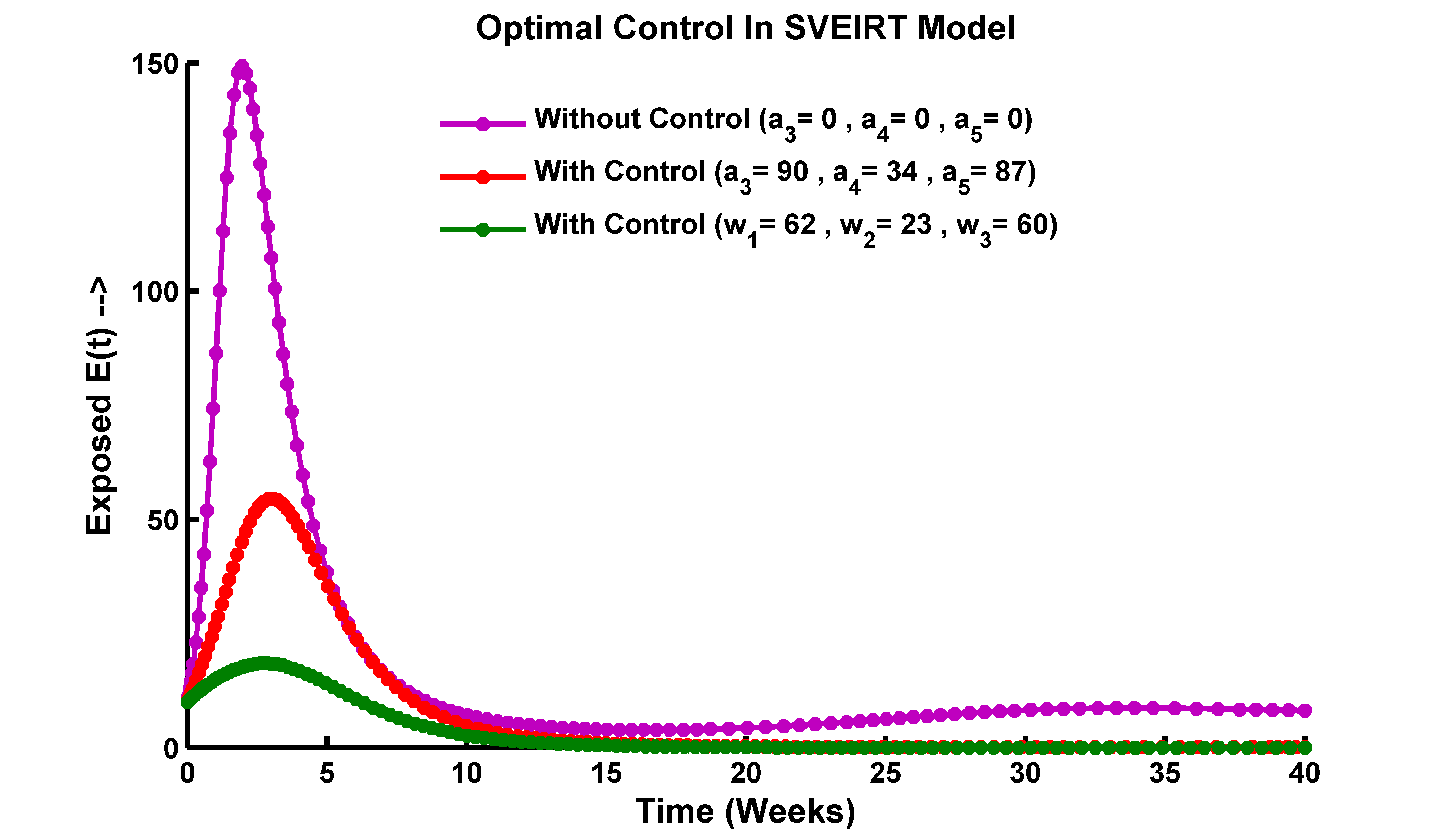}}
	\subfloat[]{\includegraphics[width=2.7 in]{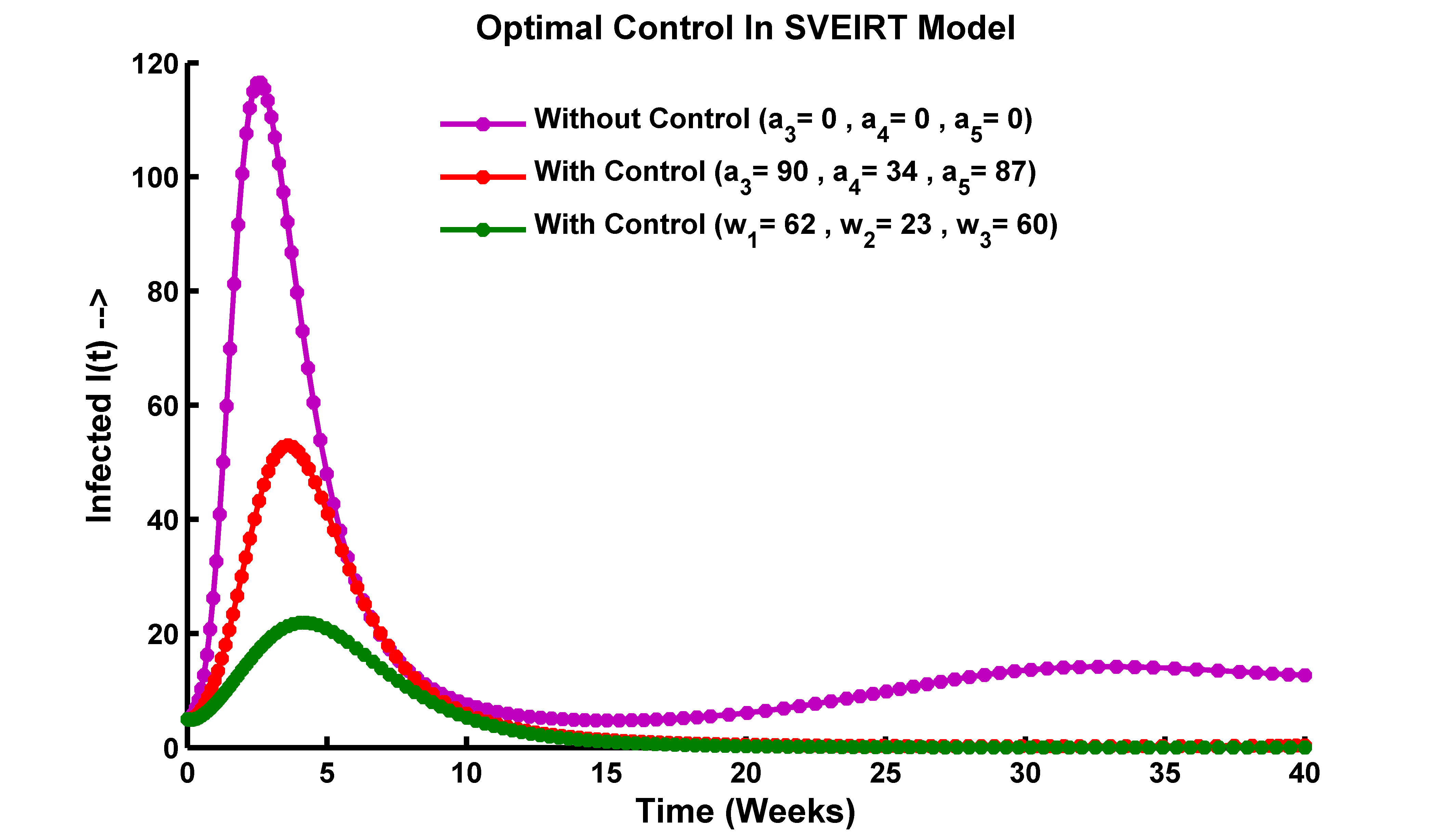}}\\
	\subfloat[]{\includegraphics[width=2.7 in]{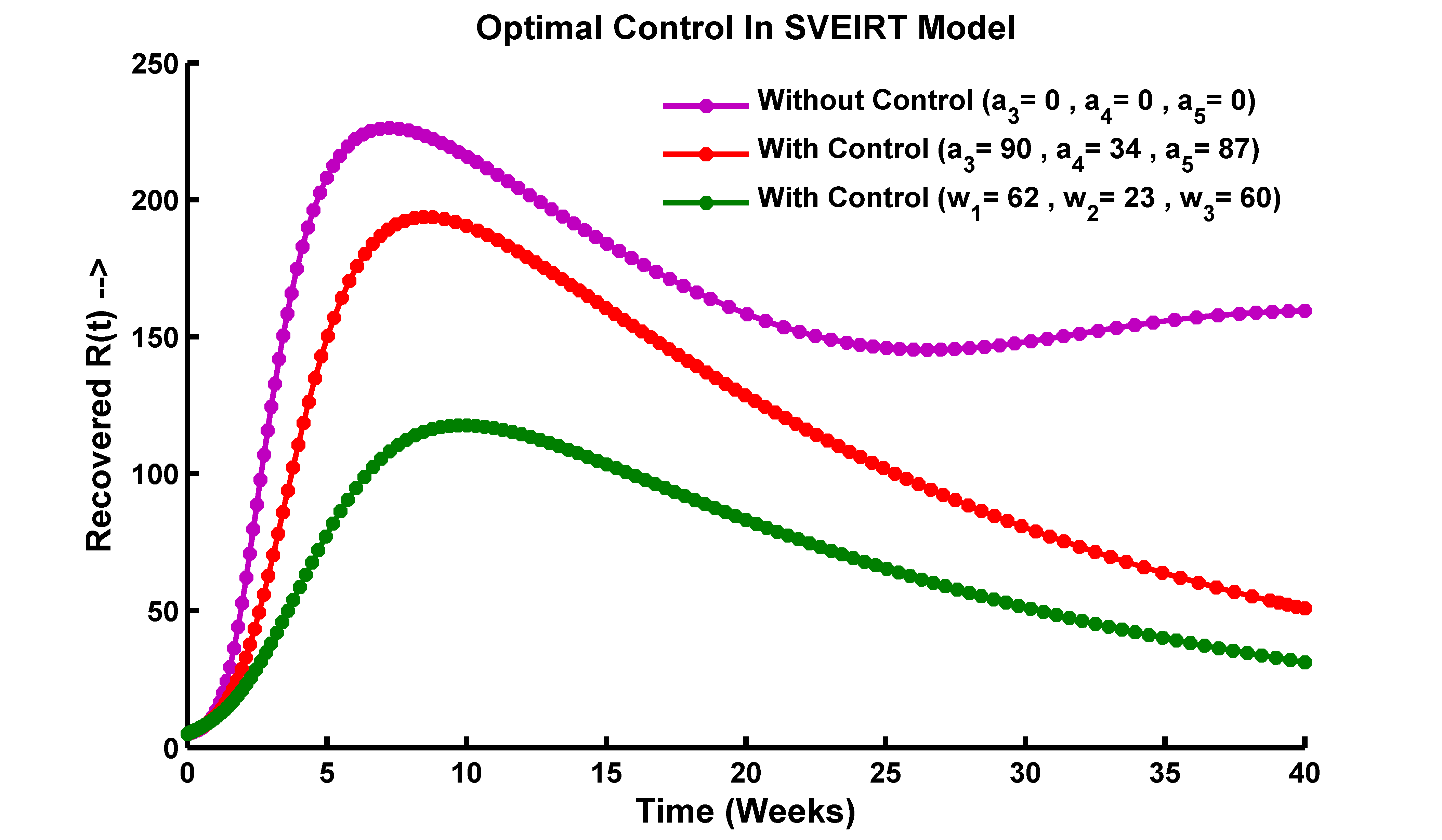}}
	\subfloat[]{\includegraphics[width=2.7 in]{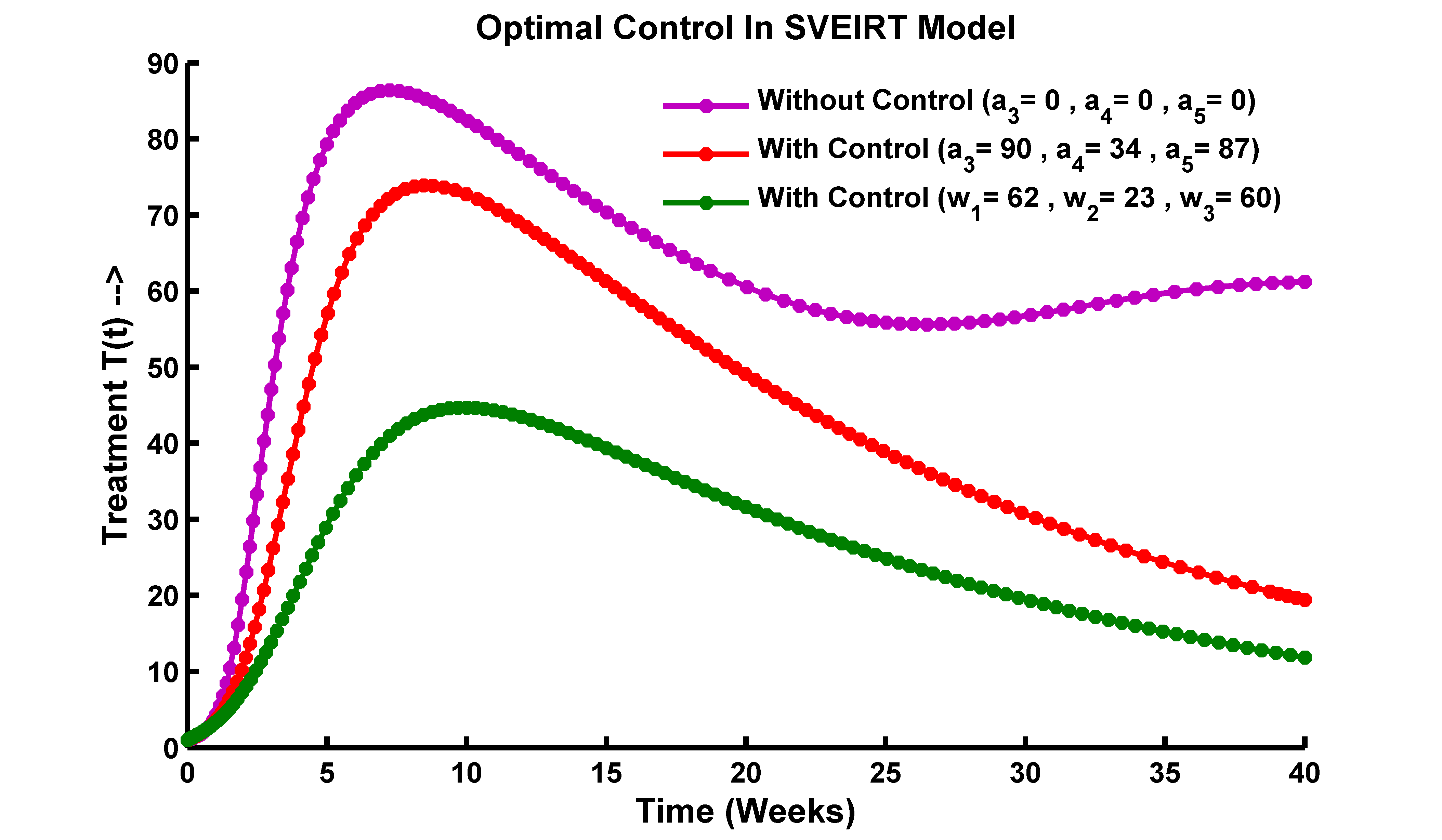}}
	\caption{Simulation of the model (a) $S(t)$ compartment,(b) $V(t)$ compartment (c) $E(t)$ compartment (d) $I(t)$ compartment (e) $R(t)$ compartment and (f) $T(t)$ compartment showing with control and without control with respect to weighted parameters $a_3,a_4,a_5$, where $a_3\in[0,62,90]$,\;$a_4\in[0,23,34]$,\;$a_5\in[0,60,87]$.}
	\label{Optimal-control vs a}
\end{figure}

Figure \ref{Optimal-control parameters} indicates that, while treatment is considered as control, it takes approximately 100\% input of treatment for 11–12 weeks before reaching their lower bound. If physical distancing is implemented at the proper rate, it takes about 100\% input of physical distancing for 9–10 weeks.
\noindent
Figure \ref{Optimal-control vs a}(a)-(f) demonstrates how control weights $a_3$, $a_4$, and $a_5$ affect exposed, infected, hospitalized, and recovered persons. In that situation, the following three situations have been taken into consideration:
\begin{enumerate}
	\item[(i)] implementation of no control strategy based on $w_1=0,\;w_2=0$ and $w_3=0$ for the case $a_3=0,\;a_4=0,\;a_5=0$
	\item[(ii)] implementation of three controls simultaneously, and in this scenario $w_1=0.45,\;w_2=0.45$ and $w_3=0.45$ for the case $a_3=90,\;a_4=34,\;a_5=87$
	\item[(iii)] implementation of three controls simultaneously, and in this scenario $w_1=0.65,\;w_2=0.65$ and $w_3=0.65$ for the case $a_3=62,\;a_4=23,\;a_5=60$
\end{enumerate}
\noindent
By examining these three scenarios, it can be observed that, in comparison to the situation where there is no control, there is a minor drop in the number of people who are exposed, infected, treated (hospitalized), and recovered when the appropriate control cost is raised. Additionally, a tiny shift is noted in the recovered class as a result of the tiny initial controls that were put in place, whose values are $a_3=62,\;a_4=23,\;a_4=60$.

From Figure \ref{Optimal-control vs a}, we observed that for scenario (ii) and (iii) when control weights are implemented in a effective rate, then susceptible and vaccinated population increases in notable rate. By applying scenario (ii) and (ii) exposed, infected and treated individuals decrease rapidly. Also for strategy (ii) implemented, there is an increase observed in recovered class than strategy (iii). In this case the implemented controls are $a_3=90,\;a_4=34,\;a_5=87$. After taking into account all of these findings, we may conclude that using mixed control strategies with low control cost weights, or high control efforts, is more successful than using single control or no control methods at all. This will assist in reducing the community's burden of sickness.
\noindent
Figure \ref{Optimal-control parameters} demonstrates the result of control profiles of implementation of $w_1,\;w_2,\;w_3$ controls simultaneously. This indicates that for $w_2=w_3=0$, control parameter $w_1$ (physical distancing, wearing masks) is implemented upto 88\% before reaching low bound. For the case $w_2=w_3=0.45$ the maximum implementation of control parameter $w_1=0.52$, while the maximum implementation of $w_1=0.22$ (in a low rate) for $w_2=0.65$ and $w_3=0.65$.  On the other hand, from Figure \ref{Optimal-control parameters}(b), we observe that maximum implementation of $w_2=0.58$ occurs when $w_1=w_3=0$, and for rising the values of $w_1$ and $w_2$ the implementation of control parameter $w_2$ (treatment rate) needed in lower rate.

The similar case is noticed for control parameter $w_3$. For the first case when $w_1$ and $w_2$ both 0, then enhanced treatment rate $w_3$ is implemented in maximum amount, where $w_3$ reaches its maximum value of 0.68 before reaching its lower bound. Moreover, the peak level of implementation of $w_2$ decreases according to the increasing of other control parameters $w_1$ and $w_2$.

More significantly, we discovered that there is a reversal relationship between control efforts and the weight of intervention costs and control parameters after considering the best control approach. The implemented levels of interventions are higher (or lower), accordingly, if the equivalent weight of cost of interventions, physical distance and treatment, is lower (or higher).

\section{Parameter Estimation and Model Validation}\label{Section-Parameter-Estimation-model-validation}
Understanding asymptotic behaviour and long-term outcomes is greatly aided by mathematical analysis of models. The values of a model's parameters have a significant impact on its results. Since models must handle illness data, a precise estimation of parameter values is required for reliable quantitative predictions within time intervals. It is preferable to use a systematic technique for the fitting when estimating several parameters. Various methods were employed to estimate the parameters in \cite{Stability Bound-5, Stability Bound-25}. The non-linear least-squares method, which is simple to apply, was used to calculate the parameters. In this least squares approach, the temporal coordinates of the data are taken to be exact, even though their corresponding y-coordinates (virions) could be distorted or noisy.

By utilizing numerical approaches, we can fit the non-linear system in \eqref{new_model} to the real influenza infected data, and then use the model to provide accurate predictions about the spread of the disease. We estimate the parameters utilizing the below-described methods in 
\cite{Stability Bound-8,Stability Bound-26,Stability Bound-27,Stability Bound-2-3}, in order to solve the data numerically.

The least squares approach is the most effective strategy to fit the suggested model using actual reported data. Fitting the model with actual infection data while minimizing the sum of square errors is the foundation of the least squares approach. If the total squares of the vertical distances between the actual data and the data predicted by the model are as little as possible, the model is said to be well-fitted. Most people refer to this distance as least squares error \cite{Hopf LEAST LHS PRCC-6, Hopf LEAST LHS PRCC-7, Hopf LEAST LHS PRCC-8}.
We use the following sum of squares error formula to fit the model:
\begin{align*}
	f(\phi,n)=\sum_{j=1}^{n}(Y_j-I(t_j))^2,
\end{align*}
where $Y_j$ represents the cumulative number of the real reported data for $j^{th}$ observation, $I(t_j)$ represents the model predicted cumulative data for $j^{th}$ observation, $n$ represents the total amount of available data, and $\phi$ represents the set of all model parameters. The formula is satisfied by the number of cumulative model predicted infected data:
\begin{align*}
	\frac{dI(t_j)}{dt}=\alpha E
\end{align*}
it is quite challenging to minimize $f(\phi,n)$ analytically. In a particular scenario, suppose the virions $K(t)$ are being fitted, with consideration of the given data $\{(t_1,\hat{K}_1),(t_2,\hat{K}_2),\cdots , (t_2,\hat{K}_n)\}$. Finding the set parameters $\theta$ that will result in the smallest feasible sum of squares error (SSE) is the fundamental challenge:
\begin{align*}
	SSE\{\min\;\theta\}=\sum_{i=1}^{n}\{K(t_i,\theta)-\hat{K}(t_i)\}^2,
\end{align*}
where $\hat{K}(t_i)$ represents the related data value at time $t_i$, and $K(t_i,\theta)$ represents the infection mixture at time $t$ with parameter $\theta$. This problem is clearly a non-linear least-squares problem since a solution depends on the parameter $\theta$ through a highly nonlinear system of differential equations. The important model parameters were estimated using the Mexico infection cases starting from October 15, 2021, and a total of 85 weeks of data from the CDC and WHO websites \cite{CDC, WHO}. We did this by using a Matlab function called \textbf{fminsearch}, which takes the least-squares error function $SSE(\theta)$ and a primary guess of the parameter value $\theta_0$. It then uses a direct search algorithm to obtain a minimum value of least-squares error. In Table \ref{tableparameter}, the calculated parameters are shown. In the same approach, by collecting total 85 weeks data infected from Colombia, Italy and South Africa the parameters are estimated for these countries which all are presented in Table \ref{table-param-values-4-countries}. In order to run the Matlab package, we have taken the starting population size into consideration as:
\begin{align*}
	S(0)= 500,\;\; V(0)=1,\;\;E(0)=1,\;\;I(0)=1,\;\;R(0)=0,\;\;\text{and}\;\;T(0)=0.
\end{align*}
where, these initial guess are estimated from \cite{Stability Bound-25}. 
Finding residuals provides justification for another method of confirming the model's fitness. Thus, by calculating the residuals, the model's fitness with the actual data is confirmed in this instance. It is defined that the residuals are
\begin{align*}
	\text{residuals}= \{Y_j-I(t_j)|j=1,2,\cdots,n\}
\end{align*}
where, $Y_j$ is the $j^{th}$ week cumulative information about infections and $I(t_j)$ represents the same week's cumulative infection data was anticipated by the model. Given a random distribution of the residuals, we can conclude that the fitness level is respectable \cite{Stability Bound-2-3}.
\begin{figure}[H]
	\centering  
	\subfloat[]{\includegraphics[width=2.7 in]{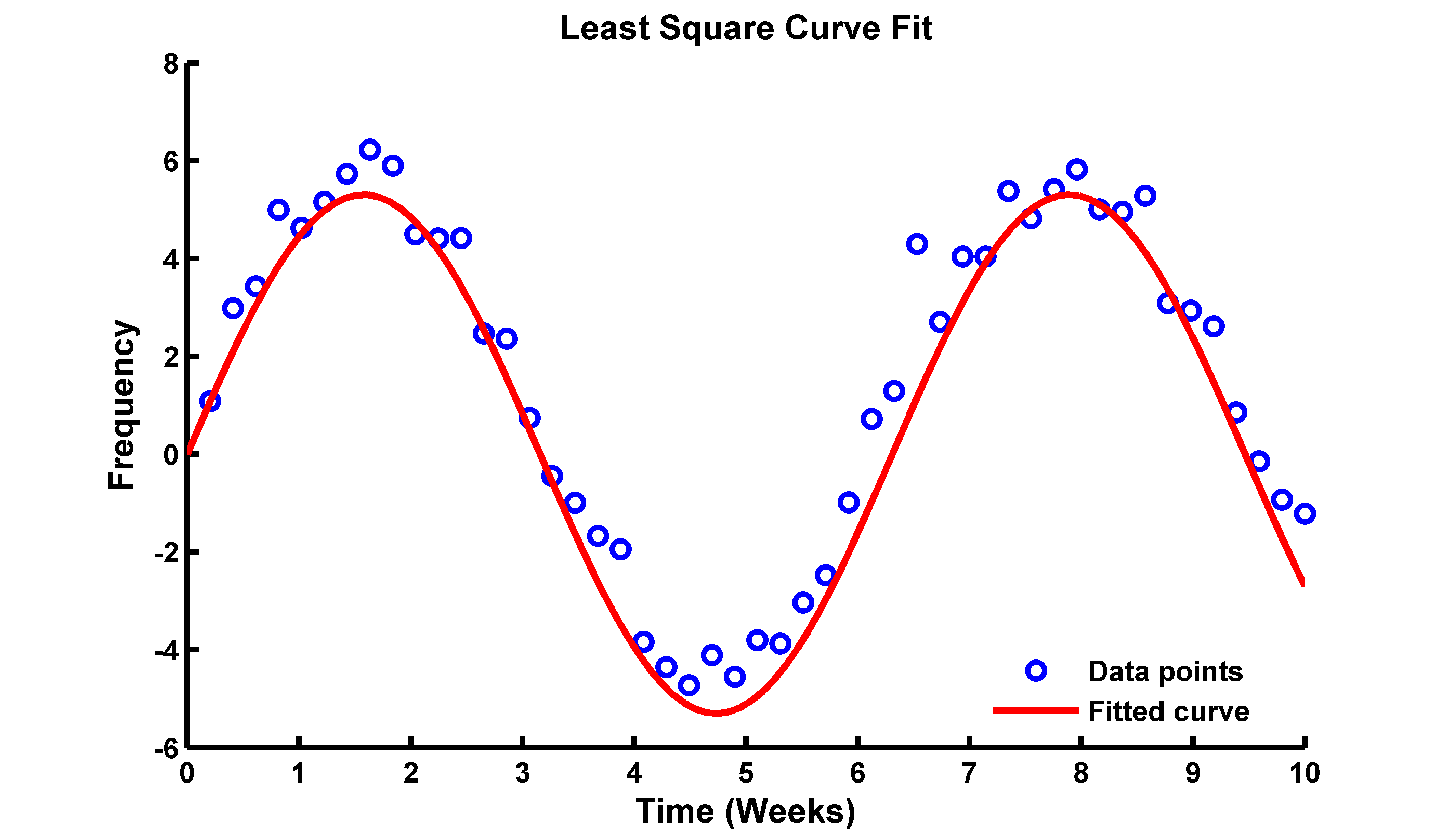}}
	\subfloat[]{\includegraphics[width=2.7 in]{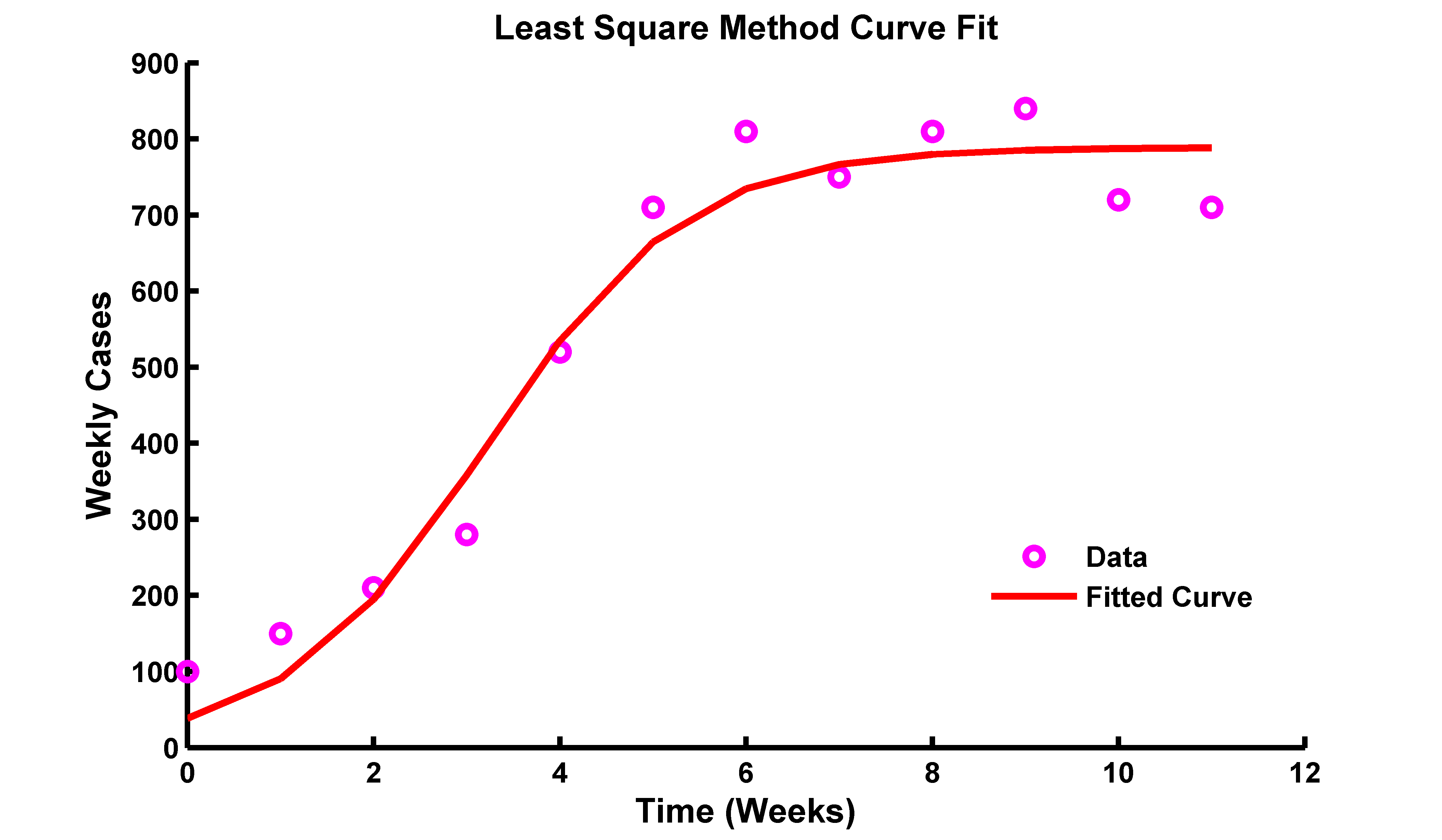}}\\
	\subfloat[]{\includegraphics[width=2.7 in]{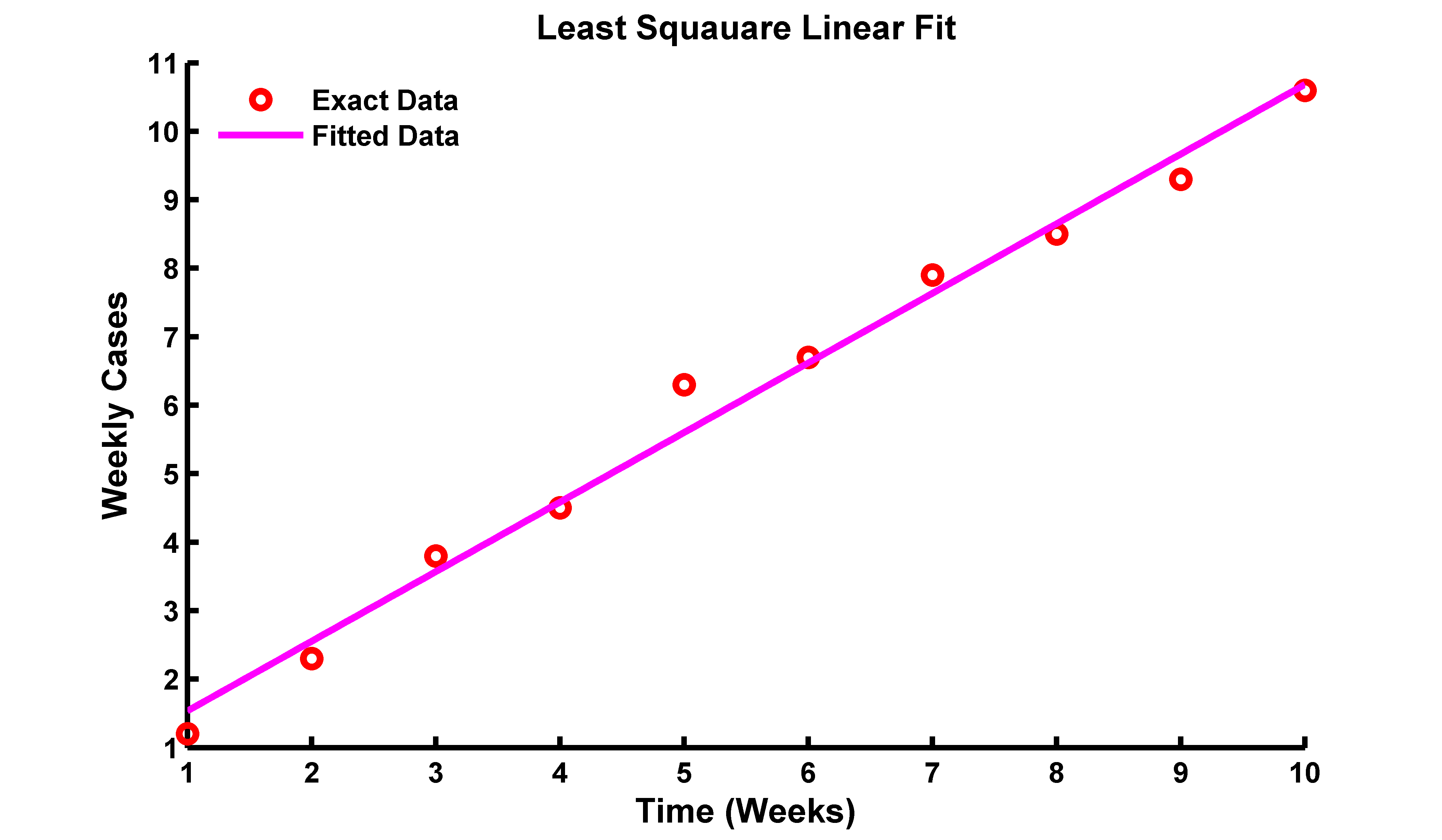}}
	\subfloat[]{\includegraphics[width=2.7 in]{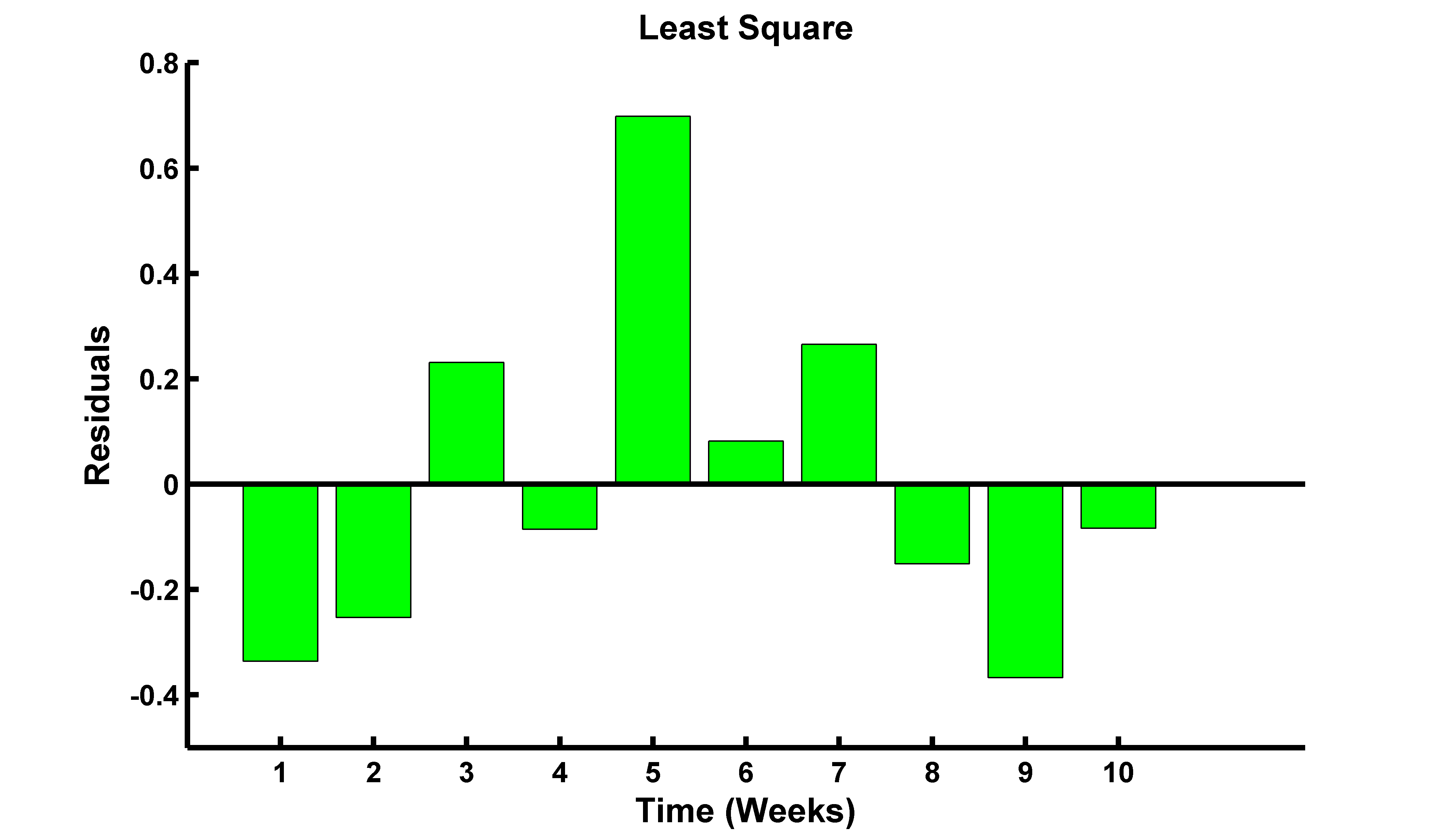}}
	\caption{(a) Least square polynomial fitting (b) Least square curve fitting  (c) Least square linear fitting  and (d) Residuals for exact vs fitted data.}
	\label{LS-fit-curve-Basic}
\end{figure}
\noindent
Figure \ref{LS-fit-curve-Basic}(a) refers to fitting a polynomial curve to the data points using the least squares method. The polynomial curve is chosen such that it minimizes the sum of the squared differences between the data points and the curve. We see that for time intervals 10 weeks, the real data and model observed data is fitted significantly. A higher-degree polynomial can provide a more flexible fit to the data, capturing more complex patterns.

For Figure \ref{LS-fit-curve-Basic}(b), in this case , the least squares method fits a general curve function which isn't just a polynomial to the data points. The curve function can be chosen based on the specific nature of the data. It allows for more flexibility in capturing the underlying trend or pattern, beyond what can be achieved by a simple polynomial fit. We see that for time intervals 12 weeks, the real data and model observed data is matched effectively.

Figure \ref{LS-fit-curve-Basic}(c) involves fitting a straight line to the data points using the least squares method. The linear fit assumes a linear relationship between the independent and dependent variables. It is the simplest form of fitting and is often used when the data appears to follow a linear trend. We have considered total 10 weeks data and these are fitted well.

The Figure of residuals \ref{LS-fit-curve-Basic}(d) represents the differences between the observed data points and the corresponding values predicted by the fitted curve observed in 10 weeks.  A small residual in weeks 3,4,6,8,10 indicates a good fit, meaning that the observed data points are close to the values predicted by the model. Conversely, a average residual in weeks in 1,2,5,9 indicates a likely fit, suggesting that the observed data points deviate slightly from the predicted values. Patterns or trends in the residuals can indicate systematic errors in the model or potential areas where the model fails to capture important features of the data. By examining the Figure \ref{LS-fit-curve-Basic} we can assess the overall accuracy and reliability of the curve fit.
\begin{table}[H]
	\begin{center}
		\caption{Model parameters values.}
		\scriptsize
		\label{tableparameter} 
		\begin{tabular}{|l|l|l|l|}
			\hline\noalign{\smallskip}
			\textbf{Notation} & \textbf{Definition} & \textbf{Value}  & \textbf{Source}  \\
			\noalign{\smallskip}\hline\noalign{\smallskip}
			$ \alpha $ & Transition rate from $E$ to $I$  &  $0.75$ week$ ^{-1} $  & \cite{Stability Bound-5,Stability Bound-13} \\
			$\Lambda$ & Recruitment rate in $S$ class & $5\times 10^2$ week$ ^{-1} $ & \cite{Stability Bound-14,Stability Bound-17} \\
			$ \beta_1 $ & Transmission rate from contact with $E$ to class $S$  & $ [0.0045,0.0055 ]$ week$ ^{-1} $ & \cite{Stability Bound-18,Stability Bound-20}\\
			$ \beta_2  $ & Transmission rate from contact with $I$ to class $S$&  $ [0.0045,0.0055 ]$ week$ ^{-1} $  & \cite{Stability Bound-21,Stability Bound-23} \\
			$\gamma$ & Recovery rate of $I$  & $ 0.65 $ week$ ^{-1} $ & \cite{Stability Bound-19,Stability Bound-24}\\
			$\gamma_1$ & Treatment progression rate of $I$ &  $ 0.25$ week$ ^{-1} $ & \cite{Stability Bound-8, Stability Bound-9} \\
			$ \mu $ & Natural death rate & $ 5 \times 10^{-2}$ week$ ^{-1} $ & \cite{Stability Bound-5,Stability Bound-8}\\
			$ \delta $ & Disease induced death rate & $0.3$ week $^{-1} $ & \cite{Stability Bound-13,Stability Bound-19}\\
			$ \lambda $ & Vaccine inefficiency rate & $0.55$  & \cite{Stability Bound-14,Stability Bound-19}\\
			\noalign{\smallskip}\hline
		\end{tabular}
	\end{center}
\end{table}

\begin{table}[H]
	\begin{center}
		\caption{Estimated parameter values for model \eqref{new_model} by analyzing Mexico, Italy,  and South Africa data.}
		\scriptsize
		\label{table-param-values-4-countries}
		\begin{tabular}{|c|c|c|c|c|c|}
			\hline\noalign{\smallskip}
			\multirow{2}{*}{\textbf{Parameter}}&\textbf{Mexico}& \textbf{Italy} & \textbf{South Africa}&\multirow{2}{*}{\textbf{Sources}} \\
			\cline{2-5}
			&\textbf{Values}&\textbf{ Values}&\textbf{Values}& \\
			\noalign{\smallskip}\hline\noalign{\smallskip}
			$\alpha$&0.75&0.67&0.78&Estimated\\
			$\beta_1$&0.0055&0.0053&0.0075&Estimated\\
			$\beta_2$&0.0055&0.0061&0.0081&Estimated\\
			$\gamma$&0.65&0.61&0.63&Estimated\\
			$\gamma_1$&0.25&0.31&0.35&Estimated\\
			$\lambda$&0.55&0.52&0.42&Estimated\\
			$\mu$&0.05&0.03&0.03&Estimated\\
			$\delta$&0.3&0.27&0.29&Estimated\\
			$\varepsilon$&0.45&0.41&0.44&Estimated\\
			\noalign{\smallskip}\hline
		\end{tabular}
	\end{center}
\end{table}

\begin{figure}[H]
	\centering  
	\subfloat[]{\includegraphics[width=2.8 in]{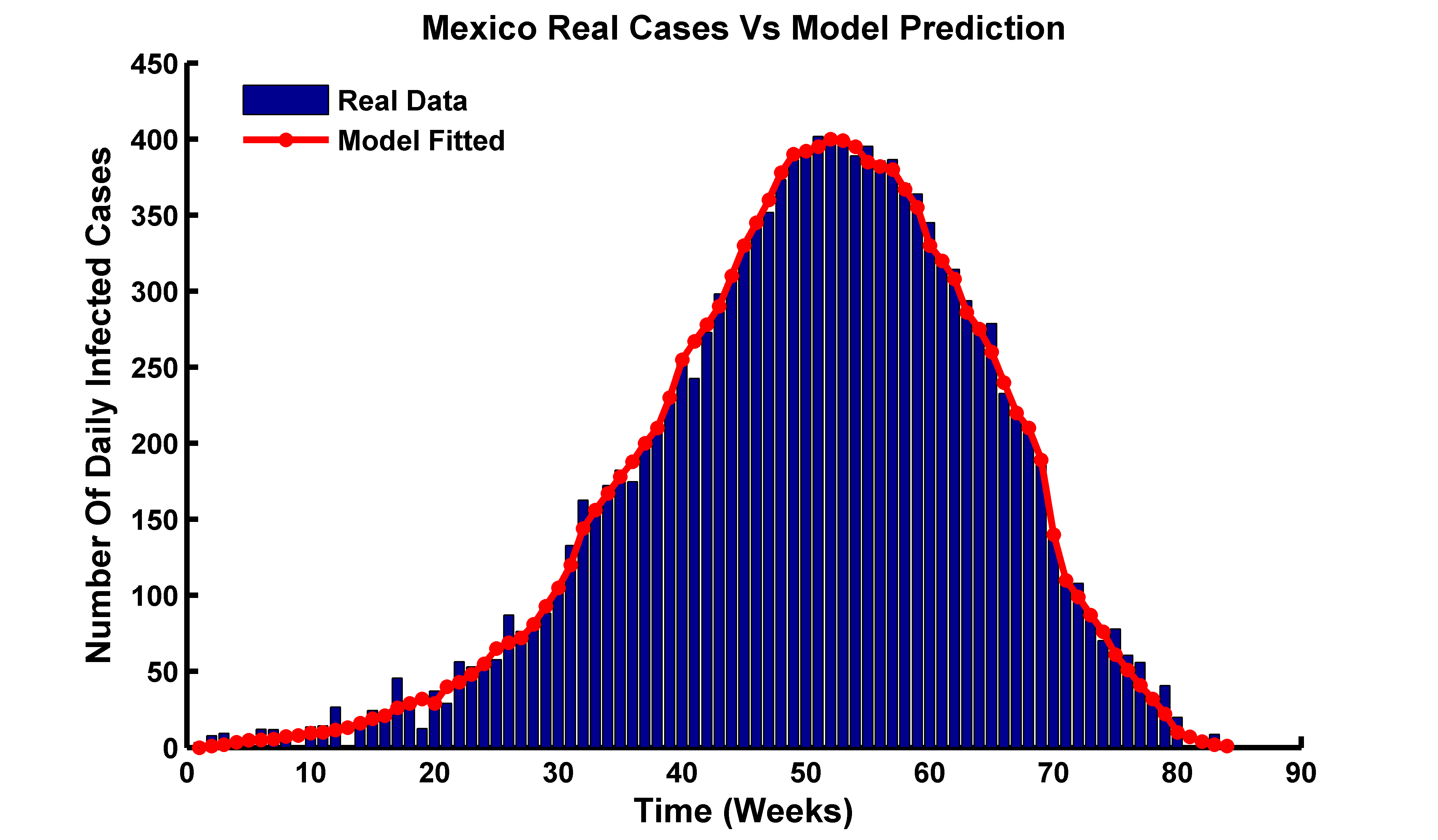}}
	\subfloat[]{\includegraphics[width=2.8 in]{LScumini.png}}\\
	\subfloat[]{\includegraphics[width=2.8 in]{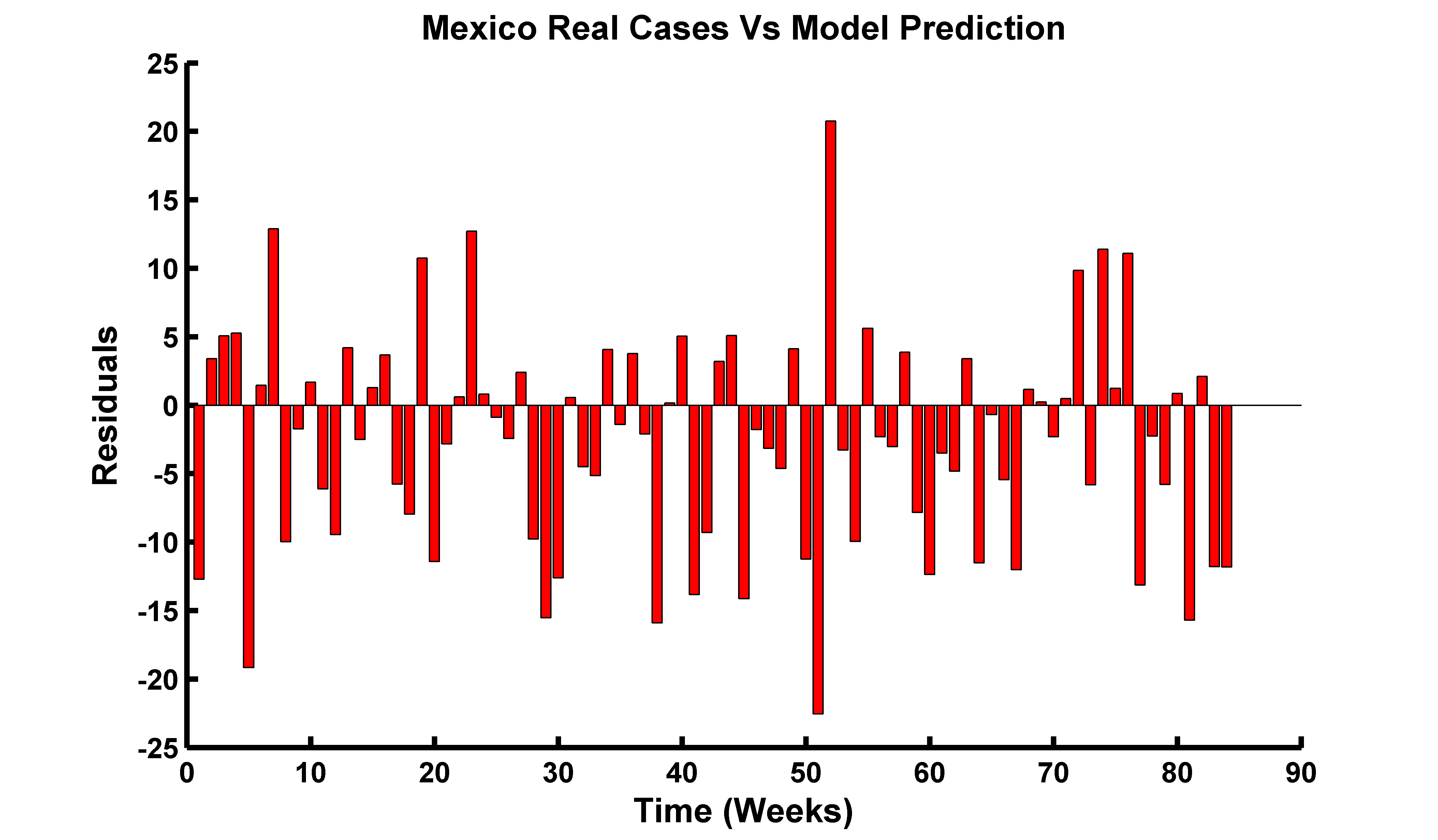}}
	\caption{(a) Fitting model data with weekly cases (b) Residuals of the fit (c) Time series of new infected cases of Influenza and (d) Weekly cumulative number of cases of Mexico data in 85 weeks starting from 15 October 2021.}
	\label{LS-fit-curve-residuals}
\end{figure}

We examine the real-world scenario of a Mexico influenza outbreak for a total of 85 weeks beginning on October 15, 2021, in order to validate the model. In Figure \ref{LS-fit-curve-residuals}(a), we display the afflicted population as a bar graph. The model is fitted to the total number of cases that have been infected, as shown in Figure \ref{LS-fit-curve-residuals}(b). The fit residuals are displayed. in Figure \ref{LS-fit-curve-residuals}(c), which demonstrates the tiny and erratic residuals.  The residuals' randomness indicates that the fitness is optimal.
\noindent
Additionally, we have confirmed that the model suggested for the current trends is based on the data from Mexico, as shown in Figure \ref{LS-fit-curve-residuals}(a). Since the model's proportion is well defined and the model's data is corroborated by actual reported data from \cite{CDC} in Mexico, it is clear from the model fitting that our model is well fitted. We can forecast the actual situations using the estimated model parameters shown in Table \ref{table-param-values-4-countries}. Figure \ref{LS-fit-curve-residuals}(a) shows that in weeks 0 to 20 the number of infected population lies around 50. After that, the outbreak occurs, and infected class increased gradually. around 40 to 60 weeks the maximum infection reported and the pick value reported around 400. This data is estimated for per 1000 populations in Mexico. After 60 weeks, the infection behaviour reverses and converges to the 0 level. That indicates the control of the disease. The actual data and the model-predicted data are substantially closer at every week. This indicates the relevance of our model \ref{new_model}.

\subsection{Estimation of $\mathcal{R}_0$ from Actual Data}
\begin{figure}[H]
	\centering  
	\subfloat[]{\includegraphics[width=2.7 in]{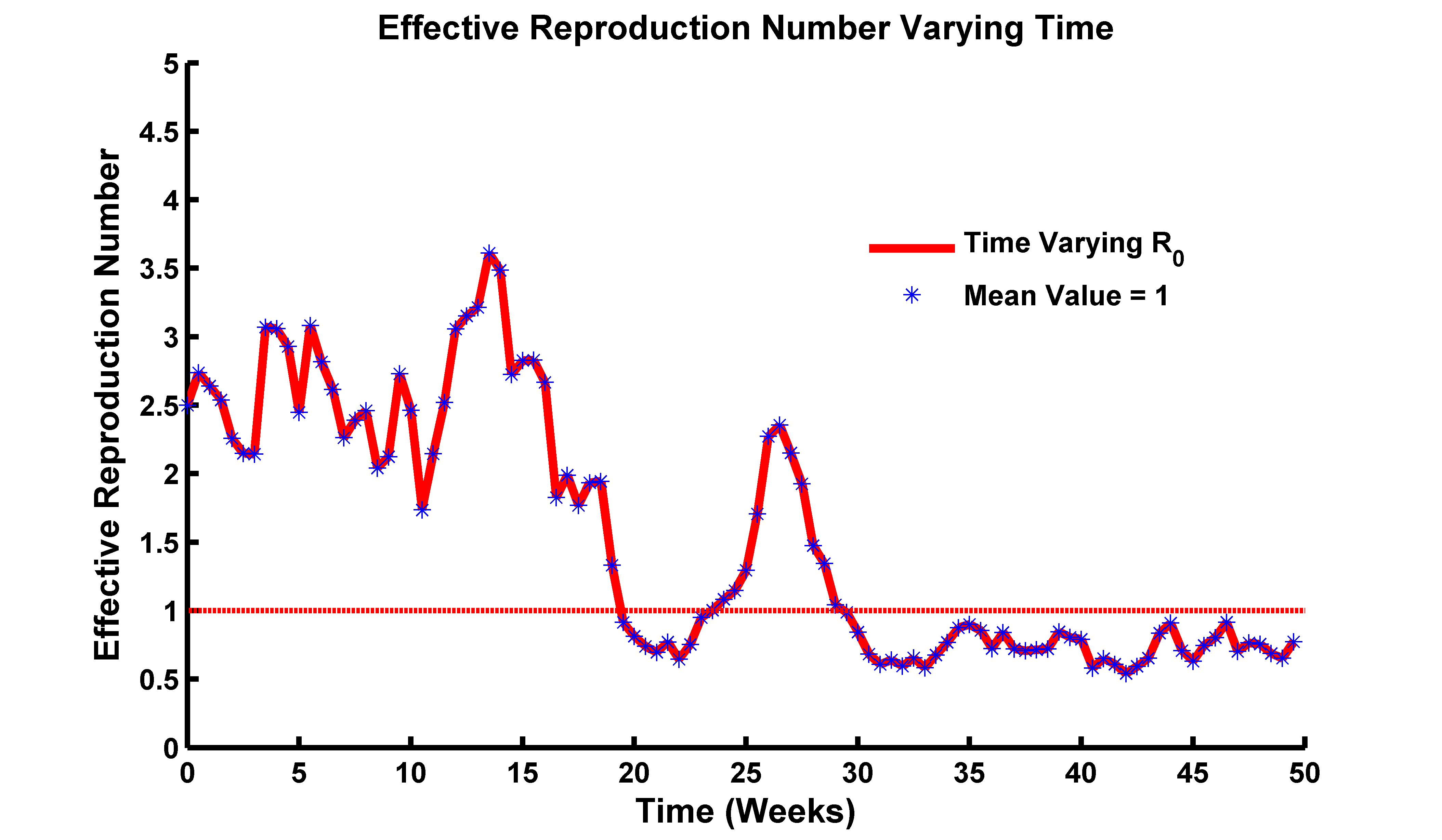}}
	\subfloat[]{\includegraphics[width=2.7 in]{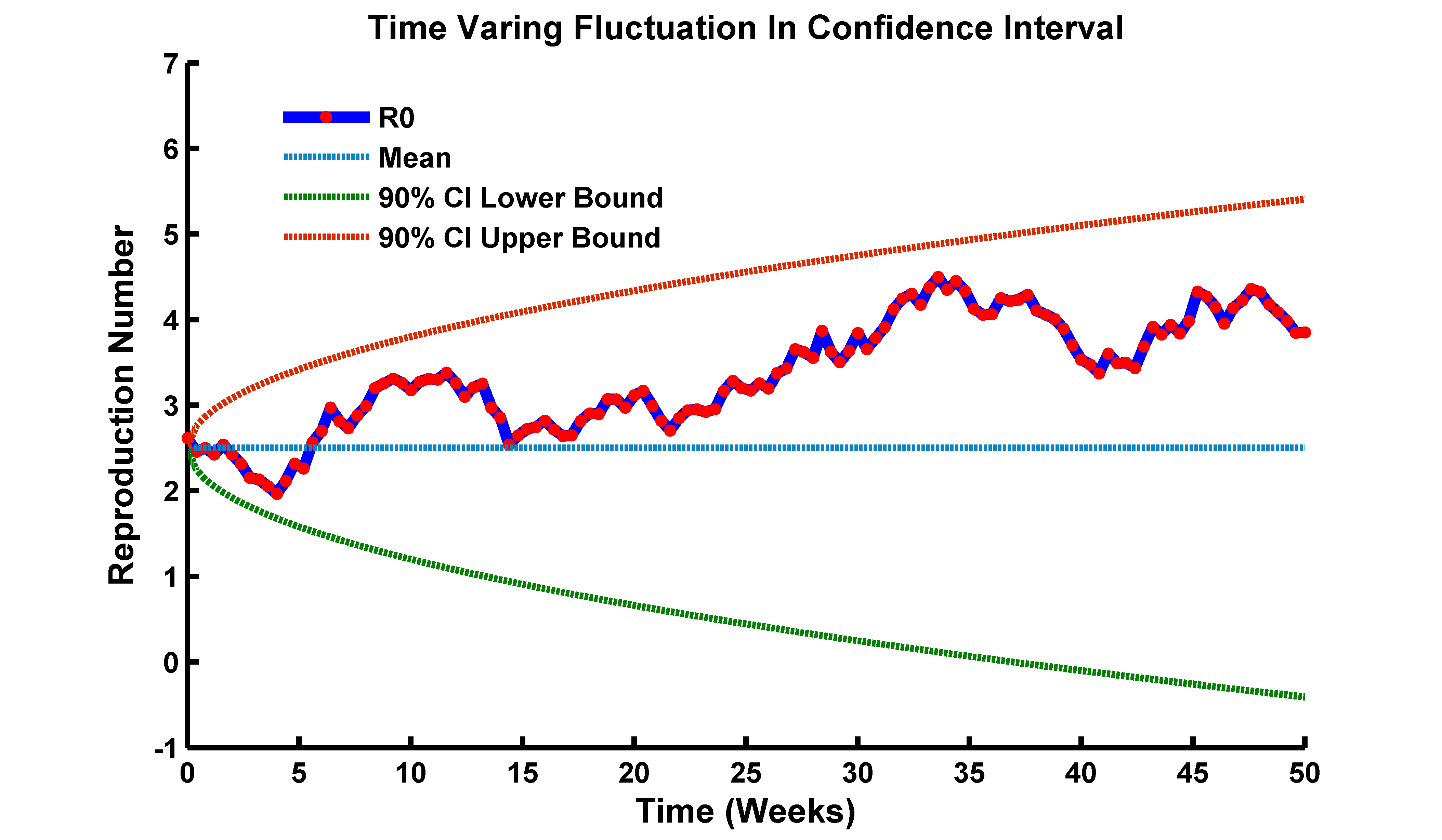}}\\
	\subfloat[]{\includegraphics[width=2.7 in]{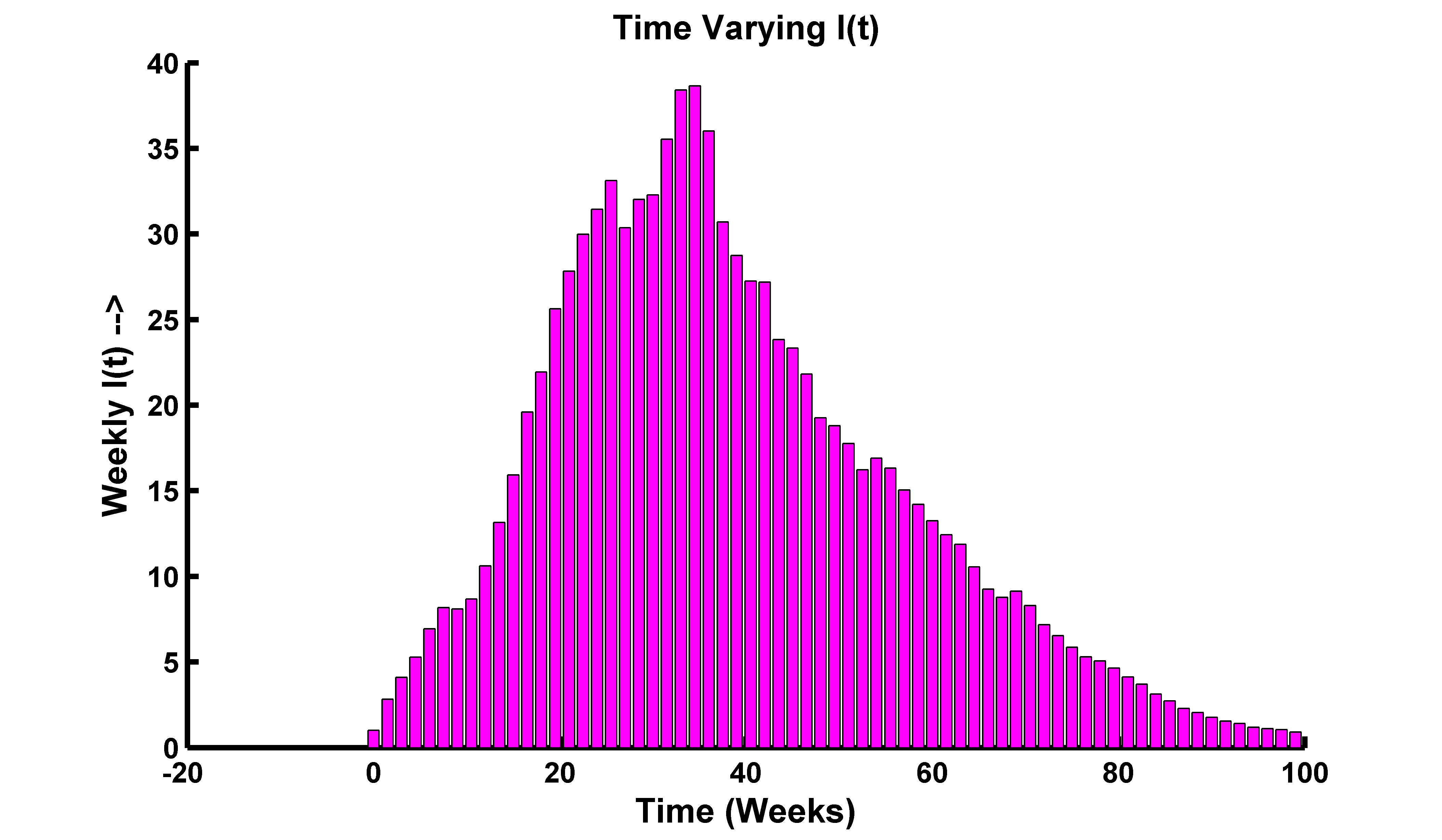}}
	\subfloat[]{\includegraphics[width=2.7 in]{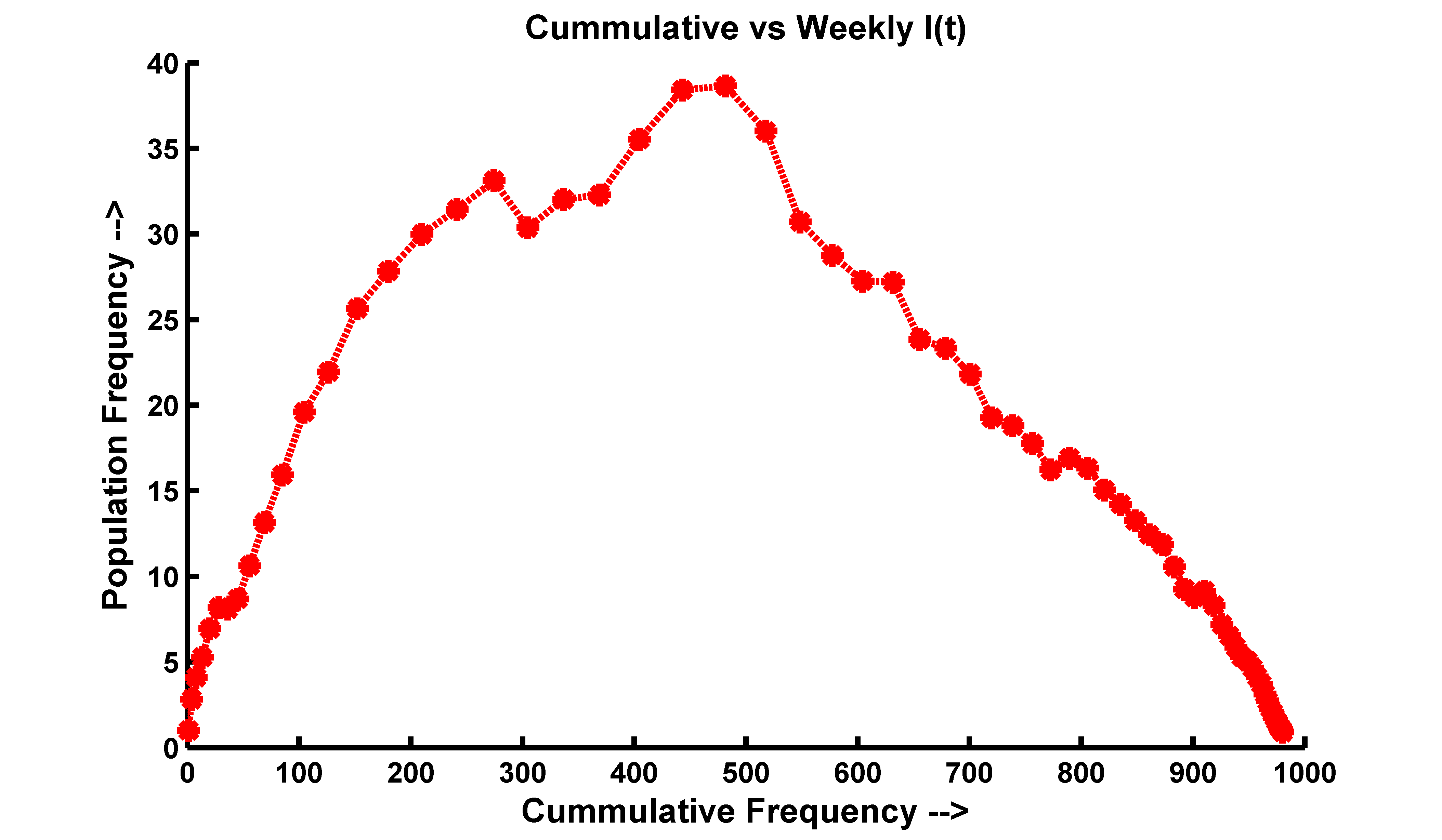}}
	\caption{(a) Time series effective reproduction number (b) Time series effective reproduction number in confidence intervals (c) Time series of new infected cases of Influenza and (d) Weekly number of cases vs cumulative number of cases from 15 October 2021 to 21 July 2022.}
	\label{effective-reproduction-number}
\end{figure}

We will determine the estimated value of $\mathcal{R}_0$ from the reported data in this section of the study, up to which the time series of the affected data remain exponential. In order to calculate $\mathcal{R}_0$ from the disease's early development phase, we employed the techniques described in \cite{Marcheva Book,Stability Bound-2-3}. At the onset of the illness, we presume that cumulative case ($Q(t)$) varies as $\exp (\Lambda t)$ that means $Q(t)\propto \exp(\Lambda t)$. The population's exposure and infection rates also fluctuate as $\exp(\Lambda t).$ Therefore
\begin{align}\label{effective-R0-E0-I0}
	E\sim E_0\exp(\Lambda t),\;\;\;\;I\sim I_0\exp(\Lambda t)
\end{align}
where, $E_0$ and $I_0$ are unchanging. The population's constant susceptibility number is provided by $\displaystyle S_0=\frac{\Lambda}{(\mu+\phi)}.$ Replacing \eqref{effective-R0-E0-I0} in third equation of model \eqref{new_model}, we get,
\begin{align}\label{effective-R0-beta_12}
	(\beta_1E_0+\beta_2I_0)\frac{\Lambda}{\mu+\phi}=(\alpha+\mu)E_0
\end{align}
Putting the value of $\beta_1$ and $\beta_2$ from \eqref{effective-R0-beta_12} in equation \eqref{new_model}, We derive the basic reproduction number expression $\mathcal{R}_0$ as shown in the format:
\begin{align*}
	\mathcal{R}_0=\frac{\Lambda\alpha\beta_2}{(\mu+\phi)(\alpha+\mu)(\mu+\delta+\gamma+\gamma_1)}+\frac{\alpha\mu\left(\Lambda+\alpha+\mu-\frac{\beta_1\Lambda}{\mu}\right)(\Lambda+\mu+\delta+\gamma+\gamma_1)}{(\mu+\phi)(\alpha+\mu)(\mu+\delta+\gamma+\gamma_1)}
\end{align*}
For estimating $\mathcal{R}_0$, first, we must determine the infection's force $(\Lambda)$. Weekly total of new cases $q(t)$ fluctuates according to the weekly total of cumulative cases $(Q(t))$ as $q(t)\sim\Lambda Q(t)$.\\
By graphing weekly new cases against weekly cumulative cases, we are able to calculate the infectious force $(\Lambda)$. From the Figure \ref{effective-reproduction-number}(a) we determine the cumulative case's threshold value, for which the number of new cases exhibits exponential increase. The linear regression curve is fitted using the least squares approach. The force of infection $(\Lambda)$ is the slope of the regression line. Based on the regression line's slope, we have $\Lambda=500\times 0.02\;\;week^{-1}.$ The calculated parameters are presented in Table \ref{tableparameter} and Table \ref{table-param-values-4-countries}. Through the examination of data from South Africa, Mexico, and Italy, we have calculated the fundamental reproduction number $\mathcal{R}_0\in [0.5,3.75]$ with  lower value 0.5 and upper value 3.75.

Thus we get a range of $\mathcal{R}_0$ from the initial phase of the Influenza infection. From the range of $\mathcal{R}_0$ lies always greater than unity. However, we should cherish unity less than the eradication of the influenza virus from the population. Consequently, the range of $\mathcal{R}_0$ indicates that the influenza virus spreads throughout the population. Furthermore, as $\mathcal{R}_0$ lies at the maximum time greater than unity, the influenza virus continues to circulate in the population \cite{Stability Bound-25}.
\noindent
From Figure \ref{effective-reproduction-number}(a), we see that $\mathcal{R}_0$ value goes pick in week 14; while in week 10 to 15 $\mathcal{R}_0$ value fluctuates high. After 20 weeks, the fluctuation of $\mathcal{R}_0$ becomes lower. In week 27, $\mathcal{R}_0=2.4$, after that it converges to 0.5. By varying time $\mathcal{R}_0$ also vary with noise.\\
Meanwhile in Figure \ref{effective-reproduction-number}(b), the time varying $\mathcal{R}_0$ is presented in confidence interval. $\mathcal{R}_0$ fluctuate according to the fluctuation of exposed and infected cases. We have considered data from 15 October 2021 to 21 July, 2022; all data collected from \cite{CDC, WHO}. The figure reveals that, among 0 to 25 weeks the fluctuation of $\mathcal{R}_0$ is low and lies in $[2,3.75]$. After 25 weeks, $\mathcal{R}_0$ varies in the range $[3.75, 5.2]$. The pick value occurs in weeks 30 to 35, also in weeks 45 to 50.\\
By analyzing the available data, when plotting time-varying infected cases (Figure \ref{effective-reproduction-number}(c)), the graph shows how the number of people who have had influenza over a given period of time (in weeks) has changed. It shows the progression of the disease from 0 to 20 weeks is very high. It provides insights into the rate at which new infections occur from week 20 to 40 rapidly. The number of cases of infected evolve over time as in week 35 to 43, the infected cases are at pick level. The number of new infections reported every interval of time, such as weekly, is usually displayed in the time-varying infected cases plot. It helps us to visualize the trends and patterns in the spread of the disease, including any spikes or declines in infection rates for total 100 weeks. After week 43, the infection rate decreases gradually and along 60 to 100 weeks the disease can be controlled effectively.\\
On the other hand, the time-varying cumulative infected cases presented in Figure \ref{effective-reproduction-number}(d). The total number of infected cases varies over time, as seen by this graphic. As time goes on (between weeks 20 and 60), the total number of cases is growing quickly. After 60 weeks, the cumulative cases grows slightly as the total infected cases reducing gradually. Thus we give insight that massive outbreak occurs along week 35 to 50 according to the available data \cite{CDC, WHO}.
\subsection{Effective Reproduction Number}
One of the most important factors in the dynamics of disease transmission is basic reproduction rate. The average number of secondary infections that an infected host causes during its lifetime is known as the basic reproduction number. Disease spreads swiftly in the early stages among the population, but it slows down after reaching its peak. Reproduction is therefore not always continuous. Our current research focuses on time-varying reproduction numbers, or numbers for reproduction determined once a week. Adequate reproduction number is the term used to describe this kind of reproduction number $\mathcal{R}_0(t)$ \cite{Marcheva Book, Stability Bound-25}. The study's ability to provide information about the disease and appropriate preventive strategies to control it will depend on the effective reproduction numbers. We employ a formula \cite{Marcheva Book, Stability Bound-25}
\begin{align}\label{effective-R0-c(t)}
	\mathcal{R}_0(t)=\frac{c(t)}{\int_{0}^{\infty}c(t-\lambda)h(\lambda)d\lambda}
\end{align}
in order to calculate the effective reproduction number, where $c(t)$ indicates fresh cases at $t^{th}$ week and $h(\lambda)$ illustrates the basic interval distribution. Permit the exposed, diseased class to go at the pace $b_1=\alpha+\mu,\;\;b_2=\mu+\delta+\gamma+\gamma_1$, respectively. Let, $b_1e^{-b_1t},\;\;b_2e^{-b_2t}$ be the sum of the distribution of generation intervals, the following formula is provided:
\begin{align}\label{effective-R0-b1-b2}
	h(t)=\sum_{i=1}^{2}\frac{b_1b_2e^{b_it}}{\prod_{j=1,\;j\neq i}^{2}(b_j-b_i)}
\end{align}
with mean $\displaystyle T=\frac{1}{b_1}+\frac{1}{b_2}$. The above formula is valid when $\Lambda>\min\{-b_1,-b_2\}.$ Using new cases and formulas from equations  \eqref{effective-R0-E0-I0},  \eqref{effective-R0-beta_12}, \eqref{effective-R0-c(t)}, \eqref{effective-R0-b1-b2}, estimating the effective reproduction number is possible. Figure \ref{effective-reproduction-number}(a) and (b) displays the values of the effective reproduction number. From the figures, value of $\mathcal{R}_0(t)$ converges to unity and oscillates around 2.5. 
Meanwhile, $\mathcal{R}_0(t)$ decreased between weeks 14 and 21 from 3.75 to 0.7.\\
As a result, we now know the weekly basic reproduction number. The reproduction number is quite high in the beginning, and as time goes on, the values of $\mathcal{R}_0(t)$ on a weekly basis drop. The value of $\mathcal{R}_0(t)$ after 30 weeks is around unity until week 50. We need to maintain its value below unity in order to stop the influenza virus from spreading \cite{Stability Bound-25}.

\subsection{Relative Influence of Parameters}\label{Subsection-Relative-Influence}
We have performed the relationship of parameters in this section for fixed values of $\mathcal{R}_0$ which are presented by contour plots. The goal of this contour plot is to visually illustrate the contours or lines of constant values for a particular parameters or relationship in particular ranges. By examining the contour lines in this study, we can observe patterns, trends, and variations in the relationship between the two parameters \cite{Stability Bound-25, Bifurcation of R0-3}.

\begin{figure}[H]
	\centering  
	\subfloat[]{\includegraphics[width=2.7 in]{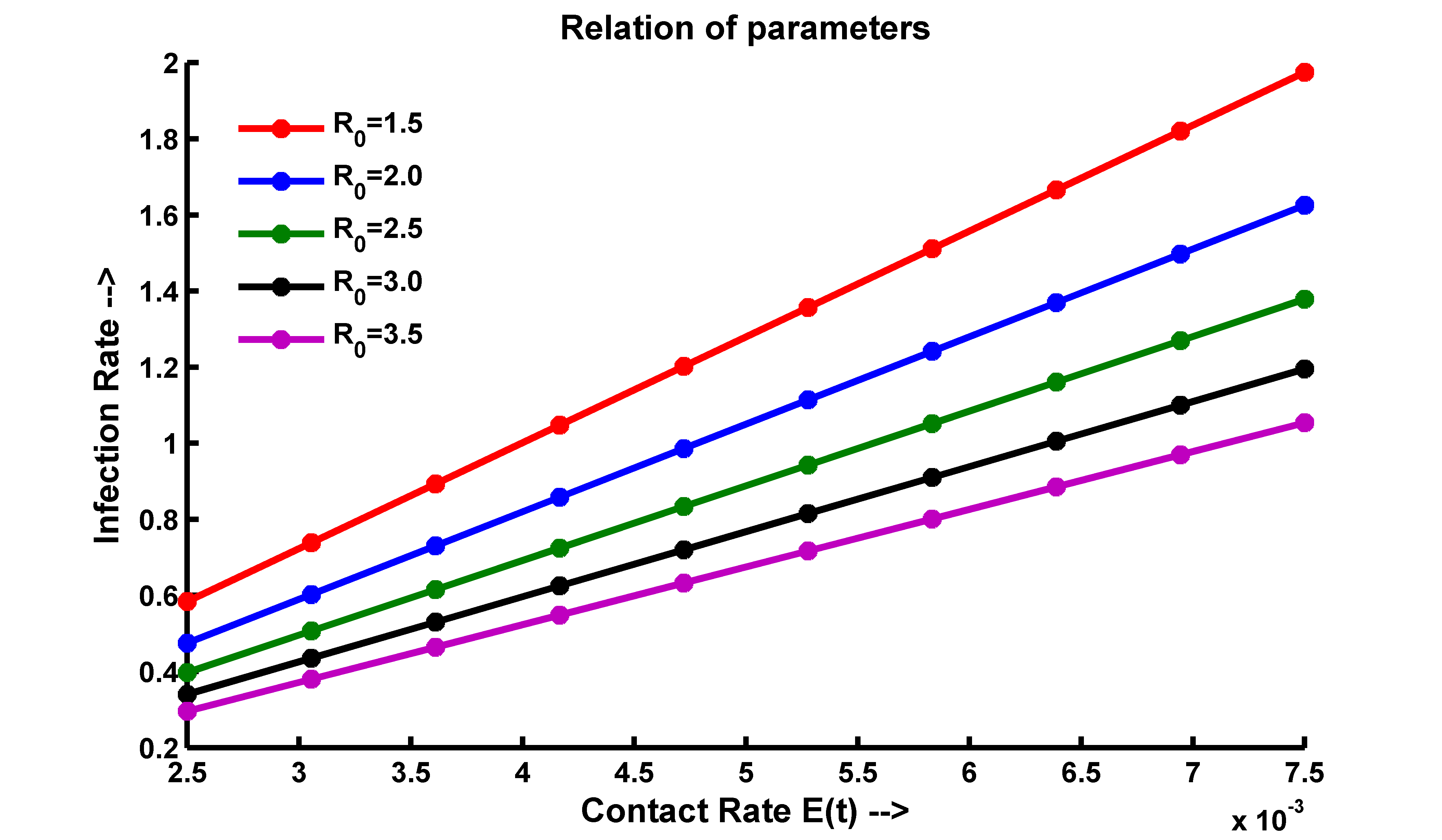}}
	\subfloat[]{\includegraphics[width=2.7 in]{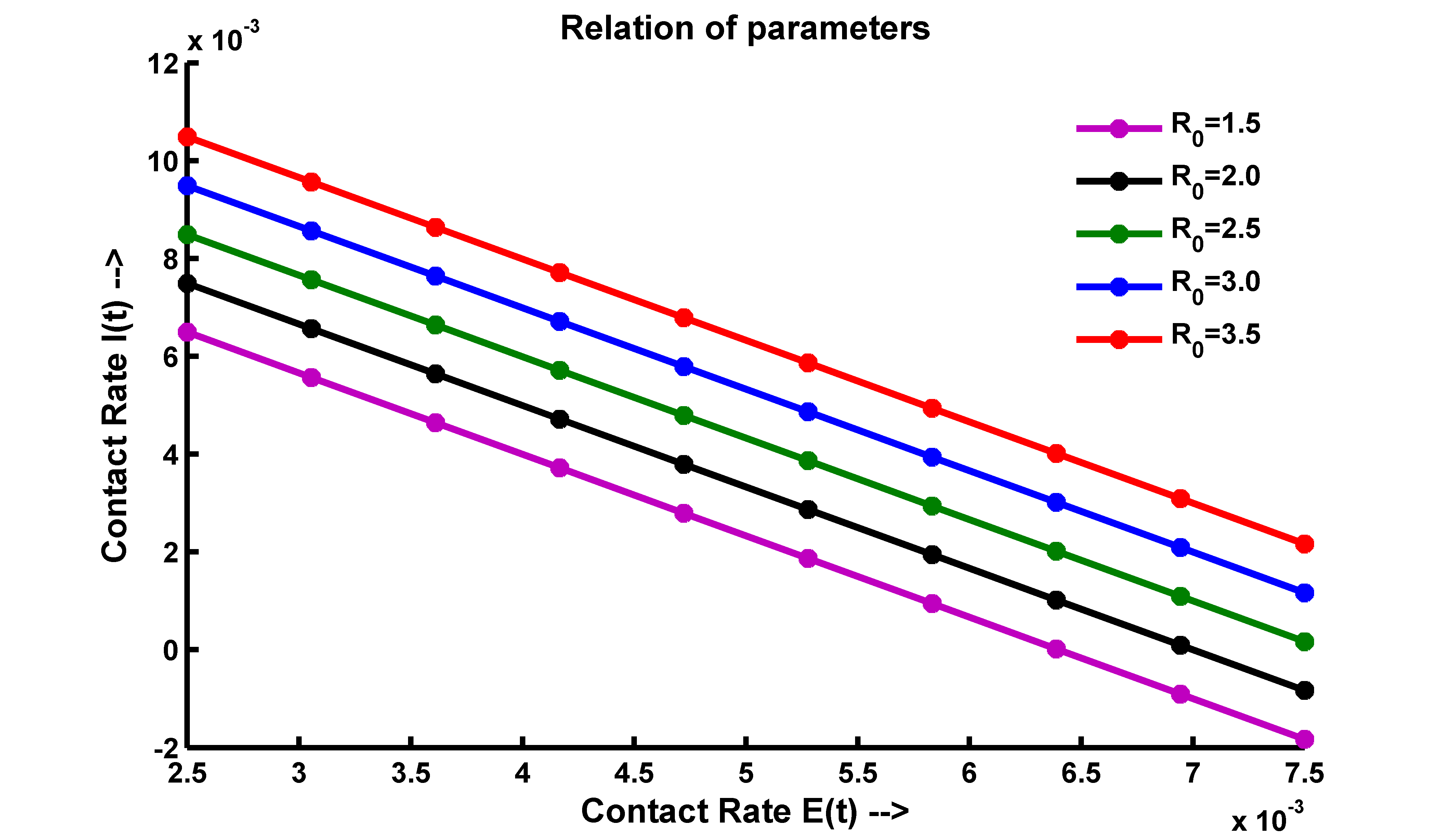}}\\
	\subfloat[]{\includegraphics[width=2.7 in]{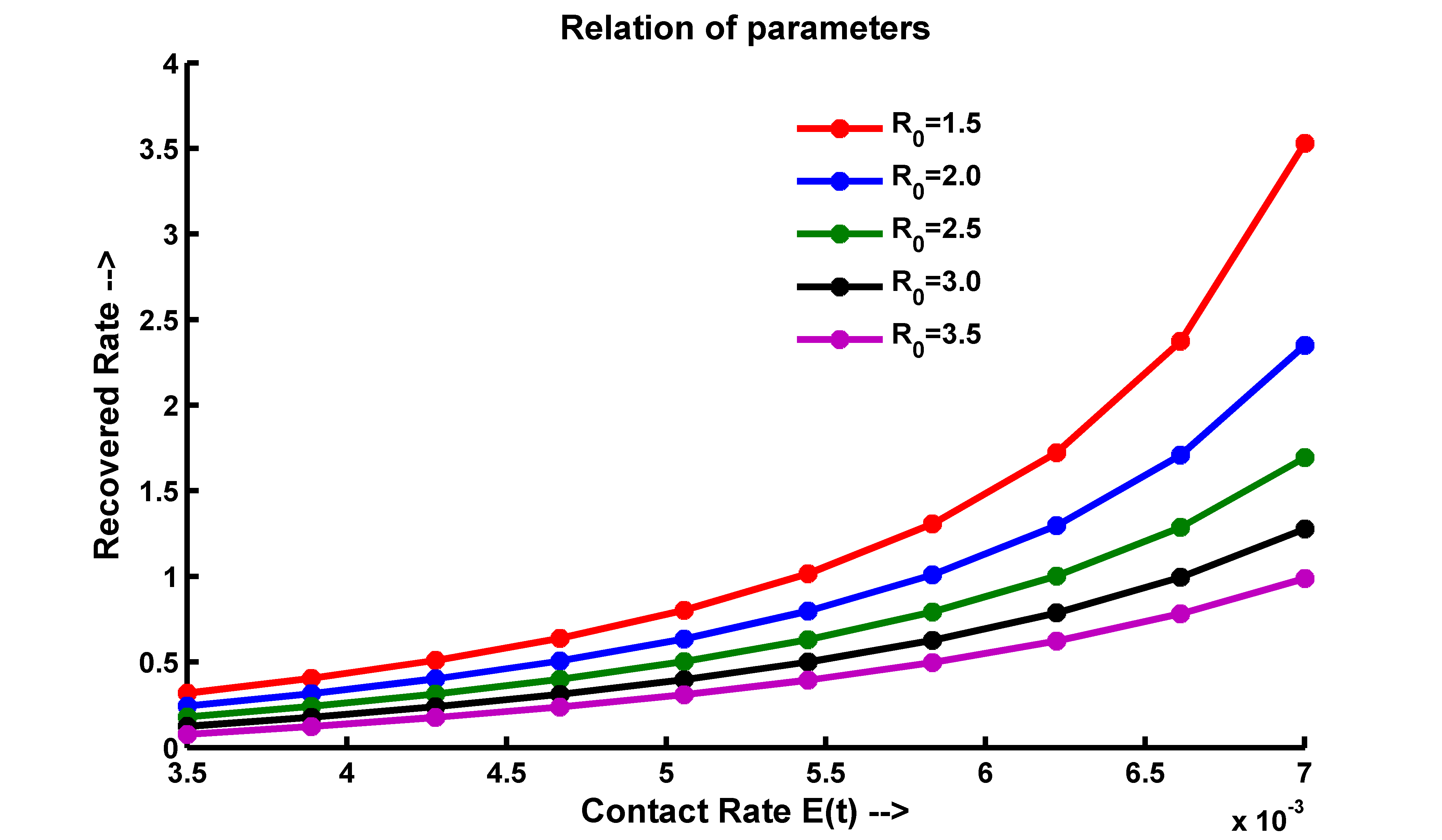}}
	\subfloat[]{\includegraphics[width=2.7 in]{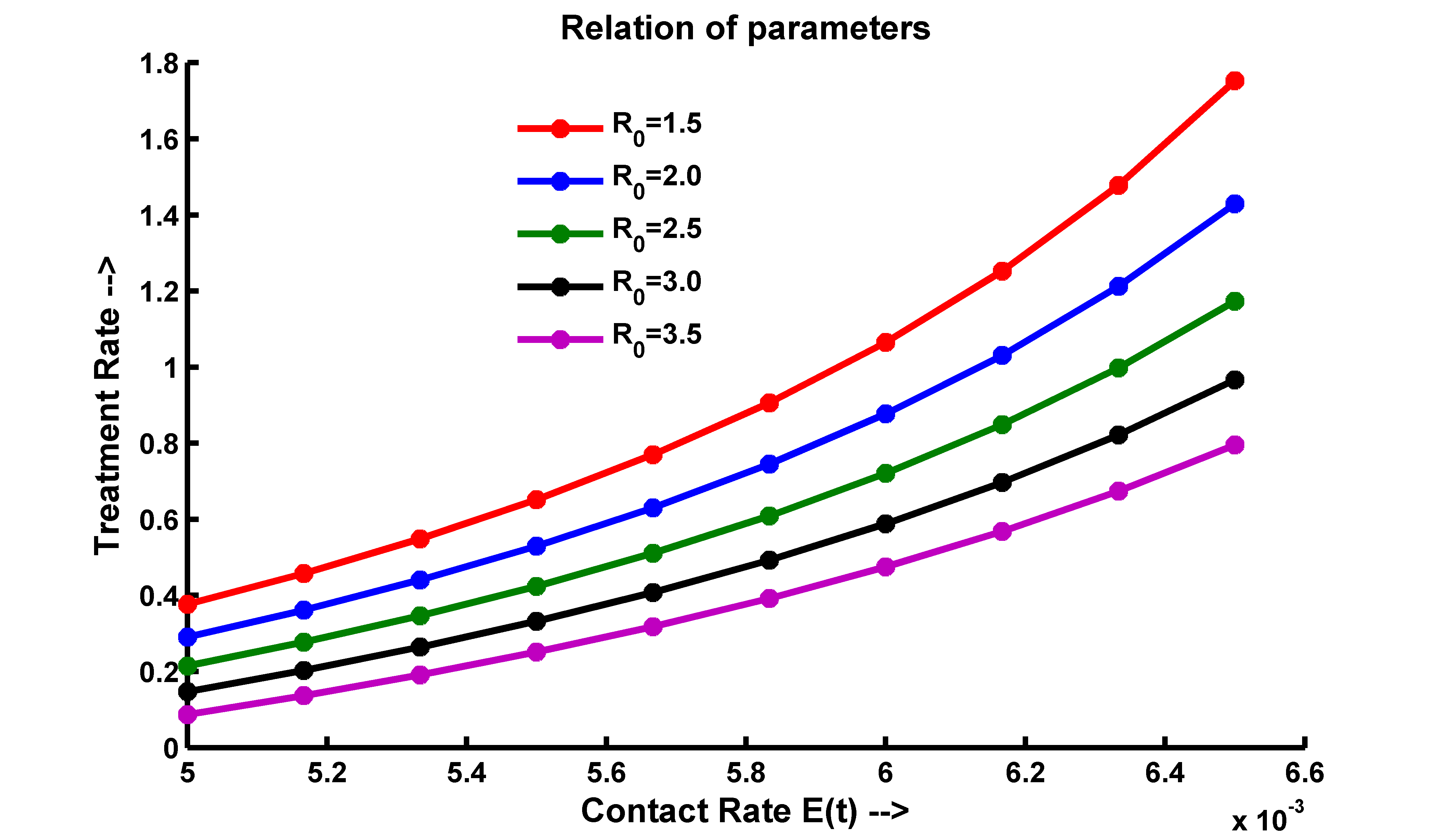}}
	\caption{For fixed range of $\mathcal{R}_0$ dynamics of relative influence of (a) $\beta_1$ and $\alpha$ (b) $\beta_1$ and $\beta_2$ (c) $\beta_1$ and $\gamma$ (d) $\beta_1$ and $\gamma_1$, where mean values of all parameters are taken from Table \ref{table-param-values-4-countries}.}
	\label{relative-influence-parameter}
\end{figure}

Figure \ref{relative-influence-parameter} reveals, some potential applications and interpretations of a contour plot in the model scenario. Figure \ref{relative-influence-parameter}(a) shows the relation of parameters $\beta_1$ and $\alpha$ for $\mathcal{R}_0\in [1.5,3.5]$. We observe that, when value of $\mathcal{R}_0$ lies in fixed range, the parameters are proportionally related. When the value of $\beta_1$ progress from 0.0025 to 0.0075, the infection rate rises rapidly from 0.35 to 0.65. The graph reveals as straight lines.

Figure \ref{relative-influence-parameter}(b) depicts that, for fixed values of $\mathcal{R}_0 \in [1.5,3.5]$, the contact rates $\beta_1$ and $\beta_2$ are reversely related i.e., negatively correlated. The contour graphs behaves like straight lines which indicates the increasing of $\beta_1$ from 0.0025 to 0.0075 results into rapid decrease of $\beta_2$ from 0.0085 and then approaches to the zero level. Thus, small changes in parameter $\beta_1$ results in significant change in $\beta_2$ to control the outbreak of influenza.

Figure \ref{relative-influence-parameter}(c) indicates, the relation of parameter $\beta_1$ and $\gamma_1$ by several contour lines. This describes that, for $\mathcal{R}_0\in[1.5,3.5]$, the parameters $\beta_1$ and $\gamma_1$ is positively correlated. The rapid growth of $\beta_1$ results into drastic impact on $\gamma$. The more contact rate increases from 0.0025 to 0.0075 results into rise of recovery rate $\gamma$ from 0.35 to 0.75. Graphs indicate that, for $\mathcal{R}_0$ values 1.2 and 2.5 the lines becomes quadratic. At the initial stage, the other three lines indicates linear relationship to control $\mathcal{R}_0$.
\noindent
Figure \ref{relative-influence-parameter}(d) depicts the influencing situation of parameter $\beta_1$ and $\gamma_1$ for fixed rage of $\mathcal{R}_0\in[1.5,3.5].$ The graphs reveals as hyperbolic which indicates a strong correlation of $\beta_1$ and $\gamma_1$ to control $\mathcal{R}_0.$ When contact rate $\beta_1$ progress from 0.0025 to 0.0075, this results a positive influence on $\gamma_1$. That means corresponding treatment rates needs to be increased to control the disease. From Figure \ref{relative-influence-parameter}(d), we get a scenario that, $\gamma_1$ rate need to increase rapidly from 0.02 to 0.55 to mitigate the disease burden.\\
Epidemiological meaning is that, hen the contact rates $\beta_1$ and $\beta_2$ and infection rate $\alpha$ increases, then corresponding recovery rates $\gamma$ and treatment rate $\gamma_1$ needed to be increased using vaccination, quarantine and proper isolation. Thus $\alpha$, $\beta_1$ and $\beta_2$ has drastic impact in the other parameters of model \eqref{new_model}.
\section{Sensitivity Analysis}\label{Section-Sensitivity Analysis}
Sensitivity analysis is used in epidemiology to determine how changes in various factors impact the transmission and consequences of infectious diseases, specifically with regard to the fundamental reproduction number $(\mathcal{R}_0)$. In this study, this approach aids in determining the critical parameters that significantly impact $\mathcal{R}_0$. Researchers can evaluate these characteristics' effects on the spread of disease by adjusting them within a reasonable range. This information can direct public health initiatives and is useful for determining which elements are most important in the disease's transmission. Epidemiological models frequently depend on input parameters like transmission rates, contact patterns, or illness duration that may be subject to uncertainty. Sensitivity analysis allows us to measure the impact of parameter uncertainty on $\mathcal{R}_0$, the threshold quantity that behaves like a crucial epidemiological parameter. We are able to gain a better understanding of the possible variety in disease outcomes and associated uncertainties by looking at the range of $\mathcal{R}_0$ values produced by altering unknown parameters. Sensitivity analysis is a useful tool for assessing how well various intervention techniques reduce $\mathcal{R}_0$. By examining how changes in parameters, such as contact rates, vaccination coverage, or treatment availability, affect $\mathcal{R}_0$, from which we can assess the potential impact of interventions and guide decision-making processes \cite{Bifurcation of R0-7, Bifurcation of R0-8}. This information is crucial for designing effective control measures and optimizing resource allocation. By analyzing this, we can identify scenarios where $\mathcal{R}_0$ surpasses a critical threshold, leading to potential outbreaks or epidemic control. This information helps us in this study to make informed decisions regarding public health measures, resource allocation, and response planning. The goal is to determine which parameters are worth pursuing in the field in order to develop a Influenza transmission model. More importantly, sensitivity analysis methods can be categorized as:\\
(1) Local Sensitivity Analysis, and (2) Global Sensitivity Analysis.
\subsection{Local Sensitivity Analysis}
To ascertain which model parameter has the most bearing on the dynamics of the model, local sensitivity analysis is used to the model's parameters. Using this approach, each input parameter is changed individually, bringing the others back to their nominal levels. The sensitivity indices of the basic reproduction number $(\mathcal{R}_0)$ to the model's \eqref{new_model} parameters have been determined in this section, to ascertain which parameters affect $(\mathcal{R}_0)$ and hence, the propagation of disease—highly, averagely, and lowly.We employ the strategy outlined in \cite{Stability Bound-25, Stability Bound-26}. 
For parameter $P_i$, the normalized sensitivity indices of $\mathcal{R}_0$ are as follows:
\begin{align*}
	I_{P_i}^{\mathcal{R}_0}=\frac{\partial \mathcal{R}_0}{\partial P_i}\frac{P_i}{\mathcal{R}_0}.
\end{align*}
The $\mathcal{R}_0$ sensitivity indices are computed using the parameter mean values shown in Table \ref{tableparameter} and the approximate parameter values shown in Table \ref{table-param-values-4-countries}. We approximated the parameter values by examining data from the CDC \cite{CDC} and WHO \cite{WHO} for the following nations: Mexico, Italy, Colombia, and South Africa. The results are presented in Table \ref{table-sensitivity index-4-countries}, where the data from the examples of Mexico, Italy,  and South Africa are shown in the second, third, fourth, and fifth columns, respectively. Generally, the value of $\mathcal{R}_0$ grows (or reduces) depending on which of the positive-sign parameters increases (or decreases) while the other parameters remain constant. For instance, raising $\beta_1$ by 10\% increases $0.6247\times 10\%$ of $\mathcal{R}_0$. However, rising $\gamma$ by 10\% decreases $0.2247\times 10\%$ of $\mathcal{R}_0$.\\
Table \ref{table-sensitivity index-4-countries} demonstrates that disease transmission probability rates from $E(t)$ to $S(t)$ are the most significant and critical characteristics i.e., $\beta_1$, disease transmission probability rates from $I(t)$ to $S(t)$ i.e. $\beta_2,$ recovery rate from $I(t)$ i.e., $\gamma$, infection rate from $E(t)$ to $I(t)$ i.e., $\alpha$ and human recruitment rate $\Lambda.$ On the other hand, mortality rate $\mu$ , treatment rate $\gamma_1$ and disease induced death rate $\delta$ have lower impact on $\mathcal{R}_0.$
\begin{table}[H]
	\begin{center}
		\caption{Sensitivity indices of $\mathcal{R}_0$ to model \eqref{new_model} parameter values, assessed at the parameter values shown in Table \ref{table-param-values-4-countries}.}
		\scriptsize
		\label{table-sensitivity index-4-countries}
		\begin{tabular}{|c|c|c|c|c|}
			\hline\noalign{\smallskip}
			\multirow{2}{*}{\textbf{Parameters}}&\textbf{Mexico}& \textbf{Italy}  &\textbf{South Africa} \\
			\cline{2-4}
			&\textbf{Sensitivity Index}&\textbf{Sensitivity Index}&\textbf{Sensitivity Index}\\
			\noalign{\smallskip}\hline\noalign{\smallskip}
			$\alpha$&-0.5622&-0.4322&-0.4472\\
			$\beta_1$&0.6247&0.6175&0.6735\\
			$\beta_2$&0.3752&0.5152&0.5552\\
			$\gamma$&-0.1951&-0.2247&-0.2747\\
			$\gamma_1$&-0.0751&-0.0635&-0.0535\\
			$\lambda$&0.0124&0.0214&0.0244\\
			$\mu$&-0.0774&-0.0624&-0.0754\\
			$\delta$&-0.09&-0.06&-0.09\\
			$\Lambda$&0.95&0.83&0.87\\
			\noalign{\smallskip}\hline
		\end{tabular}
	\end{center}
\end{table}

\subsection{Analysis of the Local Sensitivity Result}
\begin{figure}[H]
	\centering  
	\subfloat[]{\includegraphics[width=2.9 in]{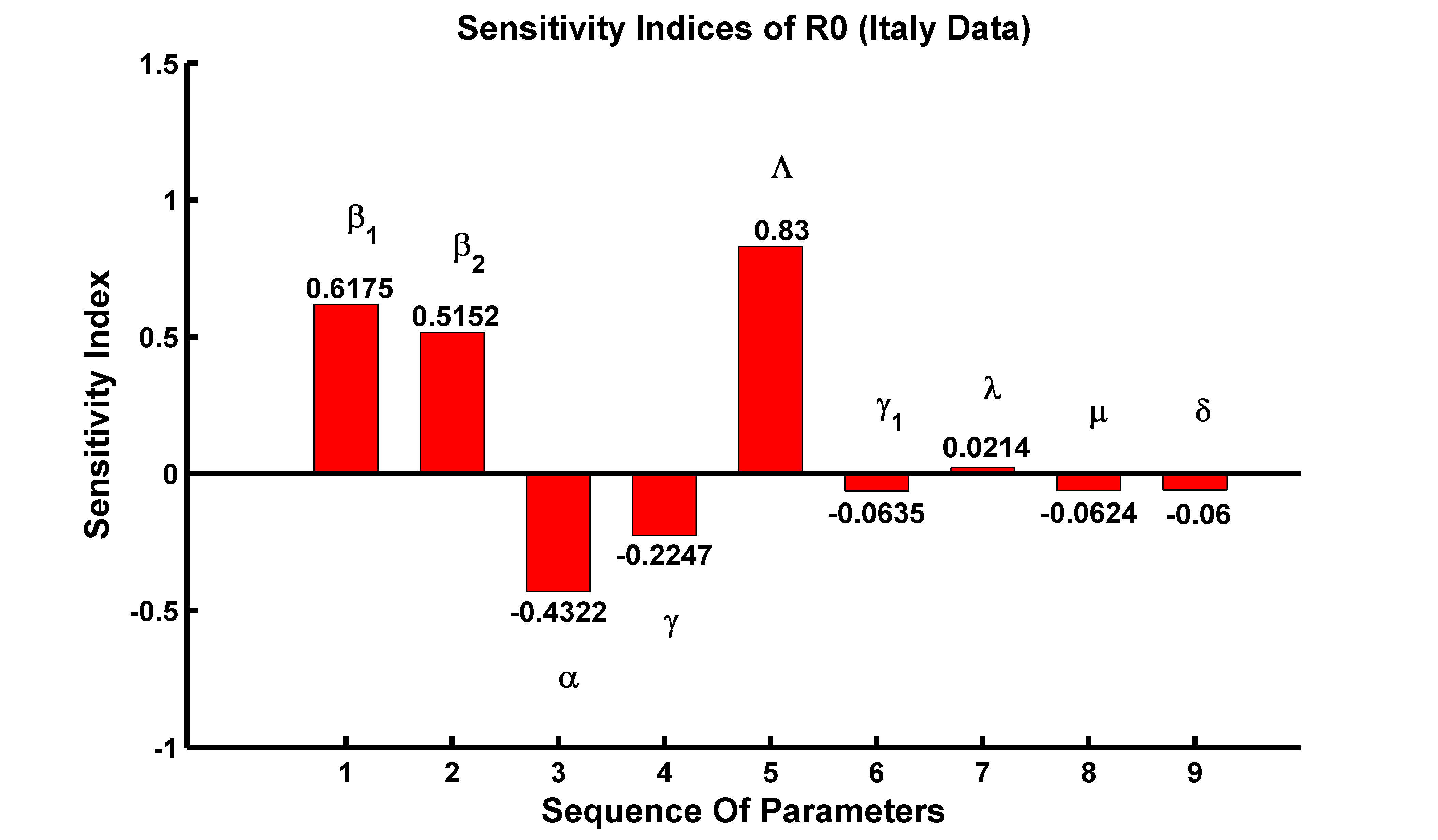}}
	\subfloat[]{\includegraphics[width=2.9 in]{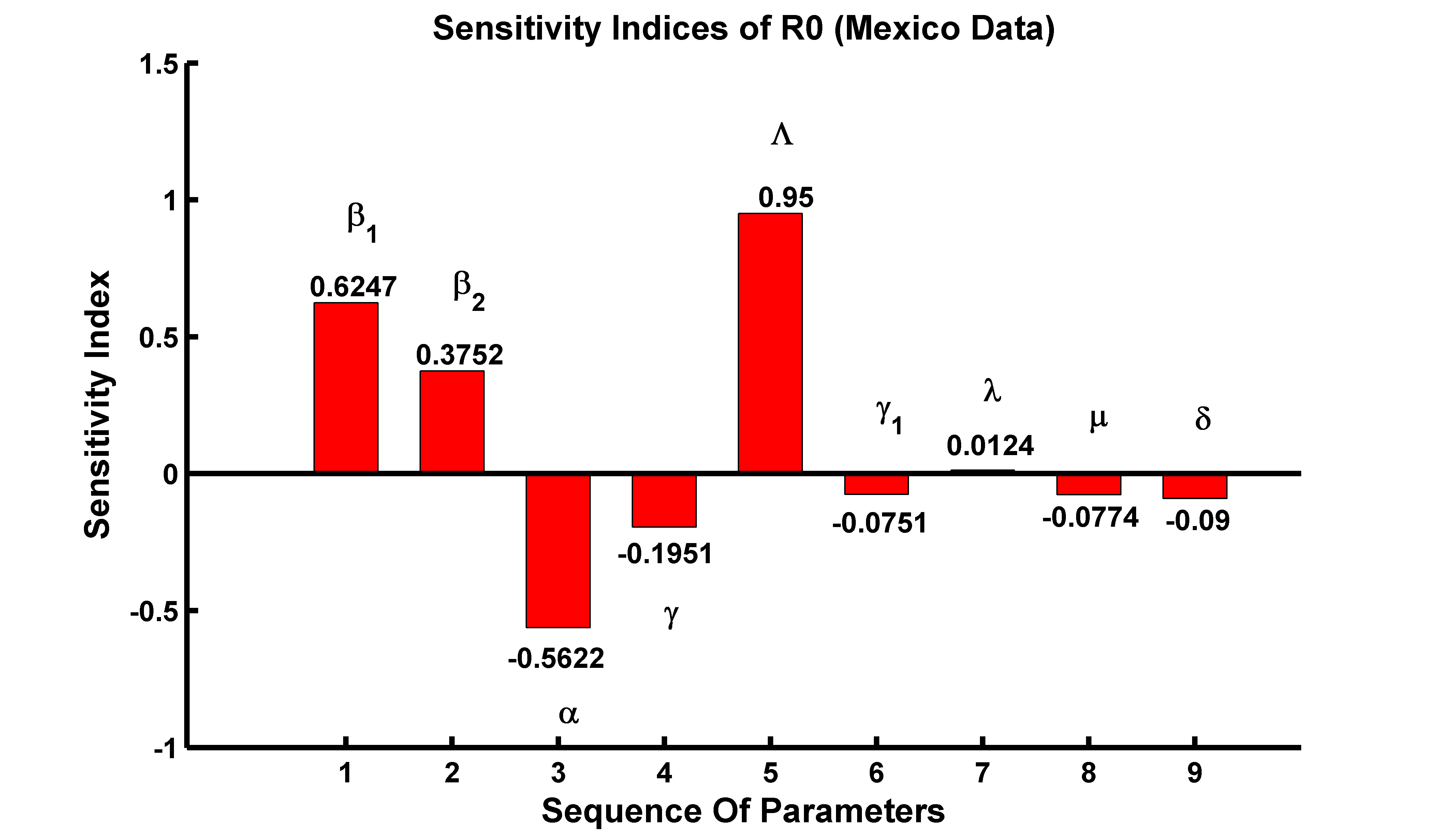}}\\
	\subfloat[]{\includegraphics[width=2.9 in]{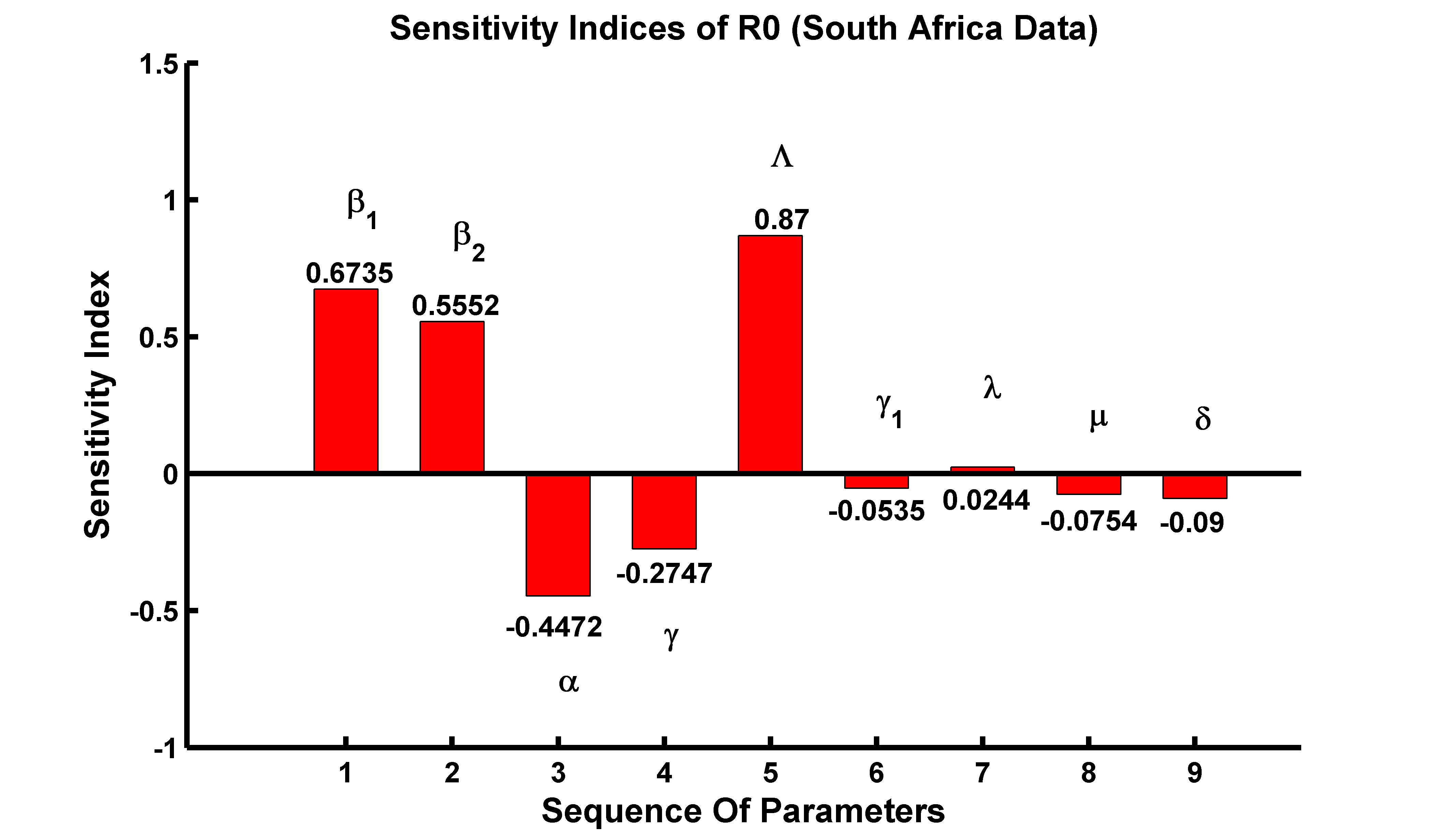}}
	\caption{Normalized sensitivity indices of $\mathcal{R}_0$ with respect to parameters for (a)   Italy data (b)  Mexico data and (c) South Africa data}.
	\label{local-sensitivity-analysis}
\end{figure}

By calculating sensitivity indices for all parameters of the three countries, presented in Table \ref{table-sensitivity index-4-countries} indicates that,
\begin{enumerate}
	\item For a particular values of parameters, from Figure \ref{local-sensitivity-analysis}, we see that the parameters having positive lower impact on $\mathcal{R}_0$ is $\lambda$ and the parameters having negative lower impact on $\mathcal{R}_0$ are $\mu$, $\delta$ and $\gamma_1$. This, reveals that by isolating the infected individuals and taking proper treatment can give significant result to mitigate the outbreak.
	\item The result from Figure \ref{local-sensitivity-analysis},  implies small changes in any of the parameters $\beta_1$,$\beta_2$,$\alpha$ and $\gamma$ contributes significantly to $\mathcal{R}_0$.
	\item When $\gamma$ and $\gamma_1$ are involved, the basic reproduction number ($\mathcal{R}_0$) has a negative effect. That means $\mathcal{R}_0$ decreases slightly when $\gamma$ and $\gamma_1$ rate increases. So mass awareness to the susceptible individuals, proper treatment to the infected individuals can be effective policy. Using vaccination properly can decrease $\lambda$ and mitigate $\mathcal{R}_0$.
	\item Table shows that the recruitment rate $\Lambda$ and $\alpha$ have higher impact on $\mathcal{R}_0$ (Figure \ref{local-sensitivity-analysis}). That means, the reduction quantity from exposed class can provide good result to control $\mathcal{R}_0$.
	\item The most crucial parameters, which have positive impact on $\mathcal{R}_0$ are $\beta_1$ and $\beta_2$. That means, $\mathcal{R}_0$ increases rapidly when $\beta_1$ and $\beta_2$ increases. Hence, proper isolation, vaccination, and quarantine strategies to the infected individuals ought to be upheld in order to lower the rate of interaction, which can stop the illness outbreak.
\end{enumerate}
\subsection{Latin Hypercube Sampling of Parameters}
\begin{figure}[H]
	\centering  
	\subfloat[]{\includegraphics[width=2.17 in]{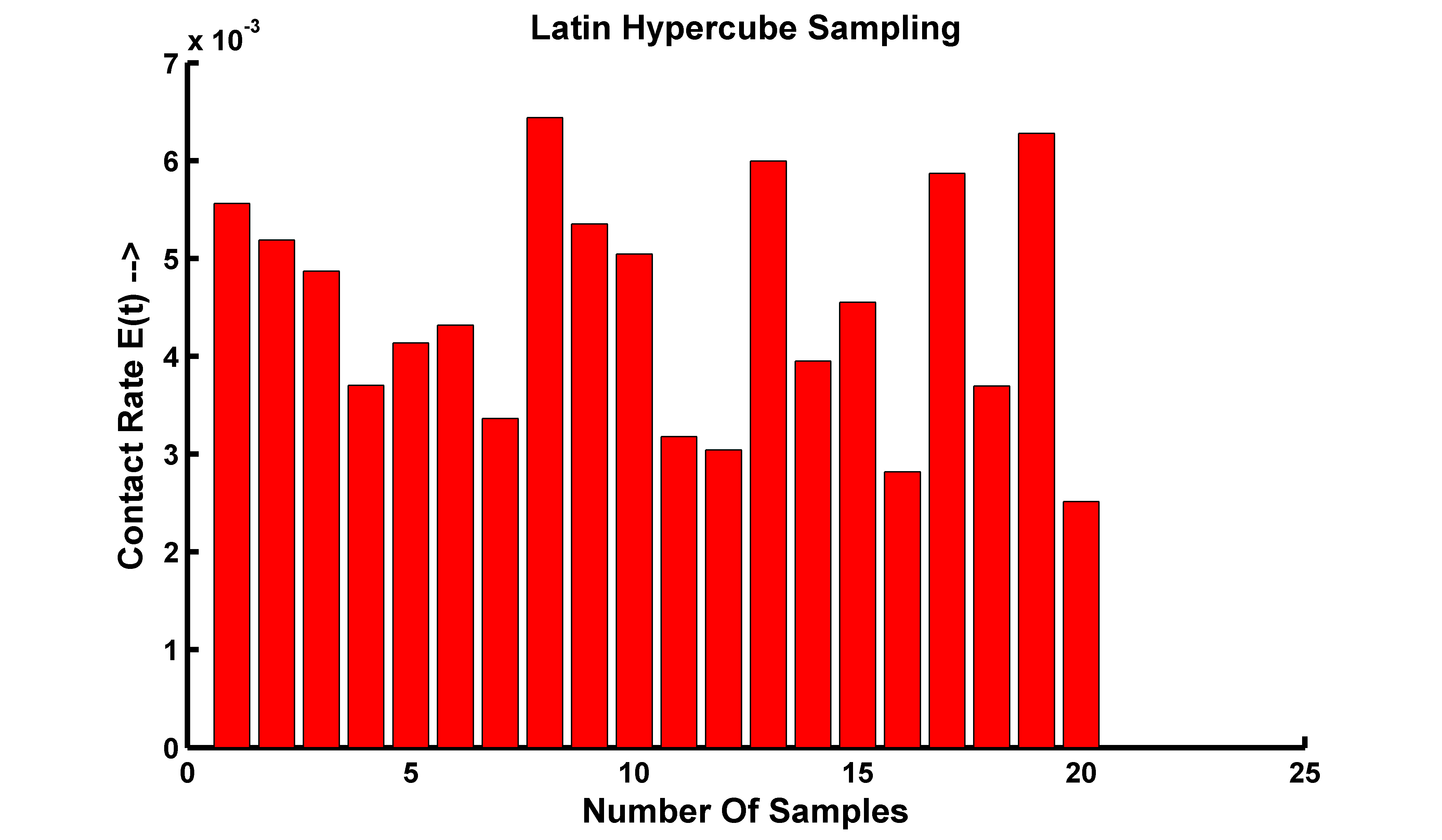}}
	\subfloat[]{\includegraphics[width=2.17 in]{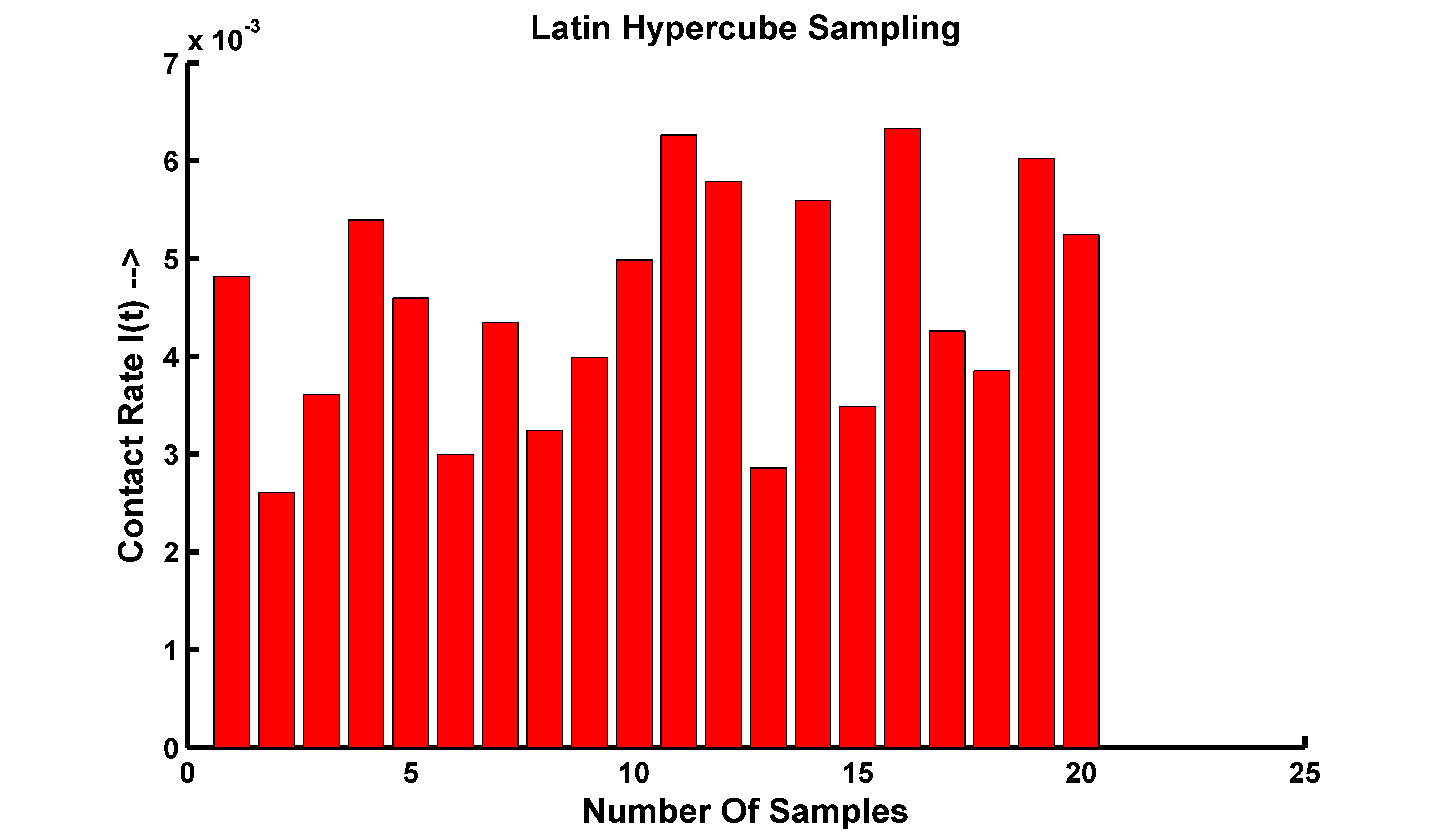}}
	\subfloat[]{\includegraphics[width=2.17 in]{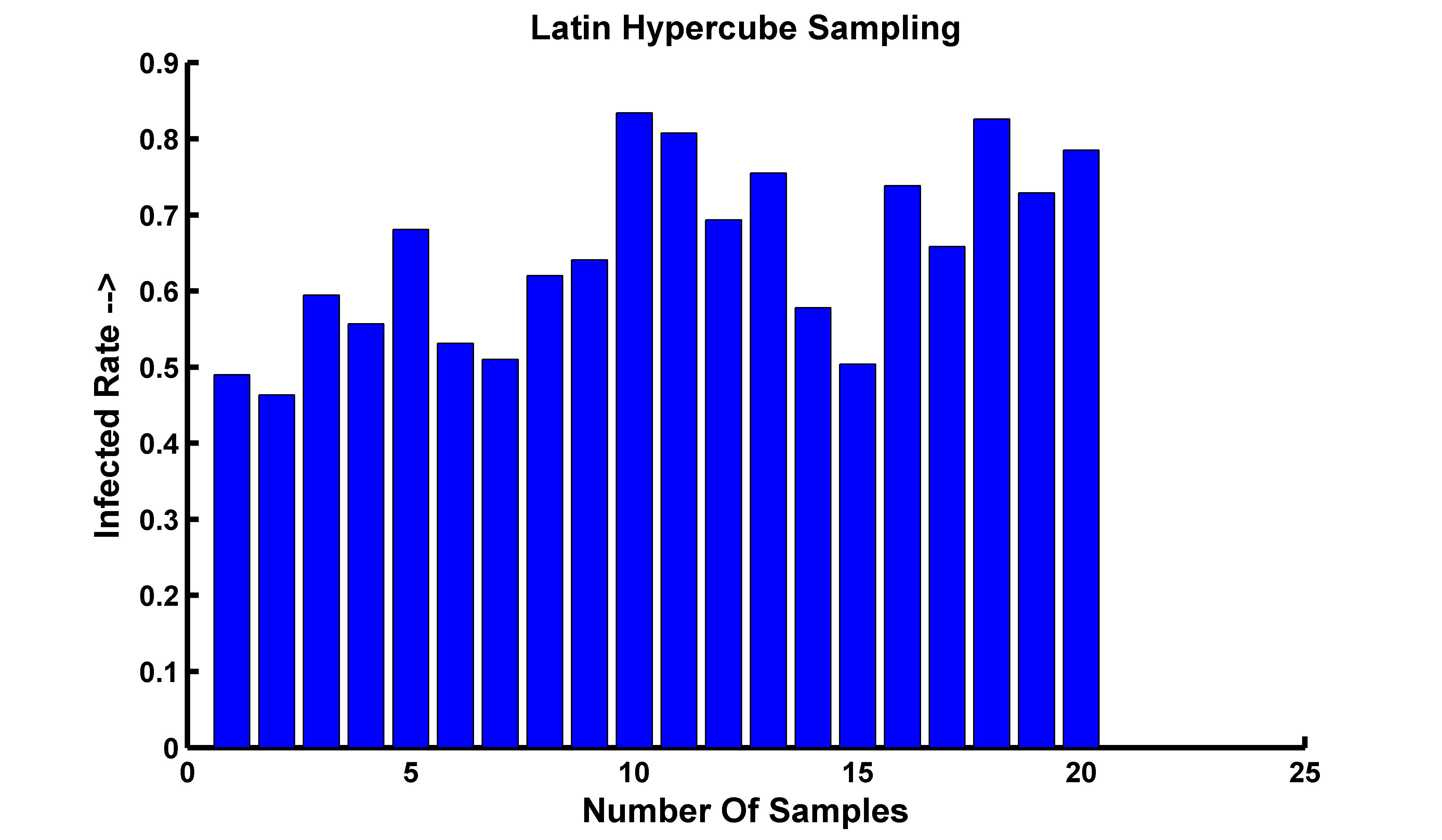}}\\
	\subfloat[]{\includegraphics[width=2.17 in]{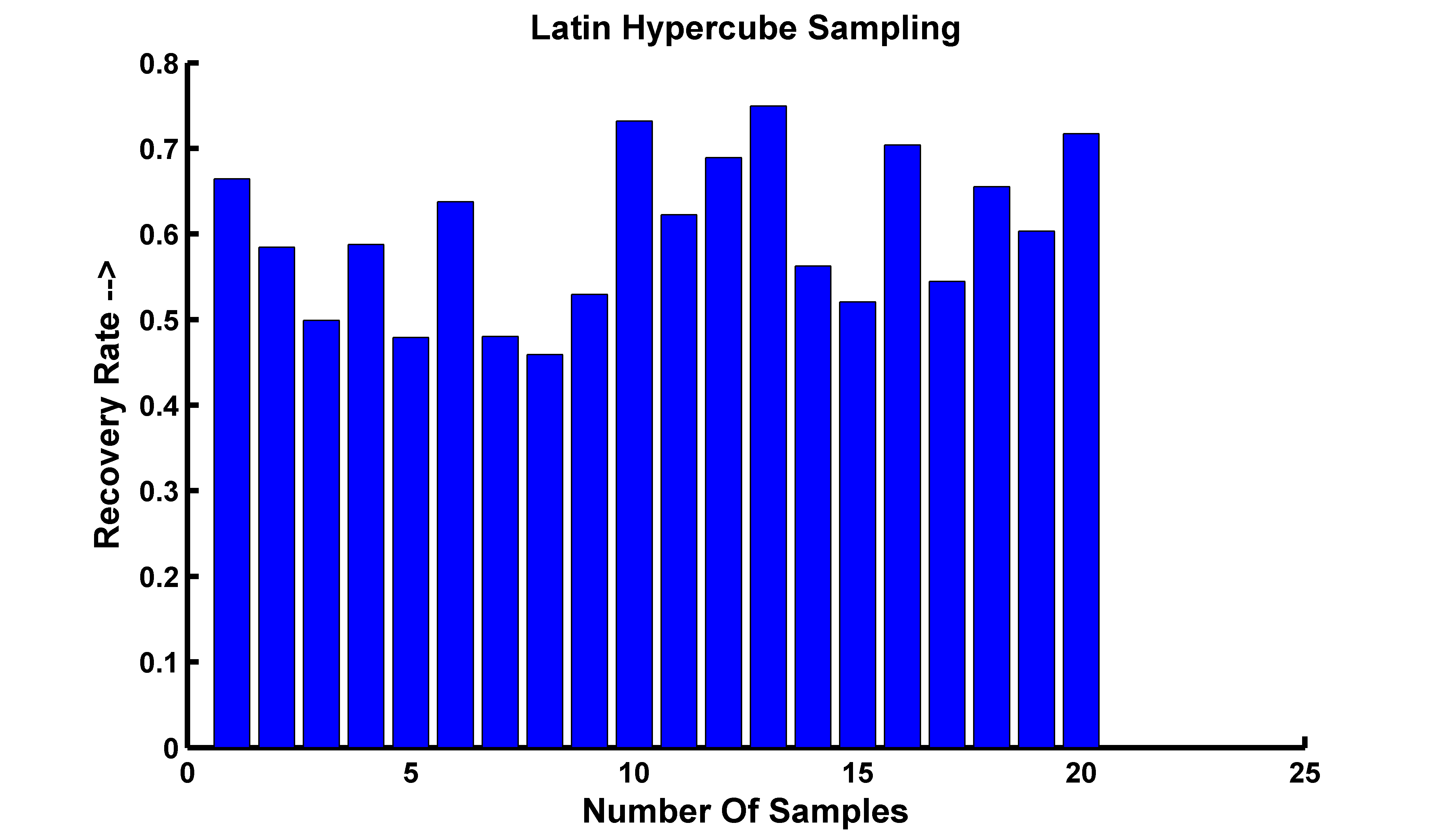}}
	\subfloat[]{\includegraphics[width=2.17 in]{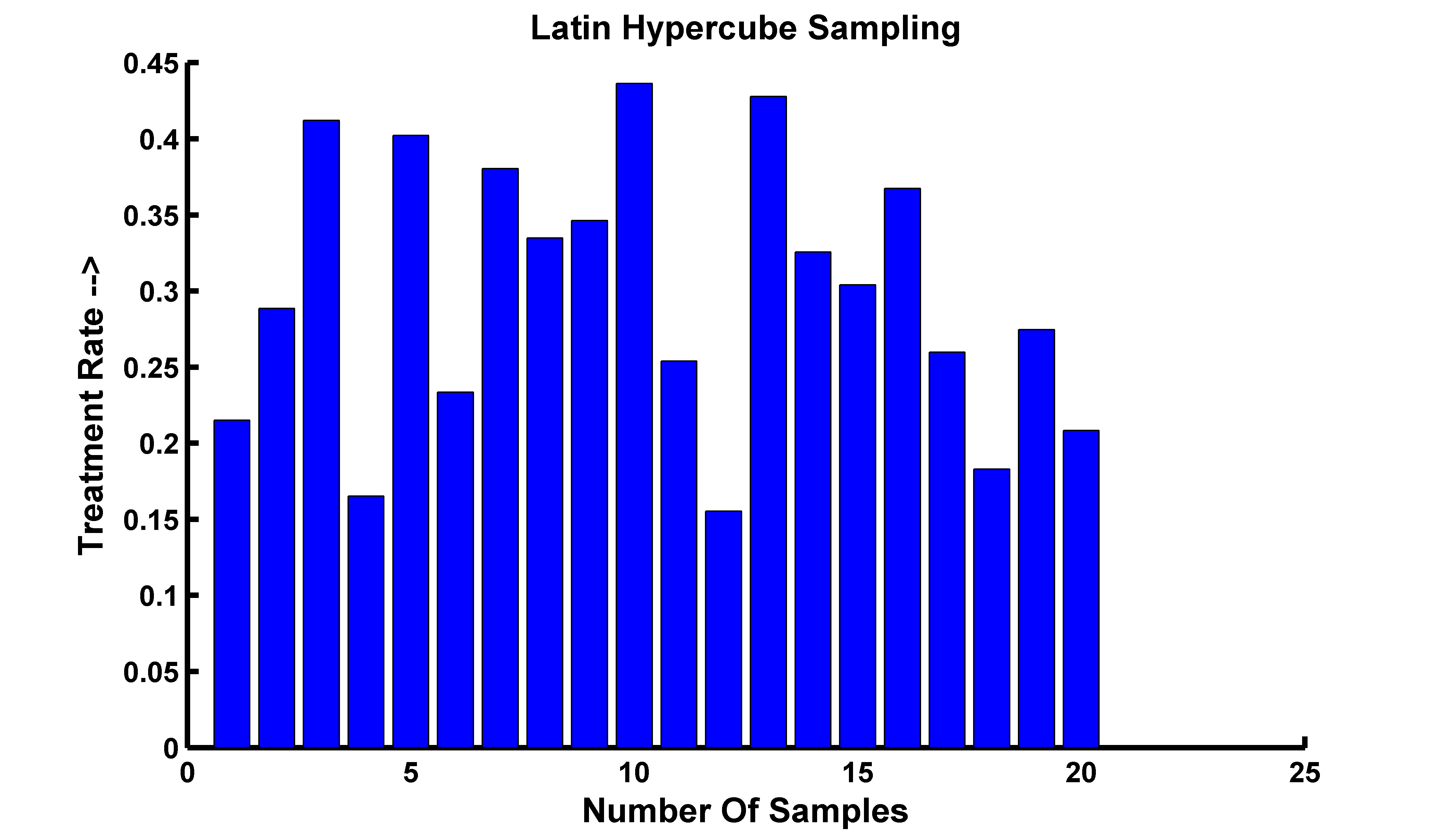}}
	\subfloat[]{\includegraphics[width=2.17 in]{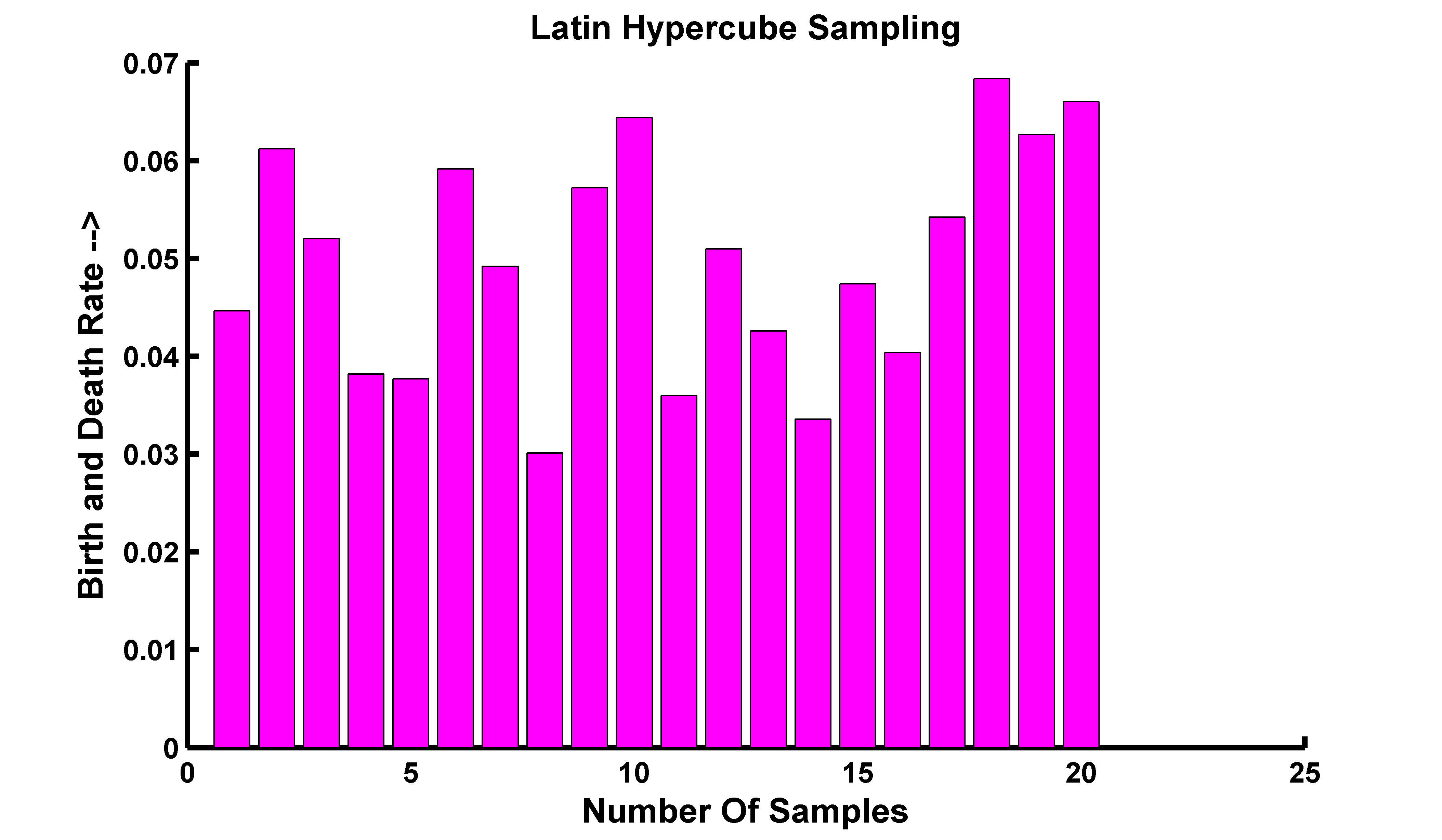}}\\
	\subfloat[]{\includegraphics[width=2.17 in]{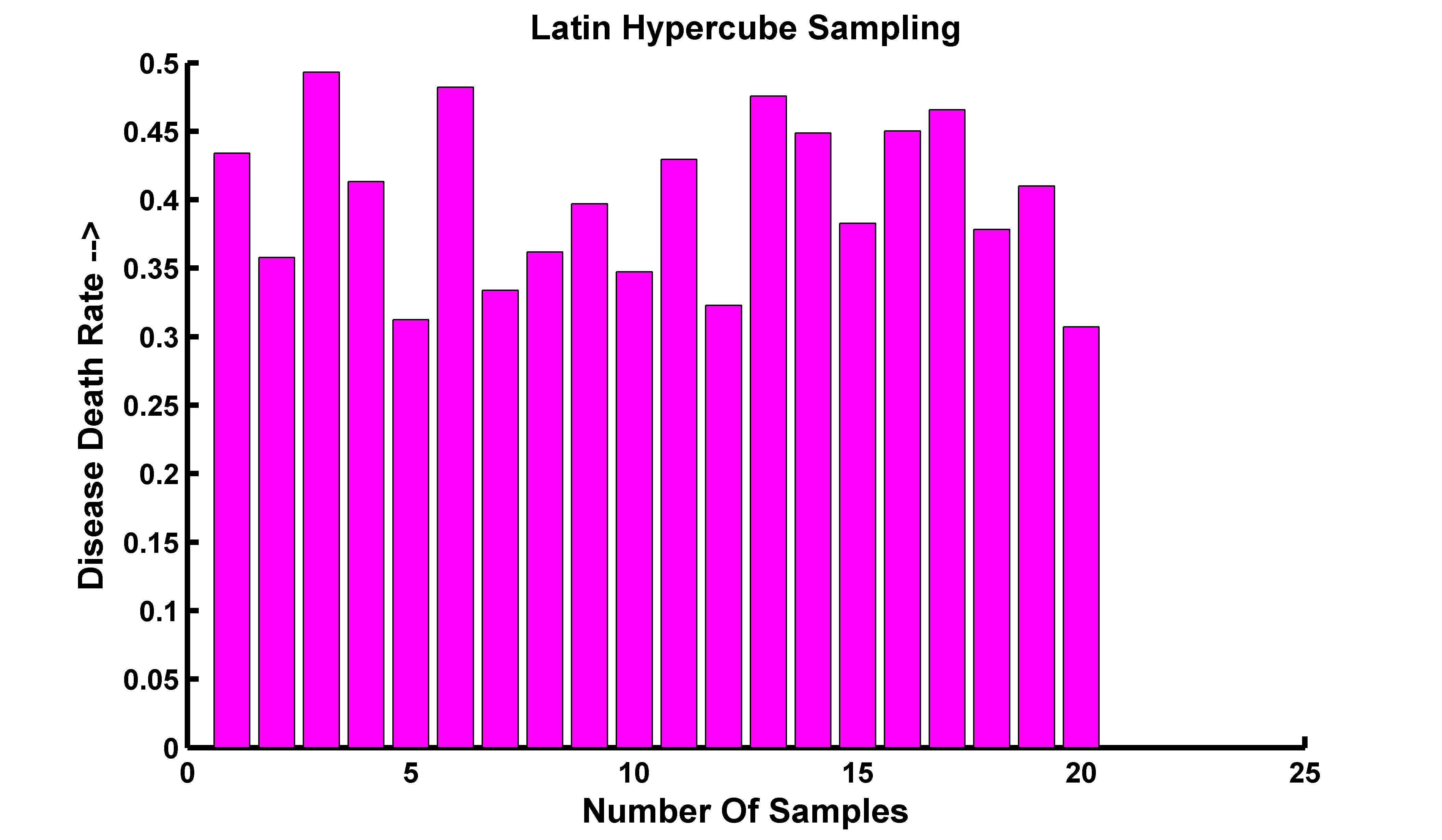}}
	\subfloat[]{\includegraphics[width=2.17 in]{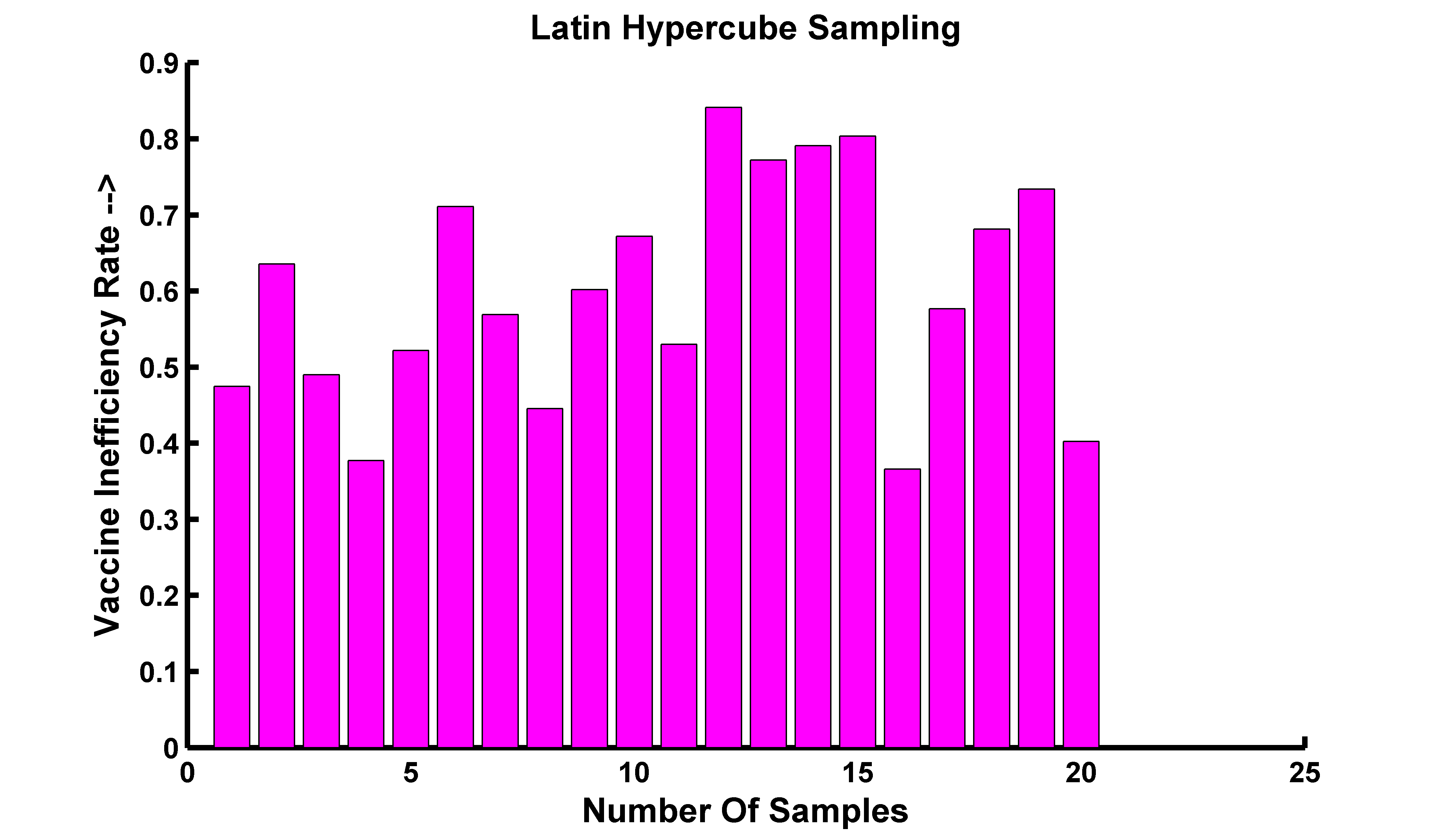}}
	\subfloat[]{\includegraphics[width=2.17 in]{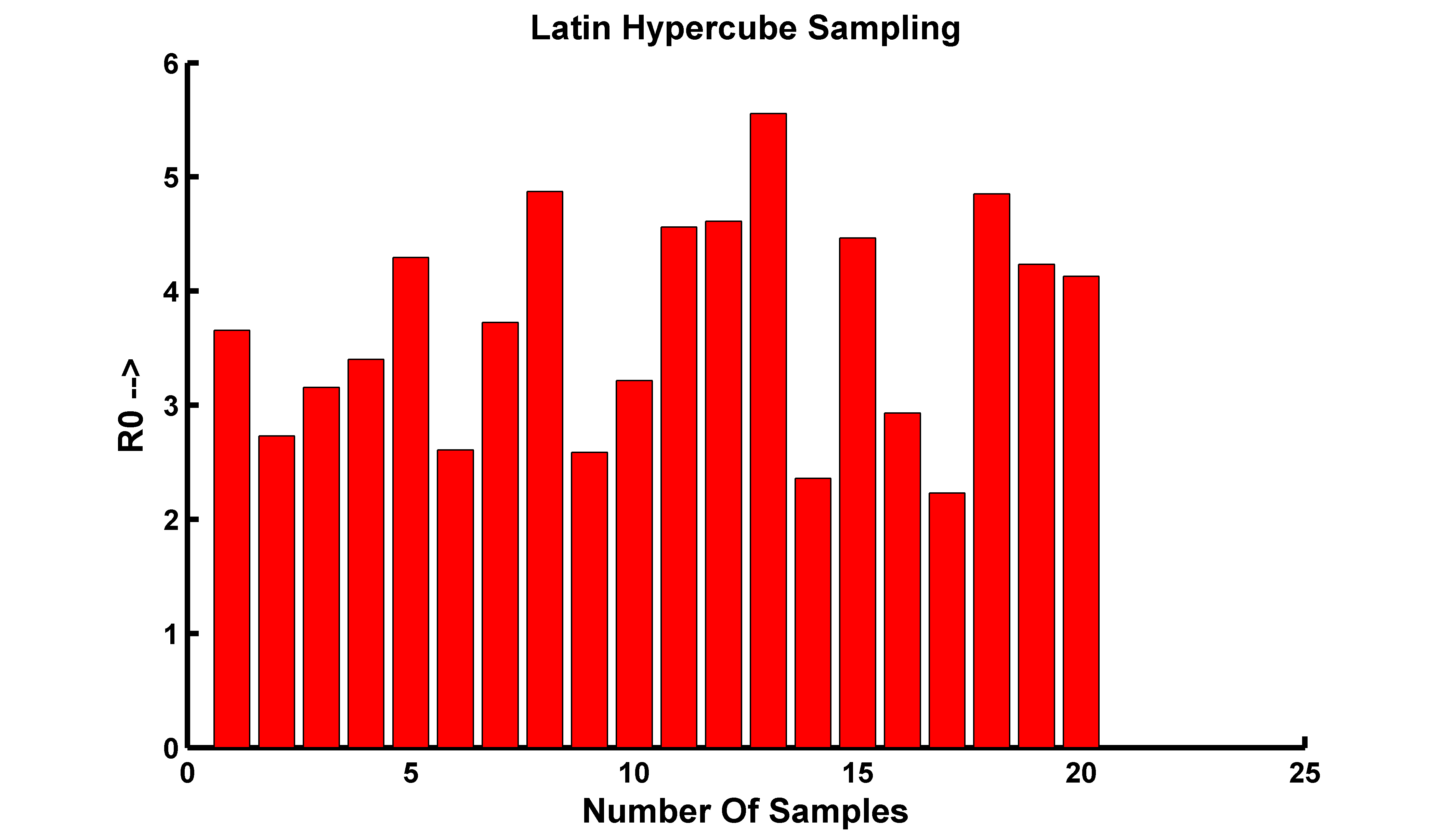}}
	\caption{Latin Hypercube Sampling graphs of parameters (a) $\beta_1$ (b) $\beta_2$ (c) $\alpha$ (d) $\gamma$ (e) $\gamma_1$ (f) $\mu$ (g) $\delta$ (h) $\lambda$ and (i) corresponding sampling of $\mathcal{R}_0$, where in the range, parameter mean values are taken from Table \ref{table-param-values-4-countries}.}
	\label{latine-hypercube-sampling}
\end{figure}

A popular method in epidemiology for sampling parameters within predetermined ranges is Latin hypercube sampling (LHS). It is intended to guarantee a more thorough and effective parameter space search in comparison to conventional random sampling techniques \cite{Hopf LEAST LHS PRCC-3, Hopf LEAST LHS PRCC-5}. The LHS method divides the parameter ranges into equally sized intervals, and then randomly selects one value from each interval to form a sample. This process ensures that the entire range of each parameter is covered while reducing the likelihood of redundant or inefficient sampling. The key steps involved in Latin hypercube sampling that we have maintained are as follows:\\
Firstly,  we have determined the minimum and maximum values for each parameter that you want to sample. Then we divided each parameter range into equally sized intervals based on the desired number of samples. Then, we generated random numbers within each interval for each parameter. These random numbers will serve as the sampling points. After that, we have shuffled the sampling points within each parameter to ensure a random arrangement. Finally, we have assigned each sampling point to its corresponding parameter.\\
In order to provide a more comprehensive exploration of the parameter space, the resultant Latin hypercube sample represents a set of parameter combinations that span the whole range of each parameter. This sampling is a valuable tool in epidemiology as it helps us efficiently explore a wide range of parameter values while maintaining a representative and non-redundant sample.\\
Latin hypercube sampling (LHS) of the contact rate $\beta_{1}$ and $\beta_{2}$ within the range of [0.0025, 0.0065] are visually represented using a bar diagram in Figure \ref{latine-hypercube-sampling}(a) and \ref{latine-hypercube-sampling}(b). The bar diagram would consist of several equally spaced bars representing different sampling points within the specified range. The number of bars would depend on the desired sample size. Each bar would correspond to a specific contact rate value randomly selected from the defined interval. The bars would be arranged in a random order to ensure a representative and non-redundant sampling of the parameter space. The bar diagram provides a visual representation of the diverse range of contact rates sampled using Latin hypercube sampling, showcasing the comprehensive exploration of parameter values within the specified range.\\
Additionally, we sampled characteristics like infection, recovery, and treatment rates within predetermined ranges using Latin hypercube sampling (LHS), which was visually displayed using bar graphs. From Figure \ref{latine-hypercube-sampling}(c), for the infection rate $\alpha$, a bar diagram can be created with bars representing different randomly selected values from the range [0.35, 0.85]. The bars would be equally spaced and arranged in a random order, providing a visual representation of the sampled infection rates.\\
Similarly,from Figure \ref{latine-hypercube-sampling}(d), for the recovery rate $\gamma$, another bar diagram can be generated with bars representing values randomly sampled from the range [0.45, 0.75]. The bars would be equally spaced and arranged randomly, illustrating the range of recovery rates sampled through LHS. In  Figure \ref{latine-hypercube-sampling}(e), we have sampled the treatment rate $\gamma_1$ using LHS and represented in a bar diagram. The bars in this diagram would correspond to values randomly selected from the range [0.15, 0.45]. They would be equally spaced and arranged randomly, visually demonstrating the variety of sampled treatment rates.
\noindent
A statistical sampling method called Latin Hypercube Sampling (LHS) is used to effectively explore a multidimensional parameter space. To create representative samples within predetermined ranges for the birth/death rate $\mu$ and disease-induced death rate $\delta$ on a bar diagram, we have employed LHS. For the birth/death rate bar diagram Figure \ref{latine-hypercube-sampling}(f), the LHS approach would involve dividing the range [0.02, 0.06] into several equal intervals. The number of intervals would correspond to the desired number of samples or data points. LHS ensures that within each interval, only one sample is selected, resulting in a diverse and representative set of samples that adequately cover the entire range.\\
Similarly, for the disease-induced death rate $\delta$ bar diagram Figure \ref{latine-hypercube-sampling}(g), the range [0.03, 0.045] would be divided into intervals, and LHS would select one sample from each interval. By distributing the samples evenly throughout the range, this procedure makes sure that the data's variability is captured.\\
In the context of vaccine efficiency rate $\lambda$ and basic reproduction number $\mathcal{R}_0$ in a bar diagram, LHS can be applied to generate representative samples within the specified ranges. For the vaccine efficiency rate $\mathcal{R}_0$ bar diagram Figure \ref{latine-hypercube-sampling}(h), the range [0.3, 0.85] would be divided into several equal intervals, and LHS would select one sample from each interval.By distributing the chosen samples evenly over the range, this procedure guarantees a varied collection of representative data points that accurately reflect the fluctuations in vaccination success rates $\lambda$.\\
Similarly, for the basic reproduction number $\mathcal{R}_0$ bar diagram \ref{latine-hypercube-sampling}(i), the range [2, 5.5] would be divided into intervals, and LHS would select one sample from each interval. This approach guarantees that the samples are evenly spread throughout the range, allowing for an accurate representation of the variability in the basic reproduction number $\mathcal{R}_0$.

\subsection{Relative Bias of $\mathcal{R}_0$ by LHS of Parameters}
\begin{figure}[H]
	\centering  
	\subfloat[]{\includegraphics[width=2.7 in]{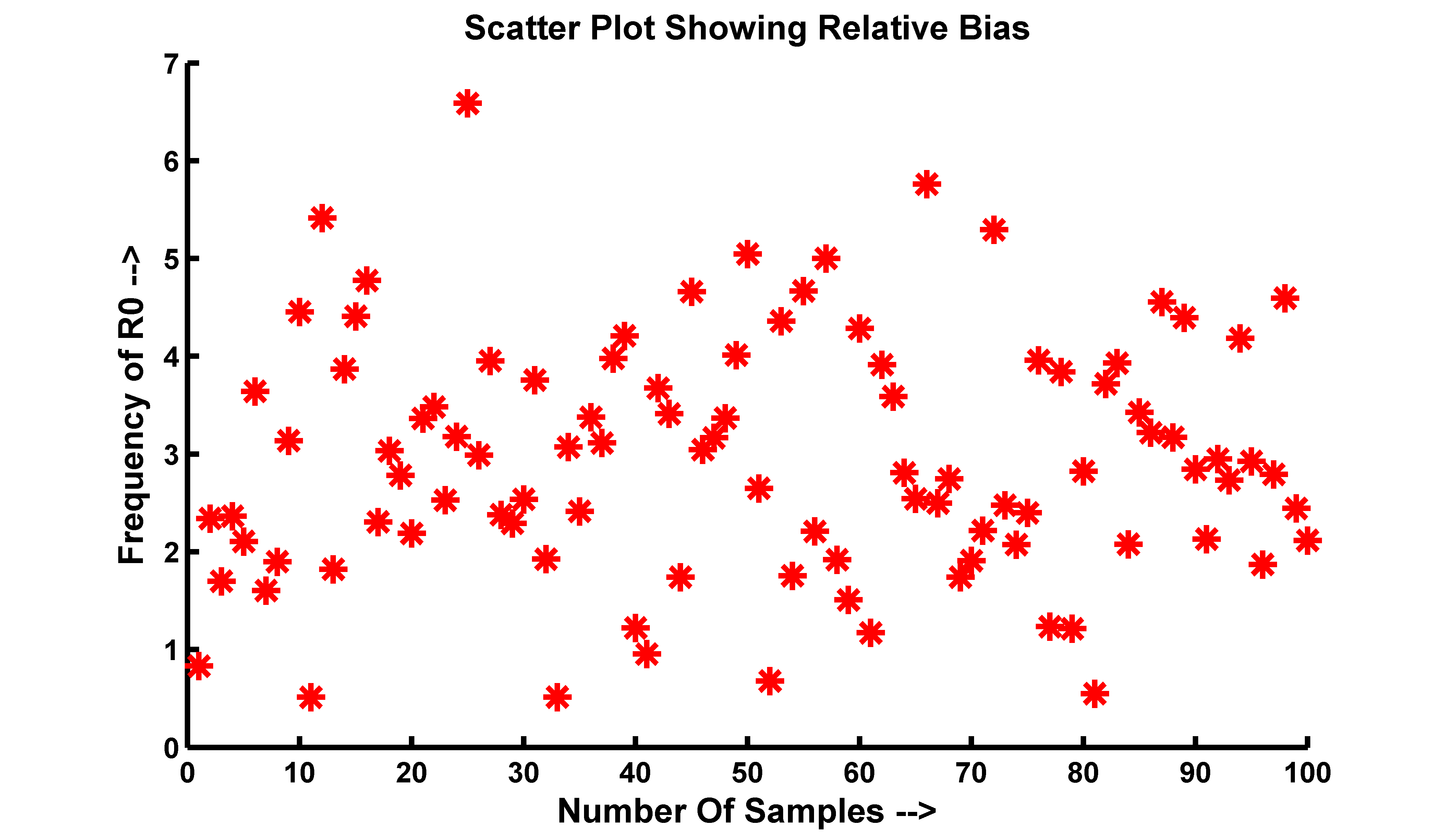}}
	\subfloat[]{\includegraphics[width=2.7 in]{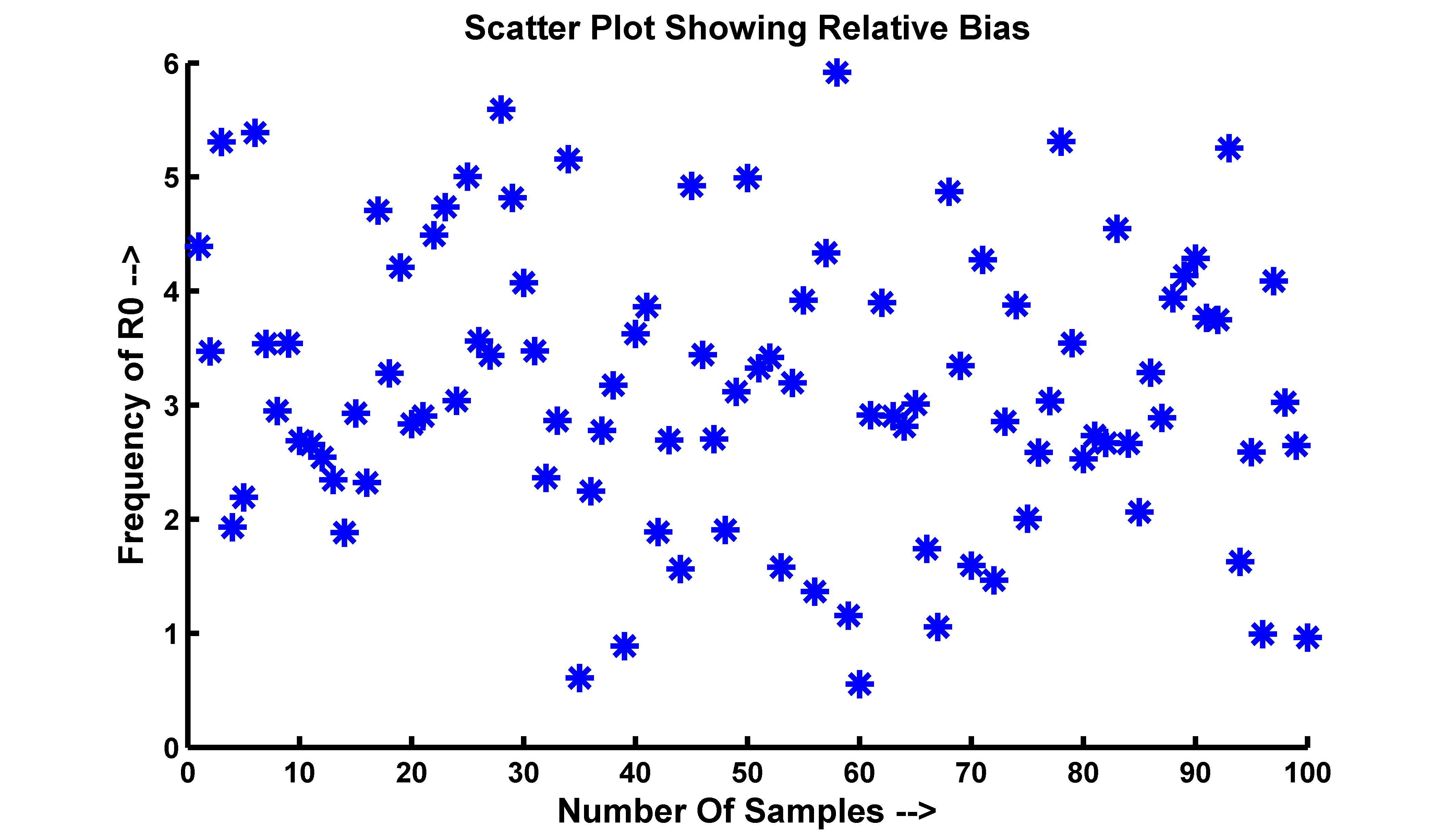}}\\
	\subfloat[]{\includegraphics[width=2.7 in]{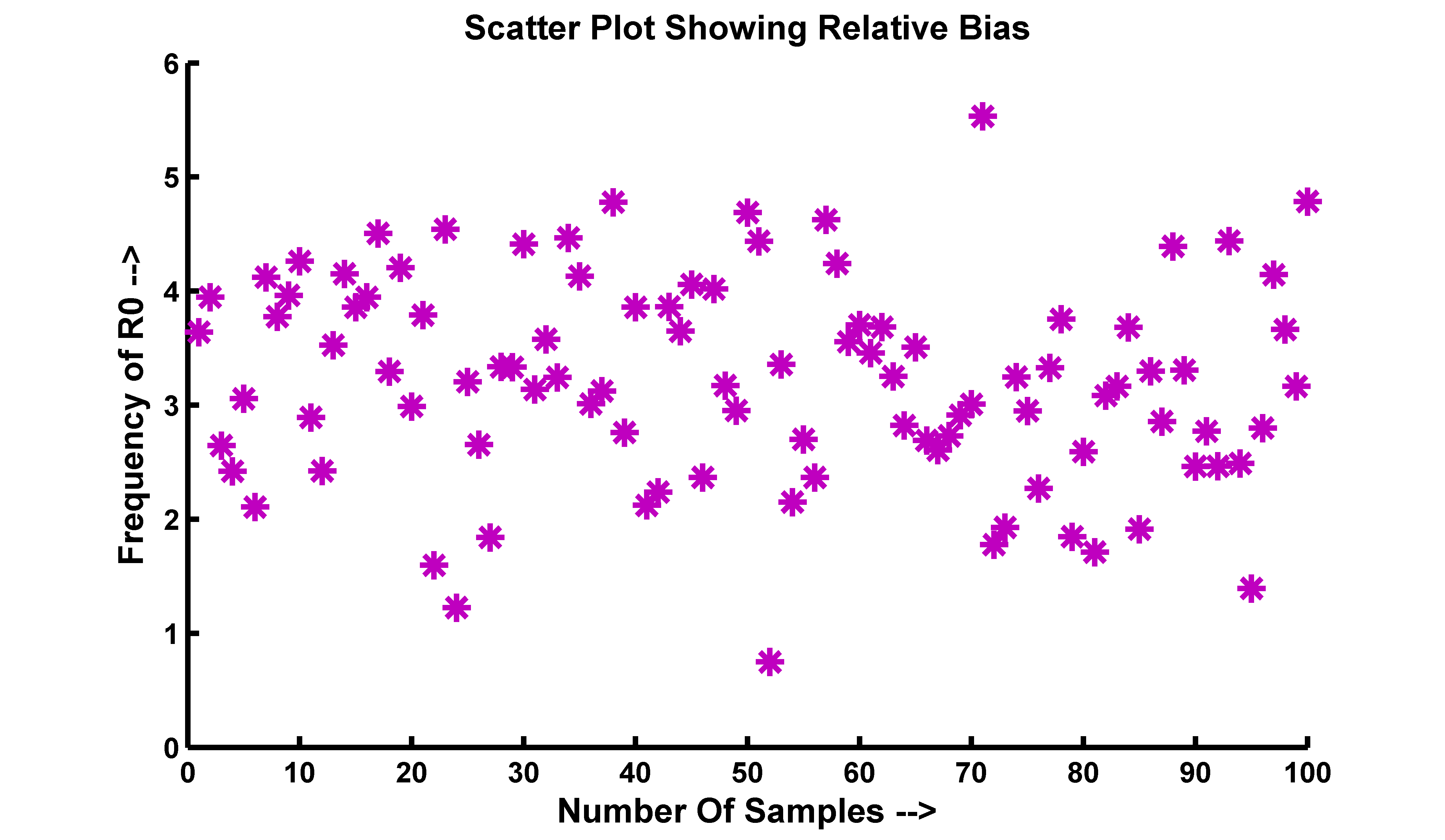}}
	\subfloat[]{\includegraphics[width=2.7 in]{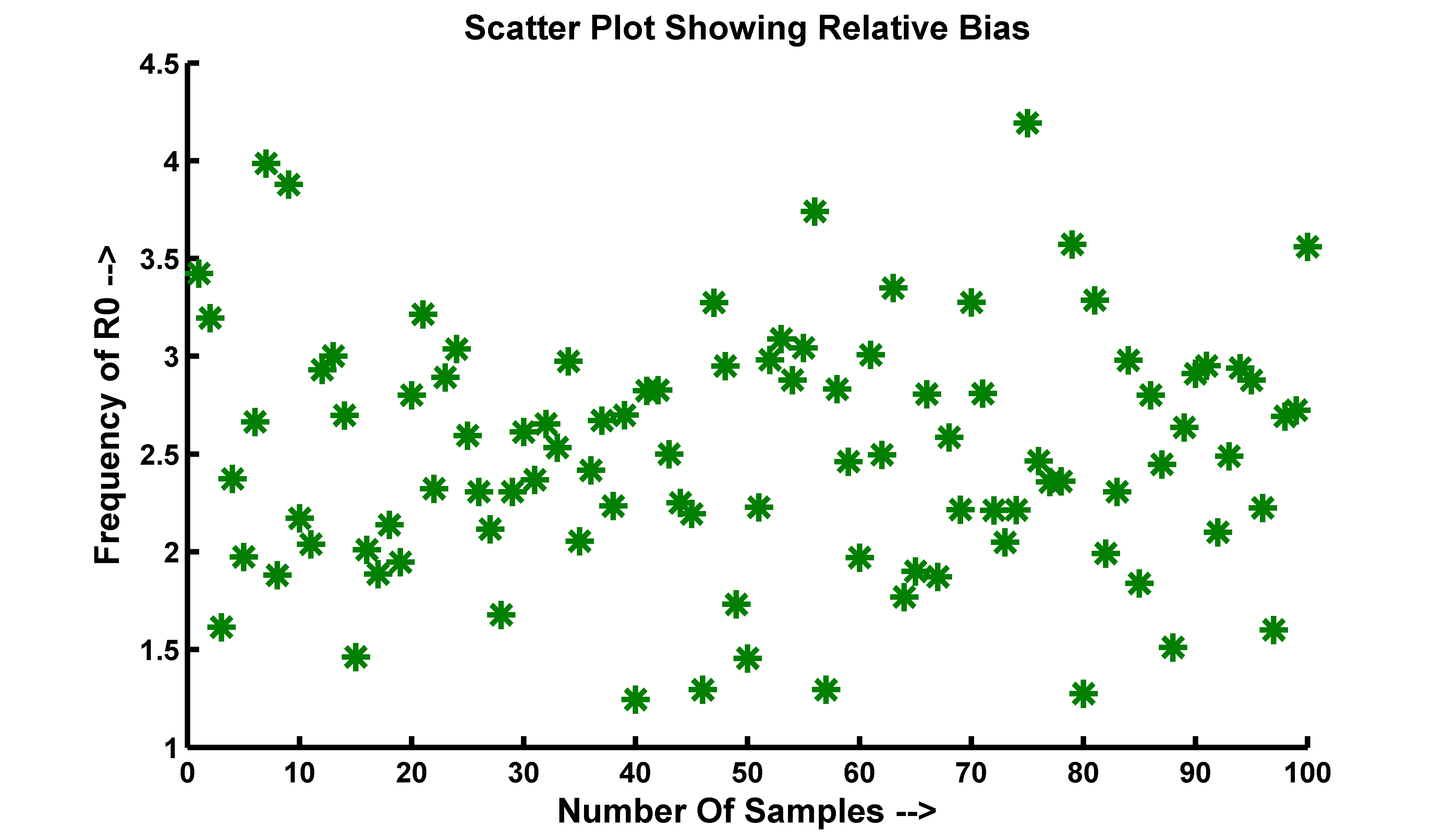}}
	\caption{Relative bias of $\mathcal{R}_0$ by LHS of parameters (a) $\beta_1$ and $\beta_2$ (b) $\alpha$ and $(1-\lambda)$ (c) $\gamma$ and $\gamma_1$ (d) $\delta$ and $\mu$ , where the mean values of parameters range are taken from Table \ref{table-param-values-4-countries}.}
	\label{relative bias of R0}
\end{figure}

Regarding Latin Hypercube Sampling (LHS) samples, the relative bias of $\mathcal{R}_0$ in epidemiology is the percentage of samples that fall into a specific range or category of the fundamental reproduction number $(\mathcal{R}_0)$. Samples of epidemiological parameters, such as $\mathcal{R}_0$, are generated using LHS, and each sample represents a distinct set of parameter values. The relative bias of $\mathcal{R}_0$ is determined by counting the number of samples that fall within a specific range or category, and then dividing that count by the total number of samples \cite{Bifurcation of R0-7, Bifurcation of R0-8}.\\
For example, if we have conducted LHS with a total of 100 samples and we are interested in the relative bias of $\mathcal{R}_0$ values between 2 and 3, we would count how many samples have an $\mathcal{R}_0$ value within that range. If 20 out of the 100 samples fall within that range, the relative bias of $\mathcal{R}_0$ would be 20/100 or 0.2 (or 20\%). The relative bias of $\mathcal{R}_0$ provides insight into the distribution of $\mathcal{R}_0$ values within the LHS samples, allowing us to understand the likelihood of different $\mathcal{R}_0$ values occurring in the population. It helps us in this thesis, to assess the potential for disease spread and design appropriate control measures based on the range of $\mathcal{R}_0$ values observed in the samples.\\
The degree of departure or inaccuracy in predicting the real value of $\mathcal{R}_0$ as a result of changes in contact rates $\beta_{1}$ and $\beta_{2}$ is referred to as the relative bias of $\mathcal{R}_0$. The Figure \ref{relative bias of R0}(a) implies that by altering the contact rates $\beta_{1}$ and $\beta_{2}$ (the frequency and degree of interpersonal contacts), the estimated $\mathcal{R}_0$ can be subject to a relative bias ranging from 0.5 to 6.5. This means that the estimated $\mathcal{R}_0$ value can differ significantly from the true value by a factor of 0.5 to 6.5 with variance 1.85. Furthermore, the maximum values of $\mathcal{R}_0$, within the given context, lie between 0.25 and 4. This implies that under different scenarios and contact rates, the highest estimated $\mathcal{R}_0$ values observed range from 0.25 to 4. These values represent the upper bounds of the potential infectivity or transmission capacity of a specific infectious disease within a susceptible population. That indicates high variation of $\mathcal{R}_0$ with the variation of contact rates $\beta_{1}$ and $\beta_{2}$.
\noindent
From Figure \ref{relative bias of R0}(b), The uncertainty or variability in predicting the fundamental reproduction number $(\mathcal{R}_0)$ in epidemiology is shown by the relative bias of $\mathcal{R}_0$ in the range of 0.85 to 6. From Figure \ref{relative bias of R0}(b), the wide range of relative bias suggests that the estimates of $\mathcal{R}_0$ vary significantly depending on factors such as the accuracy of data, the modeling approach used, and the particular traits of the illness that is being studied. Conversely, the highest values of $\mathcal{R}_0$ that occur between 2.5 and 4 represent the upper bounds of the fundamental reproduction number. These numbers indicate the likelihood that an illness may spread across a group of people. A $\mathcal{R}_0$ value of 2.5 to 4 indicates that, on average, 2.5 to 4 susceptible persons will get the illness from each sick person. Higher values of $\mathcal{R}_0$ indicate a greater potential for disease transmission, while values below 1 indicate that an outbreak is likely to die out. That, indicates that variation of infection rate $\alpha$ and vaccine inefficiency rate $(1-\varepsilon)$, $\mathcal{R}_0$ arise higher rapidly with variance 1.45.

From Figure \ref{relative bias of R0}(c), The relative bias of $\mathcal{R}_0$, ranging from 1 to 4.5, indicates the extent to which the calculated $\mathcal{R}_0$ value deviates from the true value. When considering the maximum values of $\mathcal{R}_0$, which fall within the range of 2 to 3.5, it suggests that the disease has the potential to spread relatively easily within a population. Higher values of $\mathcal{R}_0$ indicate a greater likelihood of sustained transmission and a larger outbreak size, while lower values suggest a reduced capacity for the disease to spread widely. It is important to note that $\mathcal{R}_0$ is influenced by various factors, including the recovery rates $(\gamma)$ and treatment rates $\gamma_1$ of the disease; also $\mathcal{R}_0$ fluctuates by variance 1. By altering these rates, we can affect the magnitude of $\mathcal{R}_0$ and perhaps slow the disease's progression. But this variation is lower than the Figures \ref{relative bias of R0}(a) and (b).
From Figure \ref{relative bias of R0}(d), the given information states that by varying recruitment rates $\mu$ (the rate of population entry of vulnerable people) and disease induced death rates $\delta$, the relative bias of $\mathcal{R}_0$  falls within the range of 1.5 to 4.2. The discrepancy between the estimated and true values of $\mathcal{R}_0$ is referred to as relative bias. Additionally, it illustrates that the maximum values of $\mathcal{R}_0$  are within the range of 2 to 3. This means that, under the specified conditions, the highest estimated $\mathcal{R}_0$  values would be between 2 and 3 by variation of $\mu$ and $\delta$. These estimates indicate the potential for the disease to spread rapidly within the population, with each infected individual potentially infecting 2 to 3 others on average. It also reflects that, the variation is not as high as the variation of $\mathcal{R}_0$ with respect to contact rates $\beta_{1}$, $\beta_{2}$,infection rates $\alpha$ as shown in Figures \ref{relative bias of R0}(a) and (b).

\subsection{Global Sensitivity Analysis}
In the global sensitivity method, the input parameter values vary simultaneously in their respective parameter space and can capture parameter interactions. The goal is to rank the input parameters to identify which parameter is the most influential in controlling the disease outbreak dynamics. This method can be used for both linear and non-linear models. For statistical measure, PRCC test is most commonly used.
\\[0.1cm]
\textbf{Partial Rank Correlation Coefficient}
\\[0.1cm]
Latin Hypercube Sampling (LHS)-Partial Rank Correlation Coefficient (PRCC) is one of the most effective sensitivity analysis methods for examining parameter uncertainty. Numerous mathematical model types, such as deterministic or stochastic models with continuous or discrete characteristics, can be used with this process. The degree of a linear link between an input and an output may be measured using correlation \cite{Hopf LEAST LHS PRCC-9}. Following is the calculation of a correlation coefficient (CC) between $x_i$ and $y$:
\begin{align*}
	r_{x_iy}=\frac{Cov(x_i,y)}{\sqrt{Var(x_i)Var(y)}}=\frac{\sum_{i}(x_i-\tilde{x})(y_i-\tilde{y})}{\sqrt{\sum_{i}(x_i-\tilde{x})^2\sum_{i}(y_i-\tilde{y})^2}}\;\;\;\text{for  }i=1,2,\cdots,k
\end{align*}
$r_{x_i,y}$ varies between -1 and +1; where $Cov(x_i,y)$ represents the covariance between $x_i$ and $y$. The correlation between two variables after one or more extra factors have been taken out of the picture is known as partial rank correlation \cite{Hopf LEAST LHS PRCC-10, Hopf LEAST LHS PRCC-11}. This command is intended just in the event that there is one more variable. In this instance, the standard rank correlations between the three variables may be used to calculate the partial rank correlation as follows:
\begin{align*}
	r_{xy,z}=\frac{r_{xy}-r_{xz}r_{yz}}{\sqrt{(1-r_{xz}^2)(1-r_{yz}^2)}}
\end{align*}
where the correlation between $x$ and $y$ is shown by $r_{xy}$. The following formula is used to get the partial rank correlation coefficient between the two variables:
\begin{align*}
	r=1-6\sum_{i-1}^{N}\frac{D_i}{N(N^2-1)}
\end{align*}
where $N$ is the sample size and $D_i$ is the difference between the rankings given to the corresponding pairings. PRCC is a sample-based technique to investigate the relationship between a model's parameters and $\mathcal{R}_0$ within a specified range of parameter values. Finding important factors that significantly affect $\mathcal{R}_0$ for model prediction and ranking them based on accurate forecasts are the objectives of PRCC sensitivity analysis. Every likely sample of parameters is taken at random within LHS. A result of +1 indicated a perfect positive linear connection, meaning that the parameter has a significant influence on $\mathcal{R}_0$, since the PRCC ranges from -1 to +1. Additionally, a complete negative linear relationship is shown by a value of -1, and no linear relationship is indicated by a value of 0. $\rho$(rho) in the MATLAB formula represents correlation for each LHS parameter. For instance, correlation between $\beta_1$ and $\mathcal{R}_0$ is shown by $\rho=0.85$ (Table \ref{table-PRCC-4-Countries}). In addition, provides the p-value as a second output and verify these results to see if the correlations are statistically significant. There isn't a substantial link between these parameter pairings, as seen by the huge p-values. In order to account for the impacts of a third set of variables, test for PRCC between pairs of variables $x$ and $y$ input matrices as follows:
\begin{align*}
	\mathbf{[rho,pval]=partialcorr(x,y,z)}
\end{align*}
And returns the covariance matrix of $[x,z]$ is: $\displaystyle S=\begin{pmatrix}
	S_{xx}&S_{xz}\\S_{xz}&S_{zz}
\end{pmatrix}$\\
where, $S_{xx}$ is $\rho$ (PRCC) value corresponding to $x$ and $z$.\\
By utilizing partial rank correlation coefficients (PRCC) and the associated p-values, partial rank correlation analysis allows us to ascertain the level of uncertainty that an LHS parameter contributes to the model. The size and statistical significance of a parameter's PRCC value indicate how much of an impression the parameter has on the model's prediction. The most significant and impactful parameters are those with big PRCC values ($>0.5$\; or \;$<-0.5$) and associated modest p-values ($<0.05$); these indicate that the PRCC value is statistically significant and differs from zero. The LHS parameter affects the result measure more strongly the closer the PRCC value is to $+1$ or $-1$. The qualitative relationship between the input and output variables is shown by the sign. The LHS parameter is inversely related to the result measure when it has a negative sign.
\begin{figure}[H]
	\centering  
	\subfloat[]{\includegraphics[width=3. in]{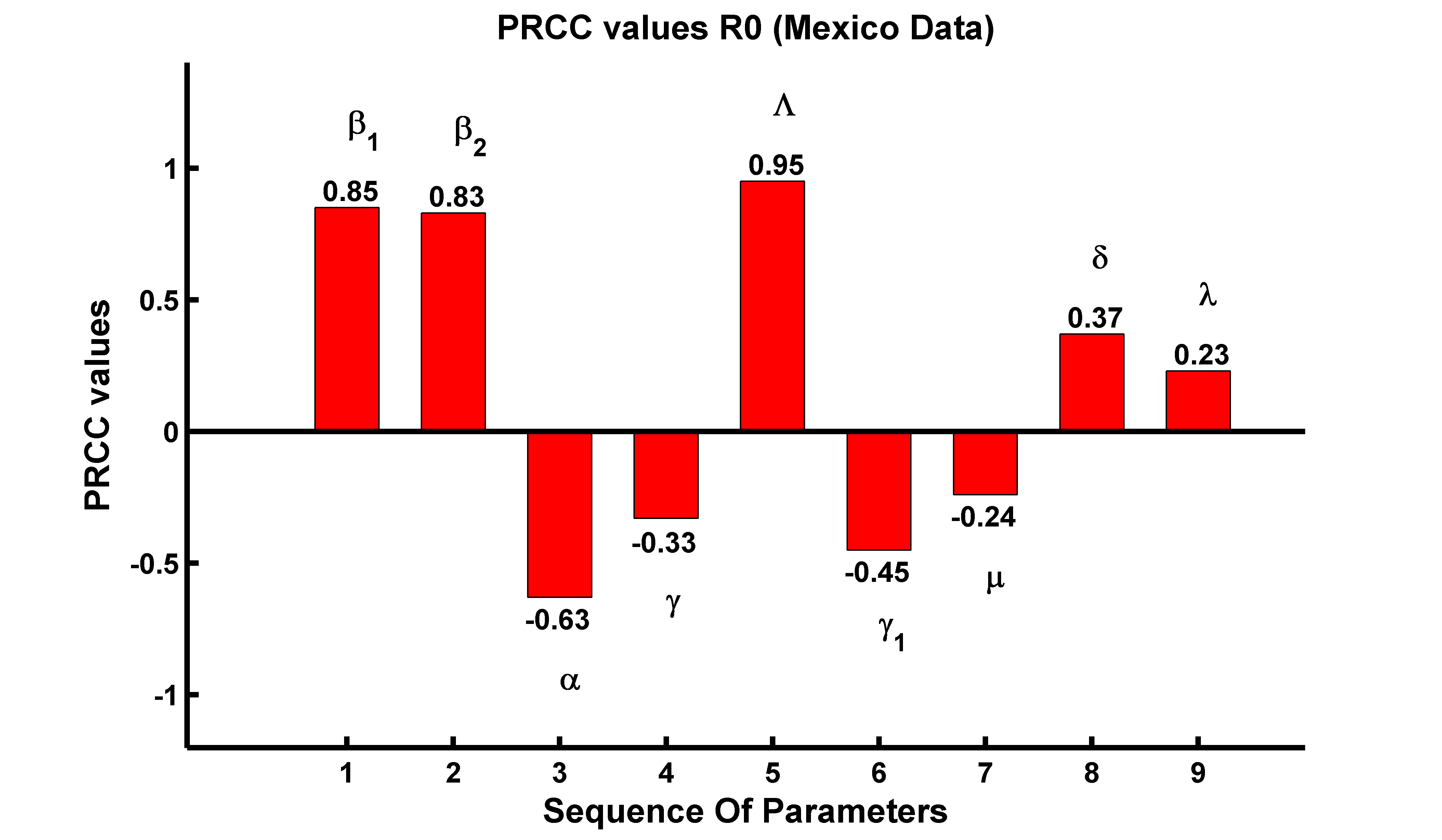}}
	\subfloat[]{\includegraphics[width=3. in]{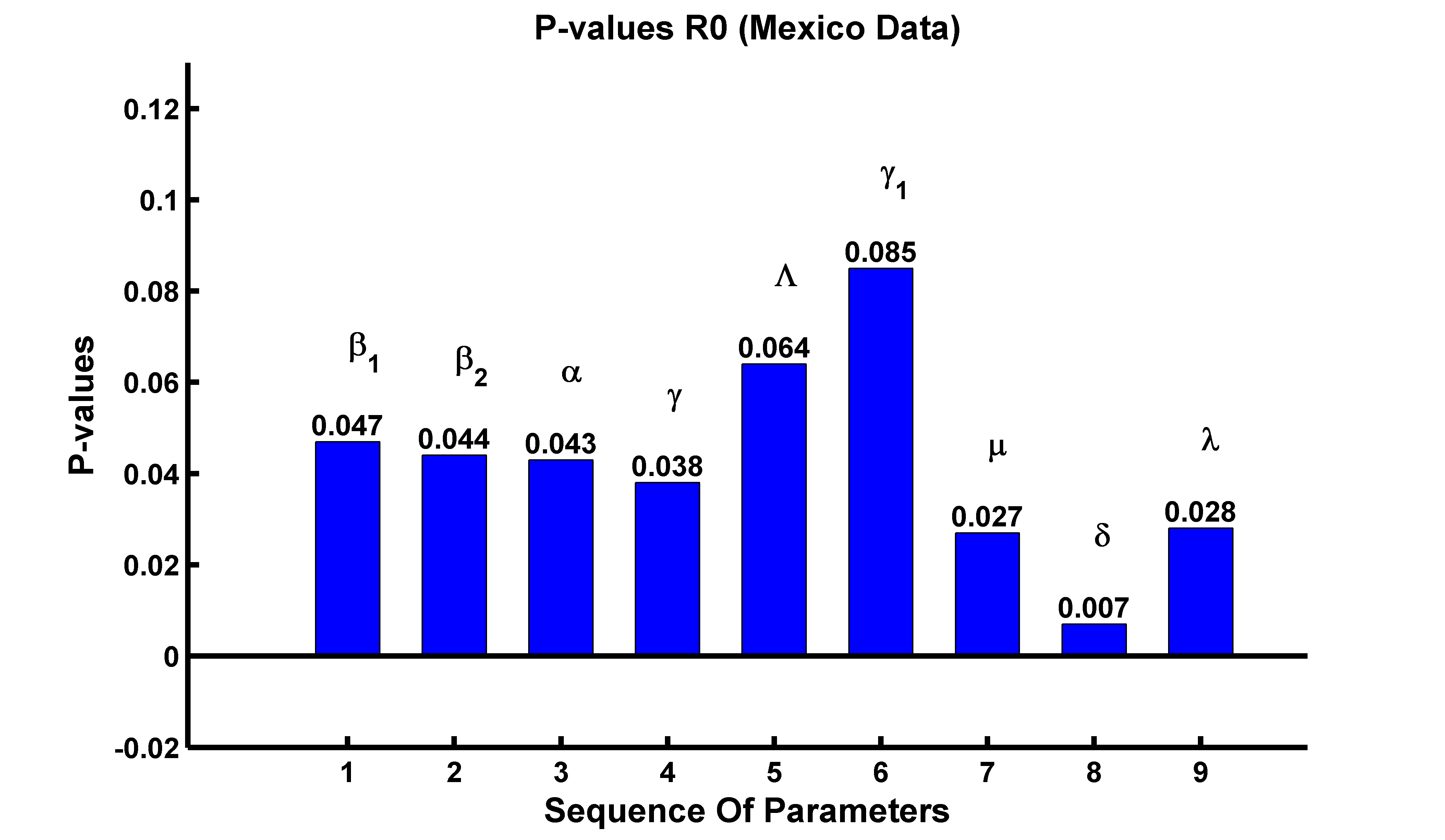}}
	\caption{(a) PRCC and (b) P-values plot for Mexico data of Influenza disease, where mean values of parameters are taken from Table \ref{table-param-values-4-countries}.}
	\label{PRCC-Pval-Mexico}
\end{figure}
\noindent
Based on our proposed model, maximum and minimum values are found for each of the nine LHS parameters displayed in Figure \ref{latine-hypercube-sampling}. It should be noted that the baseline value for each LHS parameter has been established to be somewhere in the middle of the range between the parameter's minimum and maximum values (Table \ref{tableparameter}). Meanwhile, we estimated the four set of parameter values from four countries which are presented in Table \ref{table-param-values-4-countries}. A portion of the parameter's baseline values are approximated, while others are derived from the reference shown in Table \ref{tableparameter}. In Table \ref{table-param-values-4-countries}, the best or well-fitting curves are used to estimate the lowest and maximum values of each parameter.

\begin{figure}[H]
	\centering  
	\subfloat[]{\includegraphics[width=3. in]{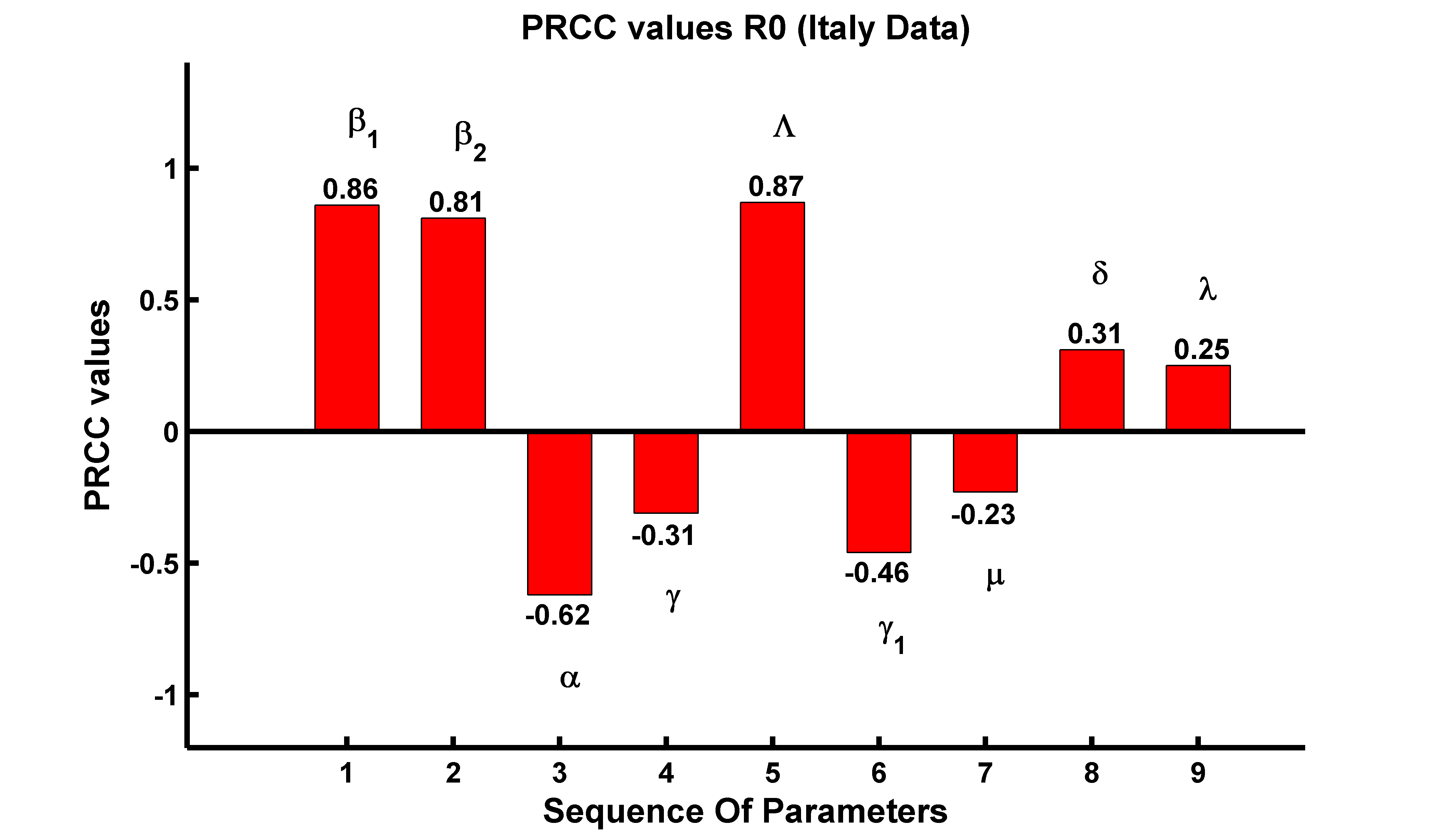}}
	\subfloat[]{\includegraphics[width=3. in]{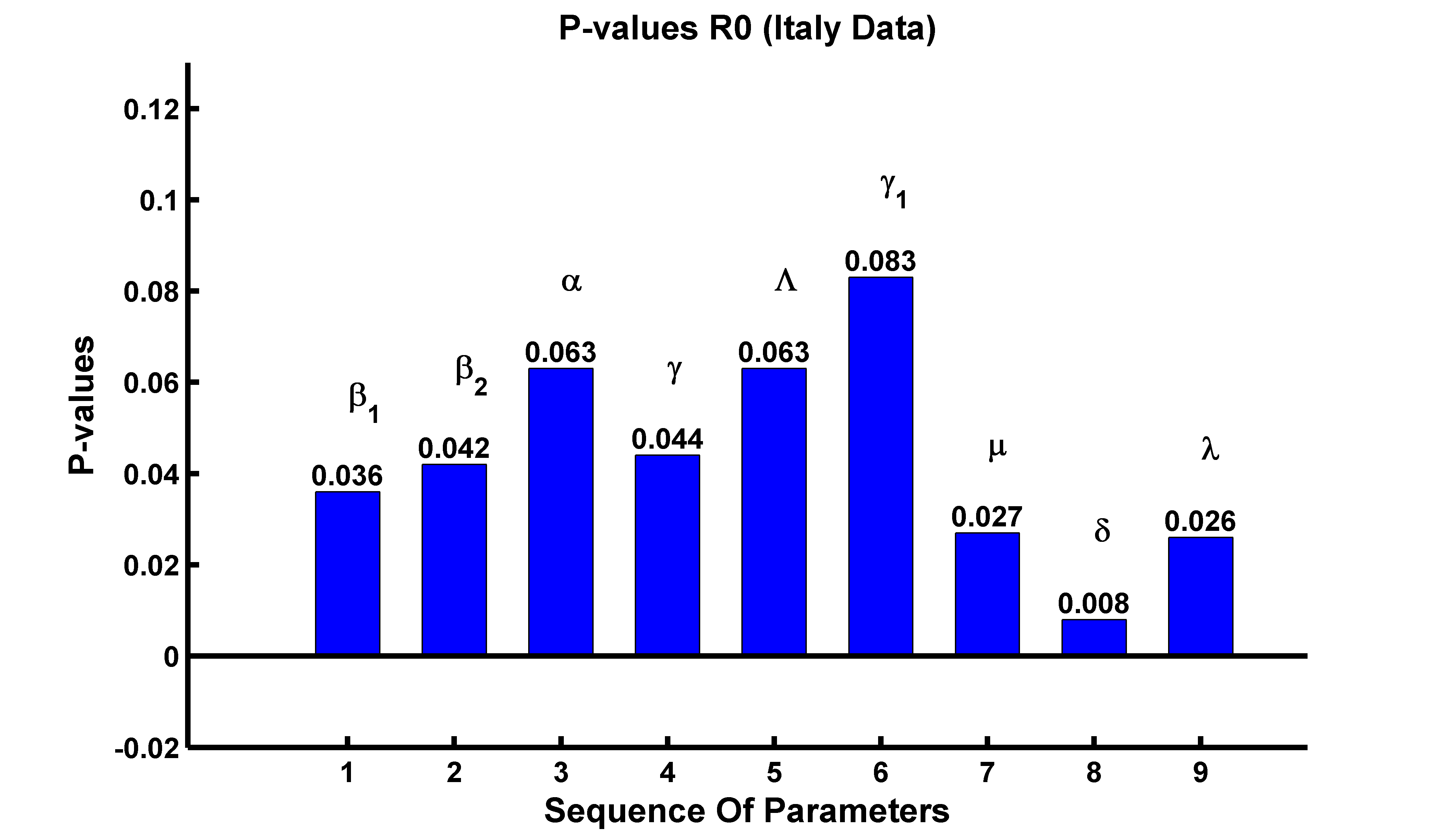}}
	\caption{(a) PRCC and (b) P-values plot for Italy data of Influenza disease, where average values of parameters are taken from Table \ref{table-PRCC-4-Countries}.}
	\label{PRCC-Pval-Italy}
\end{figure}

\begin{figure}[H]
	\centering  
	\subfloat[]{\includegraphics[width=3. in]{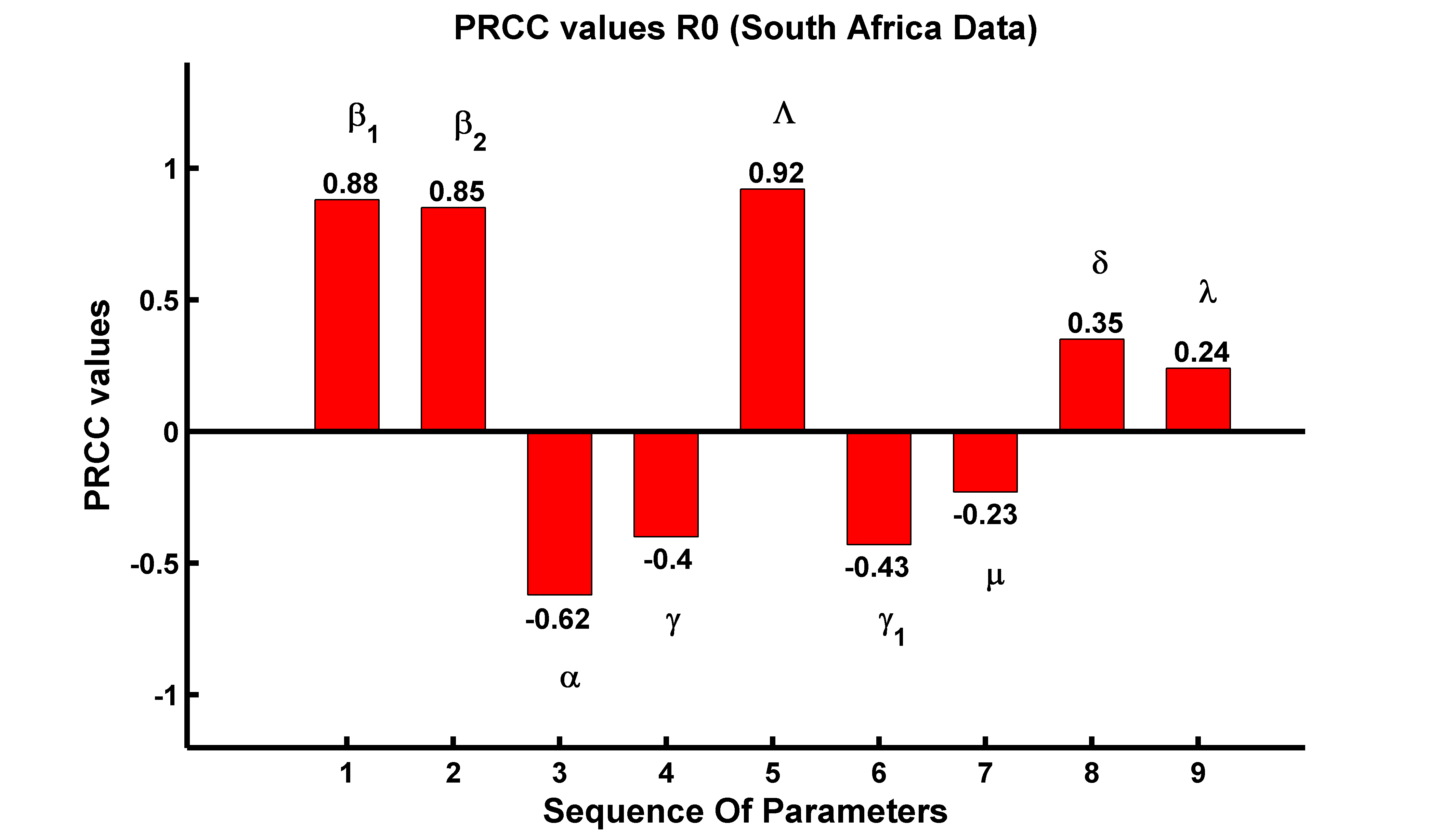}}
	\subfloat[]{\includegraphics[width=3. in]{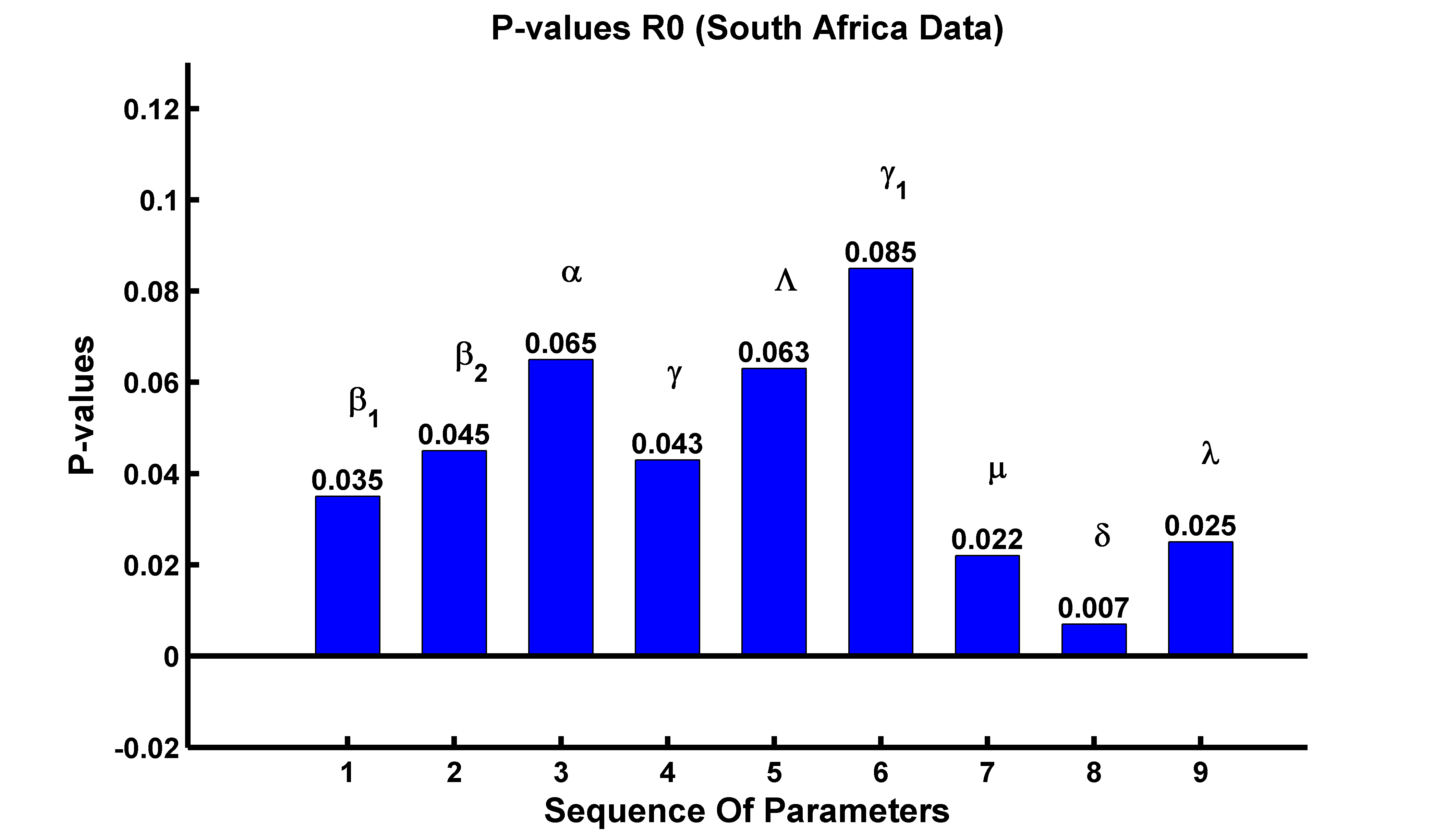}}
	\caption{(a) PRCC and (b) P-values plot for South Africa data of Influenza disease, where mean values of parameters are taken from Table \ref{table-PRCC-4-Countries}.}
	\label{PRCC-Pval-South-Africa}
\end{figure}

To ascertain the model's \eqref{new_model} resilience to parameter values which aids in identifying the most important parameters in model dynamics, we performed a sensitivity analysis in this section. The Latin Hypercube Sampling (LHS) scheme, a Monte Carlo sampling approach, is applied. It samples 100 values for each input parameter using a uniform distribution throughout the range of physiologically plausible values shown in Figure \ref{latine-hypercube-sampling} and Table \ref{table-param-values-4-countries}. For the system of differential equations described in \eqref{new_model}, 50 model simulations were carried out by randomly picking paired sampled values for all LHS parameters. Figure \ref{PRCC-Pval-Mexico}, Figure \ref{PRCC-Pval-Italy}, 
Figure \ref{PRCC-Pval-South-Africa} provides the Partial Rank Correlation Coefficients (PRCC) and associated p-values for each parameter, indicating the corresponding non-linear but unmodulated link between the model state variables and each parameter. Moreover, if the associated p-value is smaller than 1\%, the result is statistically significant.

By analyzing the data of four countries, presented in Figure \ref{PRCC-Pval-Mexico}, Figure \ref{PRCC-Pval-Italy},  Figure \ref{PRCC-Pval-South-Africa}, we can determine which parameters, $\beta_1$, $\beta_2$, $\alpha$, $\gamma$, $\gamma_1$, and $\Gamma$, have a significant impact on the dynamics of the model. The model variables are influenced in different ways by each major parameter. Using the rankings acquired for the outcome measures (i.e., the total number of exposed and infected persons) and the LHS parameters, we do a multilinear regression analysis. The regression coefficients are then obtained by doing a regression analysis on these rankings.
\noindent
As these regression coefficients provide an indication of the model's sensitivity to the LHS parameters, we then get the PRCC values to determine the degree of correlation between each LHS parameter and each outcome measure. For all exposed and infected individuals, we have PRCC plots with PRCC values and p-values in Figures \ref{PRCC-Pval-Mexico}, \ref{PRCC-Pval-Italy},  and \ref{PRCC-Pval-South-Africa}. Furthermore note that there is a high link indicated by the PRCC plots for the two outcome measures.
\noindent
\subsection{Analysis of the Global Sensitivity Result}
The most significant parameters for the model \eqref{new_model} are those that, according to the PRCC analysis, have PRCC values $>0.5$ (for direct relations) or $<-0.5$ (for inverse relations) and matching modest p-values $(<0.05)$. The above Figures \ref{PRCC-Pval-South-Africa},  \ref{PRCC-Pval-Italy}, and \ref{PRCC-Pval-Mexico} are represented in the PRCC plot. Parameters are plotted on the $x$-axis, while PRCC values are plotted on the y-axis. The PRCC and p-values of the parameters with $\mathcal{R}_0$ are shown in Table \ref{table-PRCC-4-Countries}.
\begin{table}[H]
	\begin{center}
		\caption{Using the model parameters \eqref{new_model} and accompanying p-values, PRCC of $\mathcal{R}_0$ was assessed at the parameter values provided in Table \ref{table-param-values-4-countries} in order to identify the most significant LHS parameters.}
		\scriptsize
		\label{table-PRCC-4-Countries}
		\begin{tabular}{|c|c|c|c|c|c|c|c|c|}
			\hline\noalign{\smallskip}
			\multirow{2}{*}{\textbf{Parameters}}& \multicolumn{2}{|c|}{\textbf{Mexico}}&\multicolumn{2}{|c|}{\textbf{Italy}}&\multicolumn{2}{|c|}{\textbf{South Africa}}\\
			\cline{2-7}
			&\textbf{PRCC} &\textbf{p-value}&\textbf{PRCC}&\textbf{p-value}&\textbf{PRCC} &\textbf{p-value}\\
			\noalign{\smallskip}\hline\noalign{\smallskip}
			$\beta_1$&0.85&0.047&0.86&0.036&0.88&0.035\\
			$\beta_2$&0.83&0.044&0.81&0.042&0.85&0.045\\
			$\alpha$&-0.63&0.043&-0.62&0.063&-0.62&0.065\\
			$\gamma$&-0.33&0.038&-0.31&0.044&-0.4&0.043\\
			$\Lambda$&0.95&0.064&0.87&0.063&0.92&0.063\\
			$\gamma_1$&-0.45&0.085&-0.46&0.083&-0.43&0.085\\
			$\mu$&-0.24&0.027&-0.23&0.027&-0.23&0.022\\
			$\delta$&0.37&0.007&0.31&0.008&0.35&0.007\\
			$\lambda$&0.23&0.028&0.25&0.026&0.24&0.025\\
			\noalign{\smallskip}\hline
		\end{tabular}
	\end{center}
\end{table}

The PRCC test result from Figures  \ref{PRCC-Pval-Mexico}, \ref{PRCC-Pval-Italy}, \ref{PRCC-Pval-South-Africa} and Table \ref{table-PRCC-4-Countries} demonstrates that,
\begin{enumerate}
	\item For a range of uniformly distributed parameters the most sensitive parameters for the outcome measure with $\mathcal{R}_0$ are $\beta_1$, $\beta_2$, $\alpha$, $\gamma$ and $\gamma_1$. For example, from Table \ref{table-PRCC-4-Countries}, PRCC between $\mathcal{R}_0$ and $\beta_1$ = 0.83 (Mexico), 0.88 (South Africa) and p-values are 0.047$<0.05$ and 0.035$<0.05$, respectively. Another strong positive contributor is  $\beta_2$ having PRCC 0.83. By controlling these parameters, disease transmission can be reduced more effectively.
	\item The threshold quantity $\mathcal{R}_0$ has medium negative influence with $\gamma$, $\gamma_1$, $\mu$ and $\delta$. Thus $\mathcal{R}_0$ increases (or decreases) rapidly while these parameter rates decreases (or increases). Further smaller p-values($<0.05$) for $\beta_1$, $\beta_2$, $\alpha$, $\gamma$ and $\gamma_1$ indicates the PRCC value corresponding parameters are statistically more significant. Table \ref{table-PRCC-4-Countries} shows, these parameters are likely highly contributors to uncertainty.
	\item $\mathcal{R}_0$ correlated negatively with $\alpha$, $\mu$, $\delta$ also their p-values very low. But they are medium contributor to uncertainty. That indicates quarantine and isolation can be effective policy to control strategy. Human cure rate from disease $\gamma$, treatment rate $\gamma_1$ also have great significance on controlling the disease. 
	\item As $\mathcal{R}_0$ has positive influence with $\gamma$, it indicates that, the more recruitment occurs, disease elimination process takes longer time. We can conclude that, proper vaccination, treatment policies, mitigating contact with $E(t)$ and $I(t)$ humans, are closely related to control influenza outbreak.
\end{enumerate}

\subsection{Validity of PRCC Result }
It is necessary requirement to show the relationship between model parameter and PRCC results to check parameters either monotonic or not. LHS/PRCC analysis is a statistical technique that helps identify which parameters have the most influence on $\mathcal{R}_0,$ also on the models behaviour \cite{Stability Bound-25}. In the context of checking validity, a monotonic relationship means that as the value of a parameter increases or decreases, the model's output consistency follows the same pattern without any reversal or fluctuations which provides meaningful and reliable results. From Figure  \ref{PRCC-Pval-Mexico},  \ref{PRCC-Pval-Italy}, \ref{PRCC-Pval-South-Africa}, the PRCC measures the strength or the monotonic relationship between the parameter and the output. It quantifies how much that parameter's variation contributes to the output variation while controlling the effects of other parameters.\\
In context of relationship with model phase, it allows for a systematic exploration of the parameter space and provides insights into which parameters should be prioritized for further investigation or optimization. In PRCC, we see that some parameters value +0.5 to +1; having strong positive influence on the compartments of exposed and infected class. Meanwhile, which parameters value -1 to -0.5 have strong negative influence on the compartments of exposed and infected. But $\gamma$ and $\gamma_1$ having negative correlation; increases the recovery and treatment rate which play a significant role to curtail the outbreak. The immediate summarized results are:
\begin{enumerate}
	\item \textbf{Effect of $\beta_1\;\;\text{and}\;\; \beta_2$: } 
	In PRCC, $\beta_1$ and $\beta_2$ have strong correlation 0.85 and 0.83 with $\mathcal{R}_0$. Thus, by increasing quarantine results into reduction of of $\beta_1$ and $\beta_2$, that controls $\mathcal{R}_0$.  As a result monotonic increasing behaviour of $\beta_1$ and $\beta_2$ are  verified.
	\item \textbf{Effect of $\alpha$: }  
	It effects the infected compartment as infected class rapidly increase to its peak point according to proportion of $\alpha$. From, PRCC value of $\alpha$ with $\mathcal{R}_0$ (-0.63); the monotonic behaviour of $\alpha$ are justified.
	\item \textbf{Effect of $\gamma$ and $\gamma_1$: } From the PRCC analysis, the PRCC of $\gamma$ and $\gamma_1$ with $\mathcal{R}_0$ are -0.33 and -0.45. That indicates the likely negative correlation. 
	Further, we see that, by increasing $\gamma_1$  30\%, 50\%, 80\%, reflects that infected persons gain a herd immunity by vaccination and decrease from its pick point and treatment class increases. Therefore, $\mathcal{R}_0$ can be controlled effectively. This shows that the monotonic decreasing behaviour of $\gamma$ and $\gamma_1$ are valid.
\end{enumerate}

\section{Case Study of Clinical Data vs Simulation}\label{Section-Numerical-Simulation}
In this section, we analyze the model \eqref{new_model} numerically to support the analytical results presented in previous sections.  Numerical simulations enable the incorporation of real data, such as population demographics, disease parameters, and contact patterns, to create more accurate and realistic models. This also facilitate parameter estimation by comparing model predictions with observed data, optimizing parameters to minimize the difference between simulated and actual outcomes \cite{Stability Bound-4, Stability Bound-25}. By modifying parameters or implementing different control measures within the simulation, researchers are able to evaluate the efficacy of initiatives like contact tracking, social distance, and vaccination programs. Numerical models also make it possible to predict how an epidemic would spread under various conditions, which may help decision-makers make well-informed decisions on how best to allocate resources, provide healthcare, and carry out policies. 

\subsection{A Case Study of Influenza in Italy}
\begin{figure}[H]
	\centering  
	\subfloat[]{\includegraphics[width=3. in]{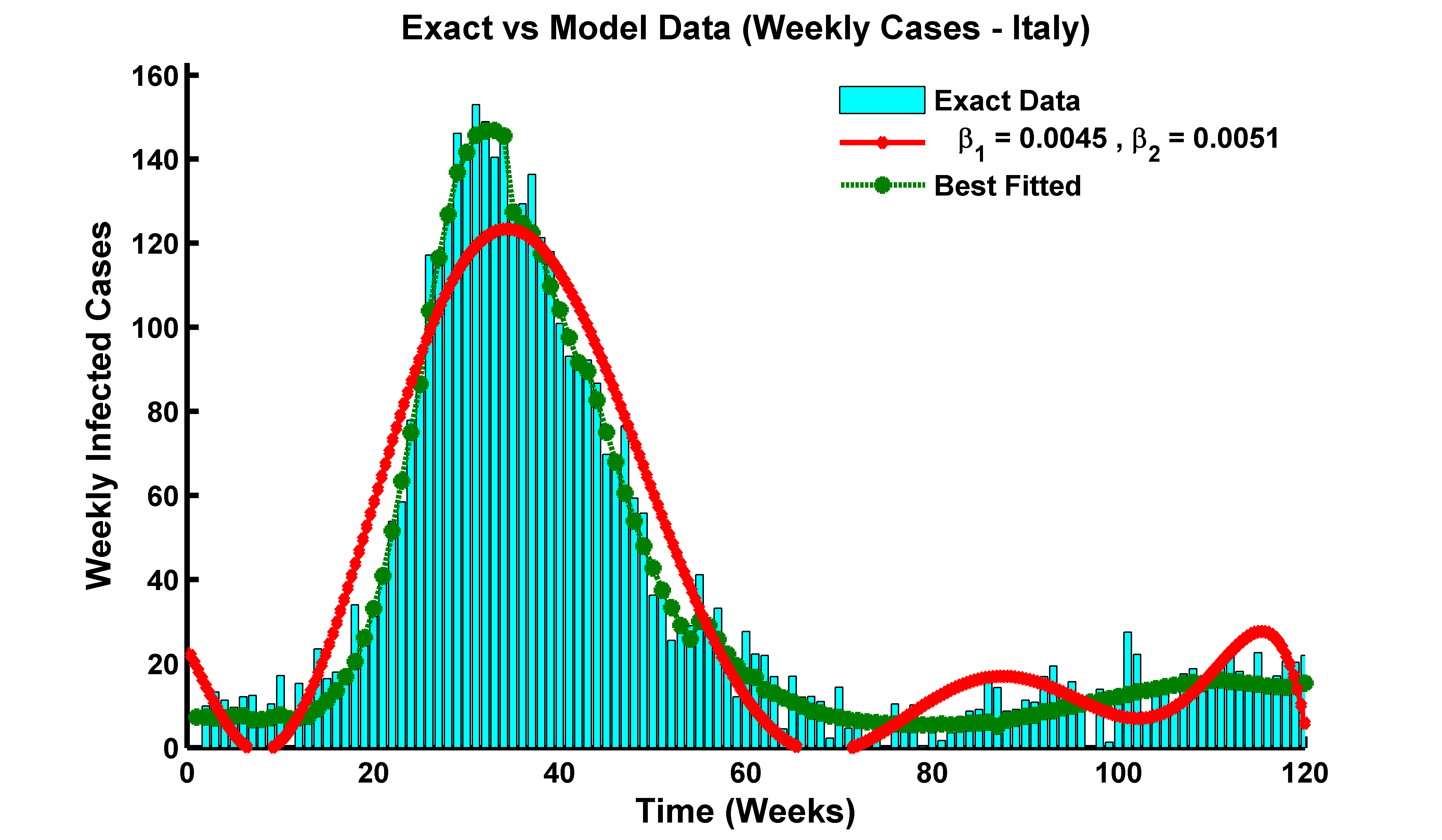}}
	\subfloat[]{\includegraphics[width=3. in]{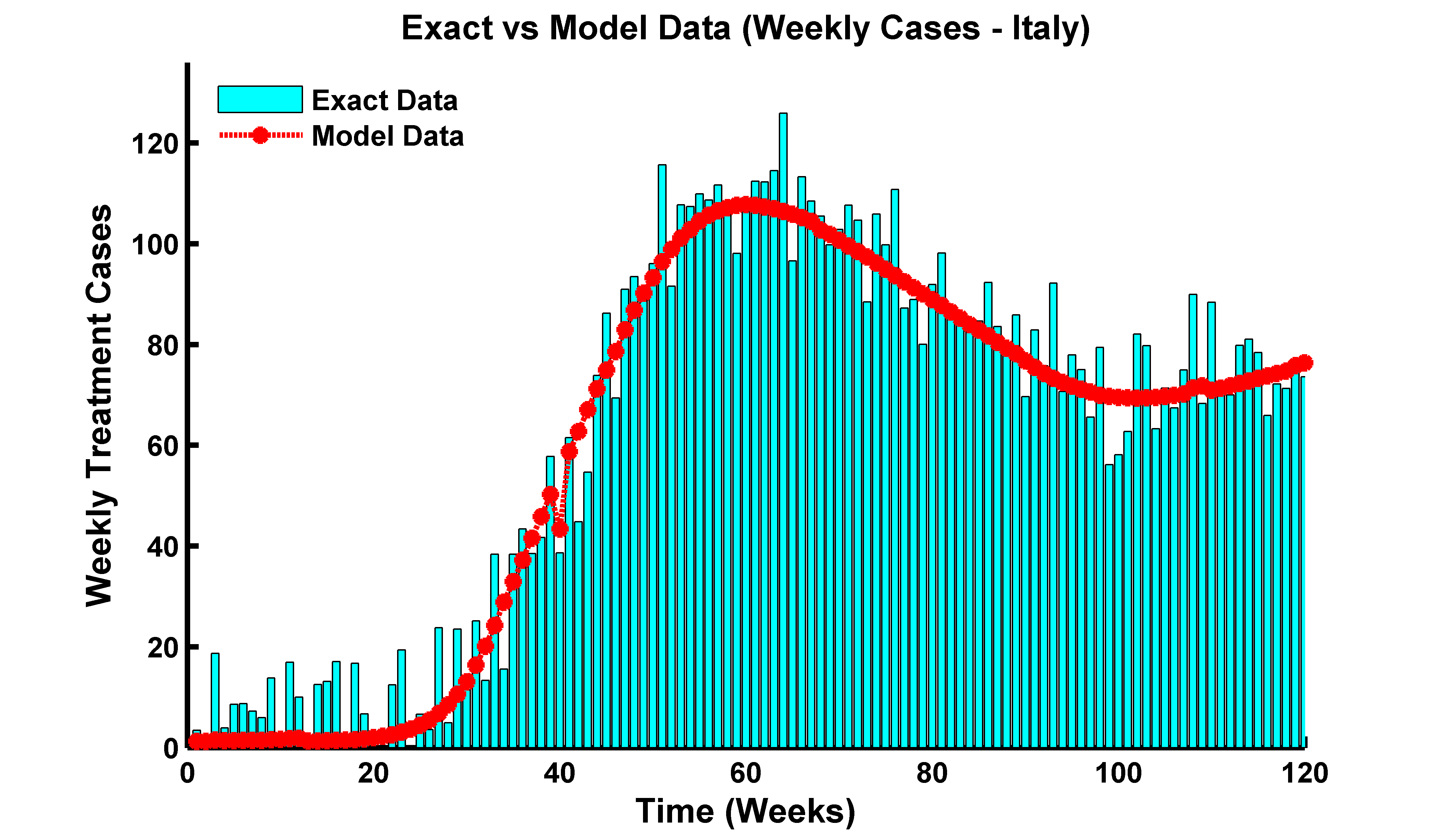}}\\
	\subfloat[]{\includegraphics[width=3. in]{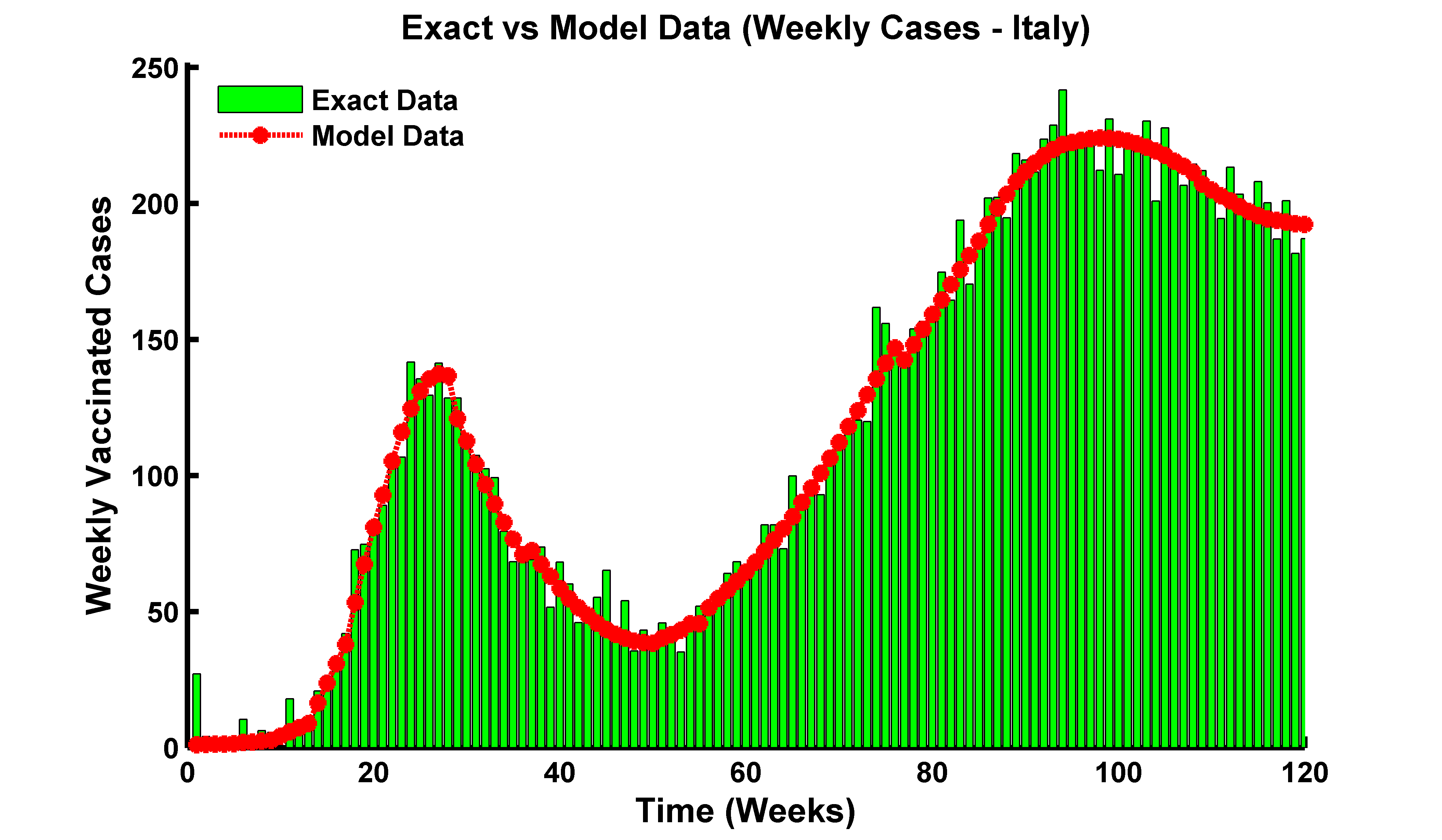}}
	\subfloat[]{\includegraphics[width=3. in]{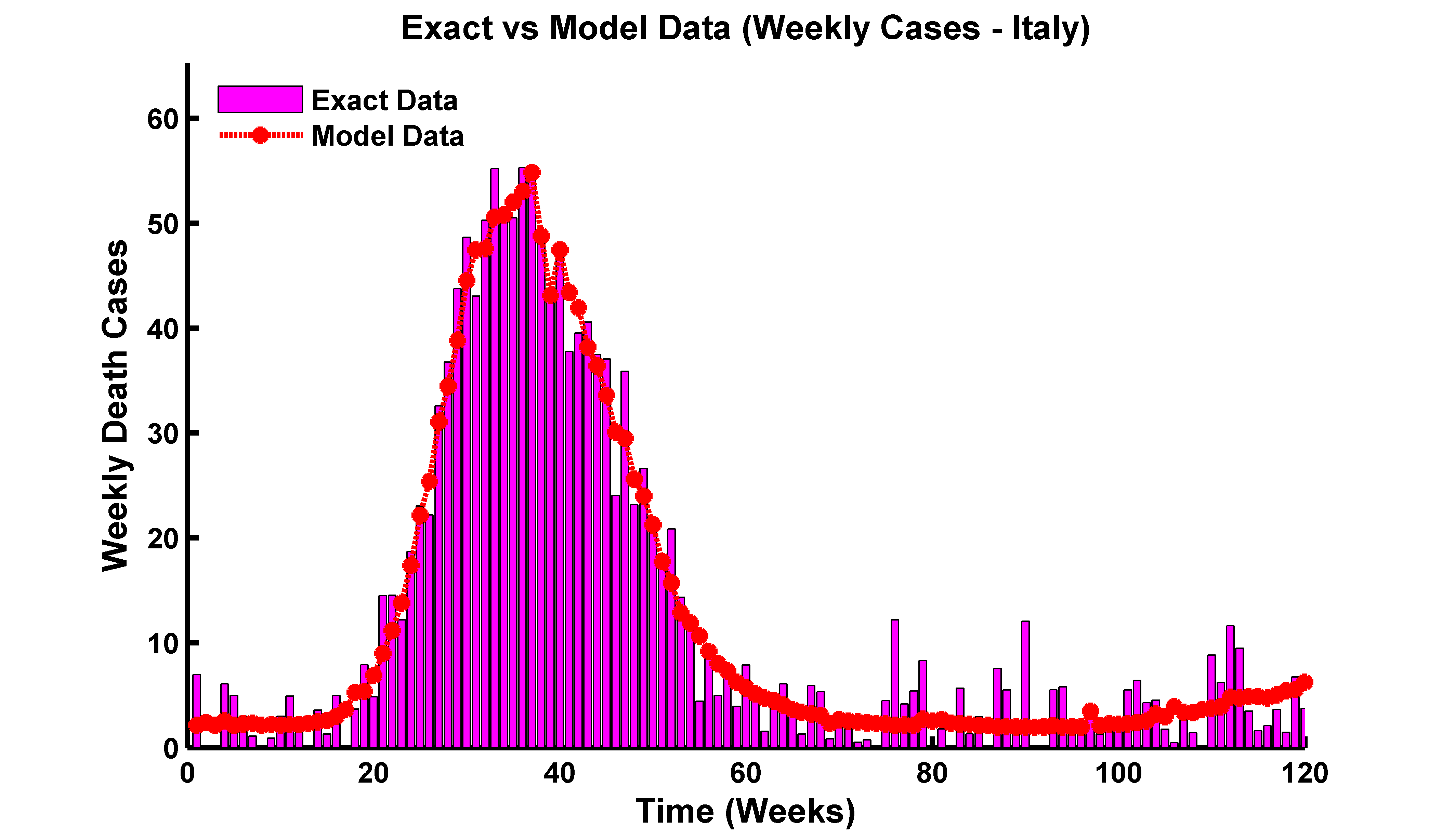}}
	\caption{Model data fitting to Italy with (a) weekly infected data (b) weekly progression of treatment to control disease burden (c) weekly vaccination programme to control disease burden (d) weekly death rates, where the fitted parameter sets are : $\beta_1=0.0075, \beta_2= 0.0081, \alpha=0.78, \gamma=0.63, \gamma_1=0.35, \lambda=0.42, \mu=0.03, \delta=0.29$.}
	\label{case-study-Italy}
\end{figure}

In this section, we have examined recent data on influenza infection cases in Italy, spanning the period from 1 October 2020 to 31 March 2023, in order to study the trend and test our model. The websites of the CDC and WHO were used to acquire all the data \cite{CDC, WHO}. To assess the model's applicability and estimate disease for the near future, we took into account a total of 120 weekly data points. The initial population is taken as $S(0)=500$, $V(0)=1$, $E(0)=1$, $I(0)=1$, $R(0)=1$ and $T(0)=0$. The model outcomes predicted with real data presented in Figure \ref{case-study-Italy}, where the weekly reported cases are depicted. The linear regression analysis with parameters in Table \ref{table-param-values-4-countries}, reveals that the model fits the weekly infected cases with excellent agreement. Thus, we demonstrate approximately 70\% accuracy to track the original confirmed cases. As Italy located in Southern Europe, it experiences a diverse range of weather patterns due to its geographical location and varying topography. Generally speaking, Italy experiences hot, dry summers and warm, rainy winters due to its Mediterranean environment. The World Health Organization reports that throughout the seasons of 2005–2006 and 2014–2015, an estimated average of about 4,800.000 infected cases were documented. The majority belonged to the 0–5 age group. 40,000 contaminated cases from the national database were analyzed, with samples from instances that happened between 2000/2001 and 2011/2012 being examined. Infections with influenza A and B viruses have been reported in Italy in many waves.11,488 influenza cases total 9,842 influenza A cases and 1,646 influenza B cases were verified during the research period. Most influenza seasons saw the co-circulation of the influenza A and B viruses, with the number of influenza B illnesses either matching or surpassing that of influenza A infections in three of the research period's seasons (2001/02, 2007/08, and 2012/13). During this time, Italy had a 1.5-fold rise in the number of samples tested for influenza viruses by PCR, from an average of 2,774 each season in 2000–2007 to an average of 4,312 per season in 2008–2012. Older children seemed to be more susceptible to influenza B, whereas young and older adults were more susceptible to influenza A(H1N1), and the elderly were more susceptible to influenza A(H3N2). Out of all the analyzed samples, the average percentage of influenza B cases was 23\% (range $<$ 1-78\%).

Beginning in 2021, there were 5.5 percent reported cases of influenza in Italy and 1.5 percent reported deaths. The authorities made the initial decision to implement the immunization program, which is still being done today in hospitals, community vaccination centers, clinics, and pharmacies. The vaccinations worked after the prescribed dosage was given. Additional vaccines that have been authorized include Fluarix, FluLaval, Fluzone, and Afluria quadrivalent. These vaccinations are safe to provide to infants as young as six months. These vaccines have an efficacy of roughly 66\%, 64\%, 67\%, and 62\%, respectively. The Figures \ref{case-study-Italy}(a) (b) (c) and (d) show the agreement between the model's prediction and the actual data. The results show how well-fitting the proposed model is.

Analyzing for total cases of population per 1000, from the Figure \ref{case-study-Italy}(a) we can see that, from 0 to 15 weeks, the weekly reported cases of seasonal flu in Italy are in the lower half. Week 20: There are about 35 reported cases. Following that, the infection grew gradually from week 21 to week 38. The highest number of infected cases (144 per 1000) were reported at weeks 32 and 38. Along weeks 30 to 50, seasonal flu cases were at an all-time high. Due to appropriate care and isolation, the illness epidemic marginally decreased after 50 weeks. After 80 weeks, the reported instances (per thousand) range from 25 to 40.
\noindent
In Figure \ref{case-study-Italy}(b), the cases for the weekly therapy are displayed. Which suggests that when more cases of infection arise, treatment policies are being raised. Between 80 and 110 treatment instances are recorded from weeks 35 to 80. Weeks later, 80, this rate has changed over time and now varies between 60 and 70 (per thousand).

In Figure \ref{case-study-Italy}(c), with our suggested approach, the ongoing vaccination cases in Italy are projected. Our model predicts that immunization rates are highest between 22 and 33 weeks. A preventive policy was applied to 140 populations (per 1000) on average. After then, the number of vaccination cases increased along with the number of infected patients. Between weeks 65 and 120, there were documented vaccination cases totaling between 150 and 220, which contributed to the outbreak's efficient suppression.

On the other hand, Figure \ref{case-study-Italy}(d) describes the recent death cases in Italy. According to data analysis, from October 1, 2021, to week 27–45, the majority of fatality cases take place. The death rate ranges from 0 to 20 (per 1000) from weeks 0 to 25. After 60 weeks, the number of deaths decreased and now ranges from 0 to 15 (per 1000).
\subsection{A Case Study of Influenza in South Africa}
\begin{figure}[H]
	\centering  
	\subfloat[]{\includegraphics[width=3. in]{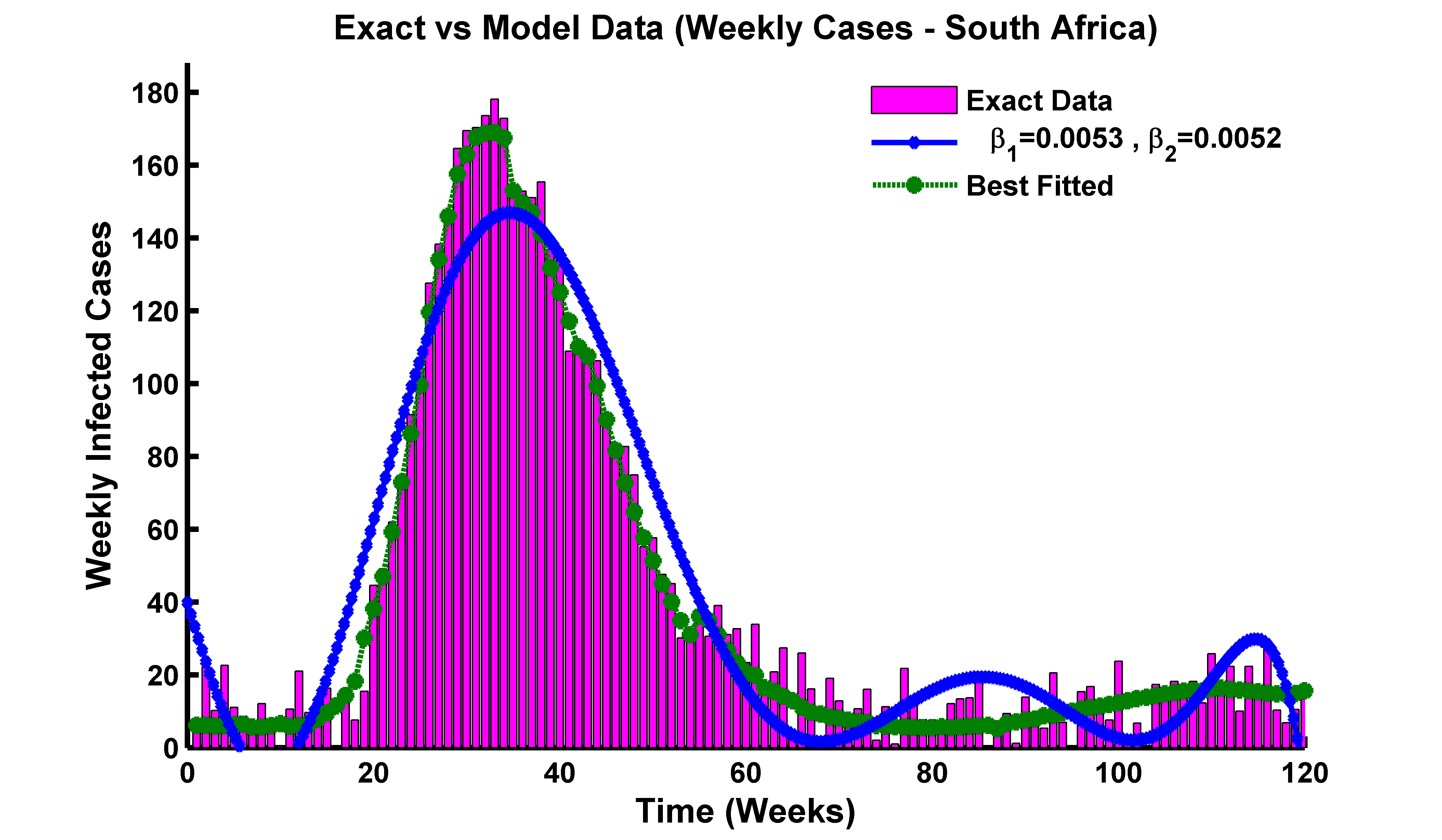}}
	\subfloat[]{\includegraphics[width=3. in]{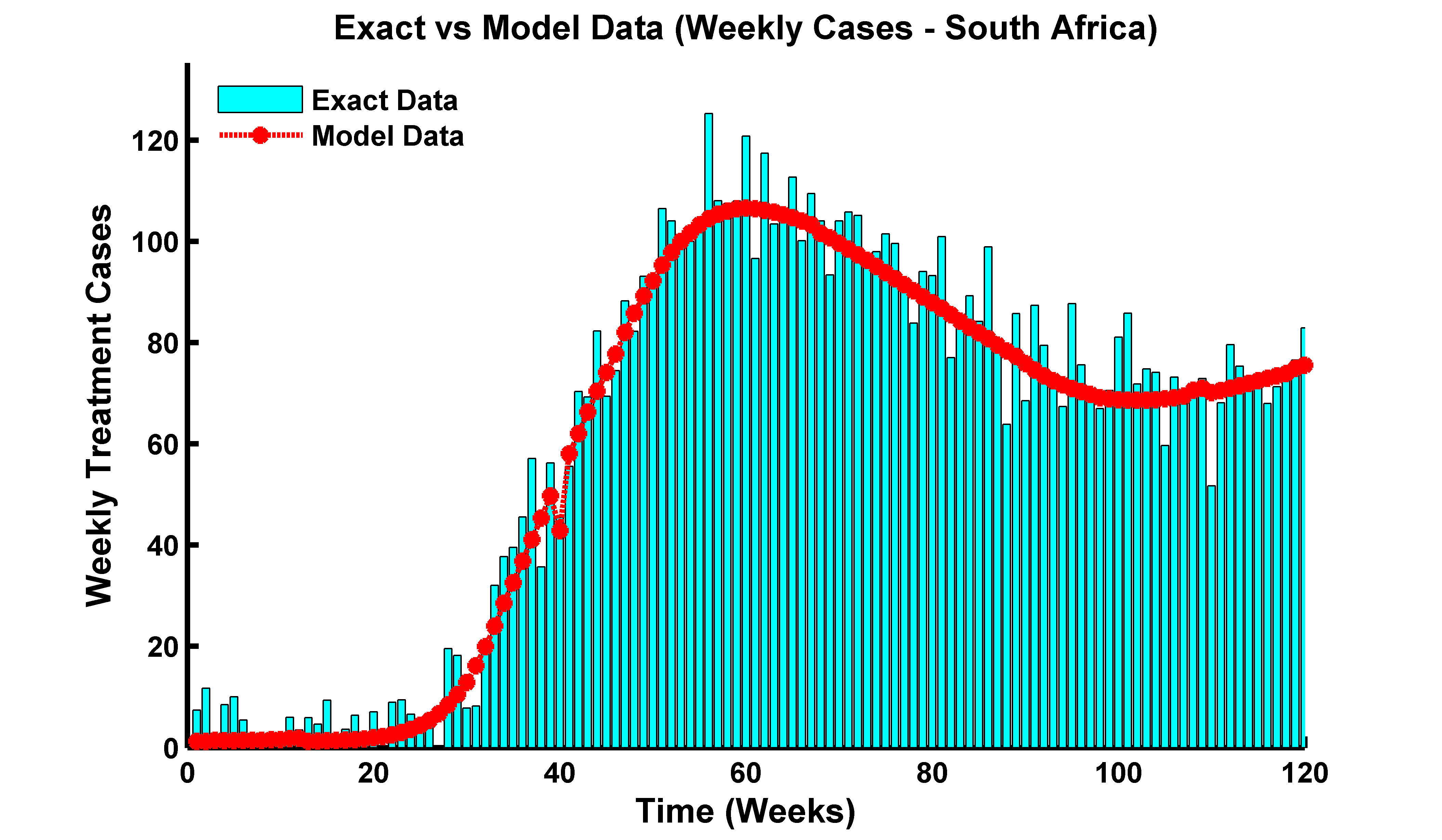}}\\
	\subfloat[]{\includegraphics[width=3. in]{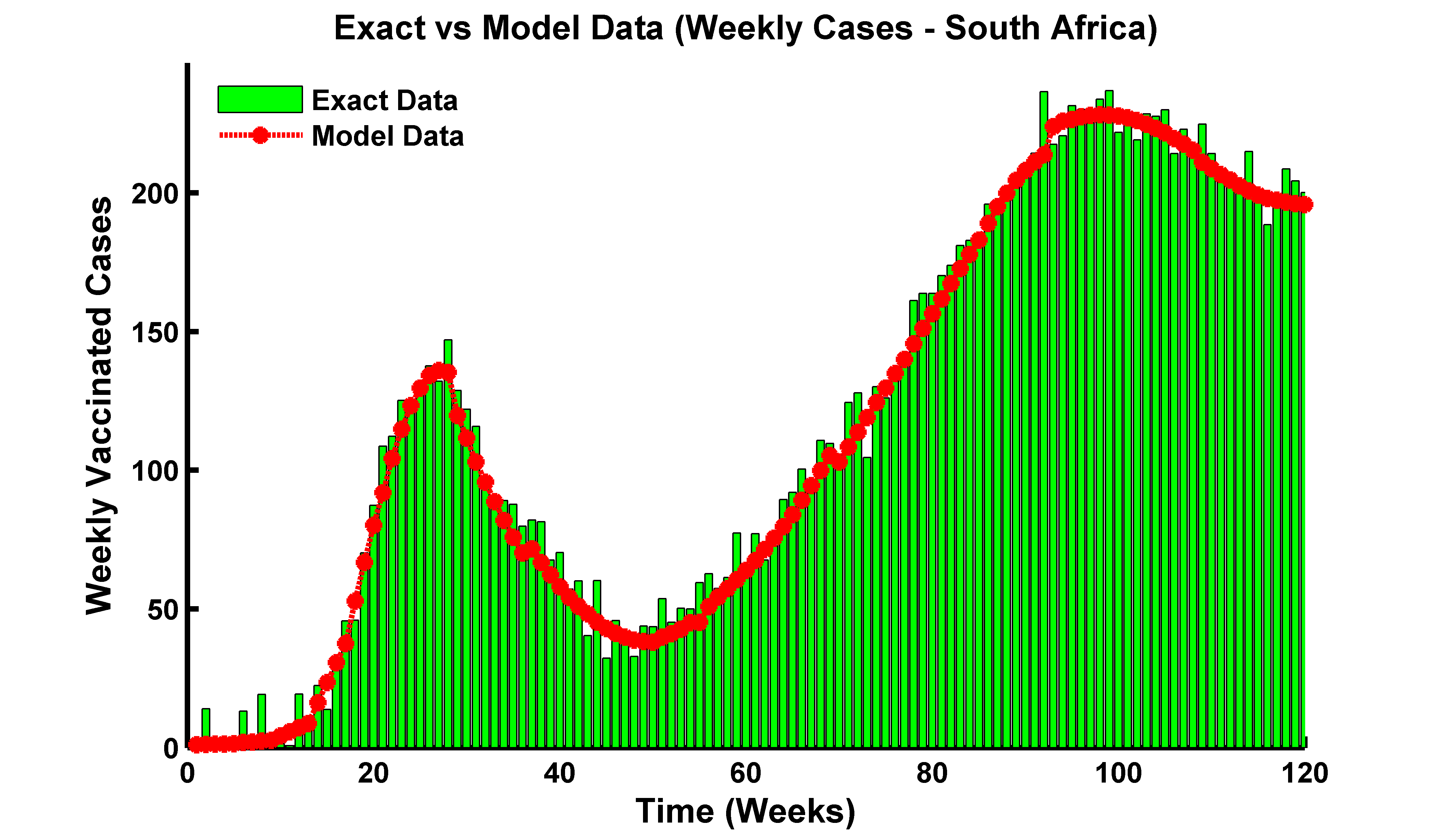}}
	\subfloat[]{\includegraphics[width=3. in]{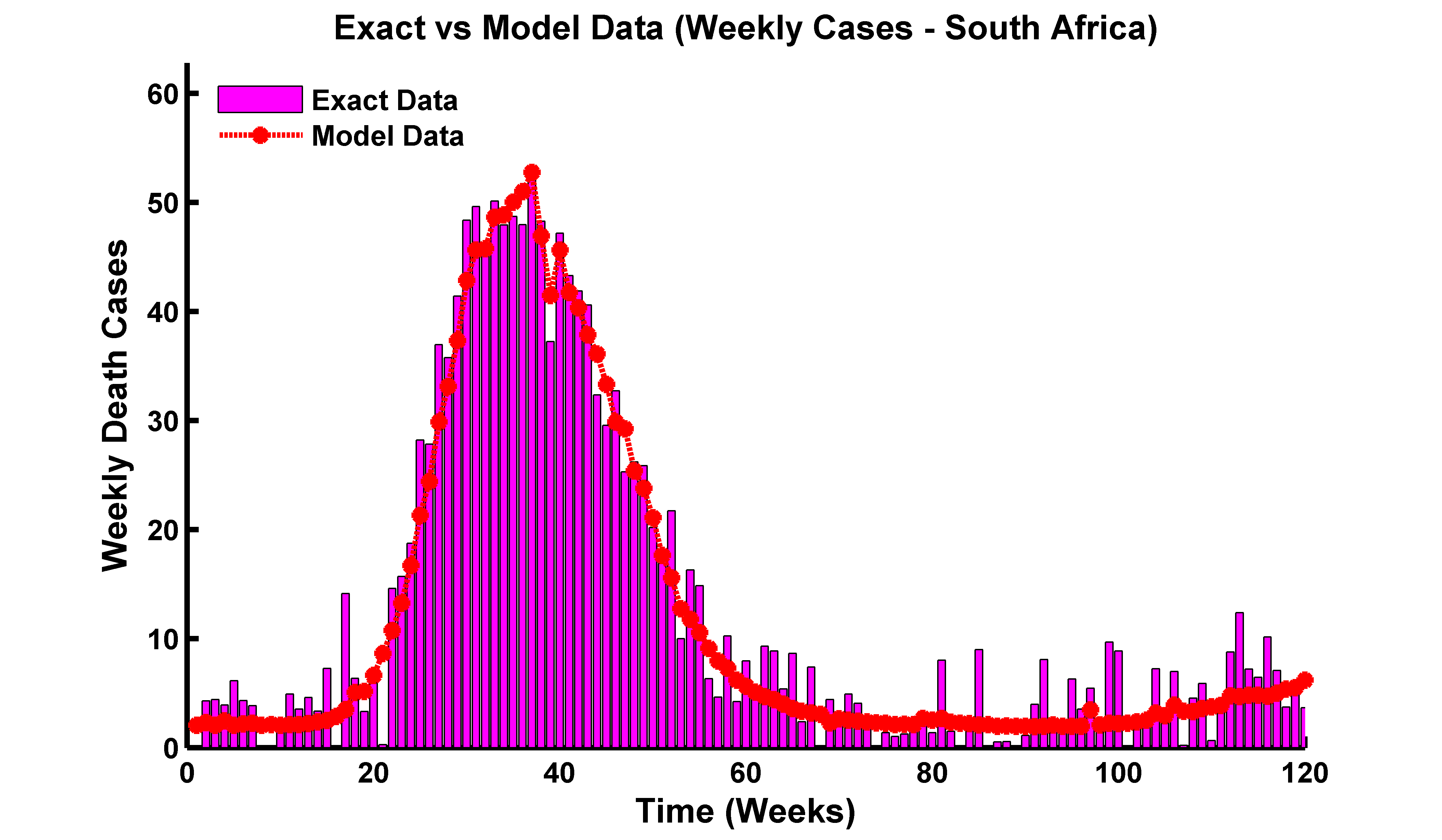}}
	\caption{Model data fitting to South Africa with (a) weekly infected data (b) weekly progression of treatment to control disease burden (c) weekly vaccination programme to control disease burden (d) weekly death rates, where the fitted parameter sets are : $\beta_1=0.0053, \beta_2= 0.0061, \alpha=0.67, \gamma=0.61, \gamma_1=0.31, \lambda=0.52, \mu=0.03, \delta=0.27$.}
	\label{case-study-South-Africa}
\end{figure}

In order to study the recent trend with our suggested model, we have examined recent data of South Africa for influenza infection cases, starting from 1 October 2020 to 31 March 2023. The CDC and WHO websites were used to acquire all of the data \cite{CDC,WHO}. We have taken into account a total of 120 weeks' worth of data in order to assess the model's applicability and estimate disease for the near future. The initial population is taken as $S(0)=500$, $V(0)=1$, $E(0)=1$, $I(0)=1$, $R(0)=1$ and $T(0)=0$. The model outcomes predicted with real data presented in 
Figure \ref{case-study-South-Africa}, where the weekly reported cases in South Africa are depicted. The linear regression analysis with parameters in Table \ref{table-param-values-4-countries}, suggests that there is remarkable agreement between the model and the weekly infected cases. Thus, we demonstrate approximately 76\% accuracy to track the original confirmed cases. As South Africa is located at the southernmost tip of the African continent; generally, South Africa experiences summer from December to February and winter from June to August. During summer, temperatures can range from warm to hot, with coastal areas experiencing pleasant temperatures and inland regions being hotter. Winter temperatures are generally milder. CDC reported that: the first time the impact of influenza in South Africa; the illness causes 40,000 hospital admissions and about 10,000 fatalities annually. In the meanwhile, the 2022 flu season in South Africa began the week of April 25 and the number of cases seems to be rising over the last few weeks.

Beginning in 2021, there have been documented occurrences of influenza in South Africa of 5.8\% and deaths of 1.7\%. The authorities made the initial decision to implement the immunization program, which is still being done today in hospitals, community vaccination centers, clinics, and pharmacies. After receiving the recommended dosage, the immunizations were effective. Additional vaccines that have been authorized include Fluarix, FluLaval, Fluzone, and Afluria quadrivalent. These vaccinations are safe to provide to infants as young as six months. These vaccines have an efficacy of roughly 66\%, 64\%, 67\%, and 62\%, respectively. The Figures \ref{case-study-South-Africa}(a)-(d) illustrate the model prediction which correlates with the actual data. The outcomes show that the suggested model is applicable.

Analyzing for total cases of population per 1000, from the Figure \ref{case-study-South-Africa}(a) we see that, weekly reported cases of seasonal flu in South Africa is in small portion starting from 0 to 20 weeks. At week 20, reported case is 15.From week 21 to week 44 after that, the infection gradually got worse. The highest number of infected cases (120 per 1000) were reported in South Africa between weeks 33 and 40. Along weeks 30 to 50, disease conditions are really bad. Due to appropriate care and isolation, the illness epidemic marginally decreased after 50 weeks. After 75 weeks, the reported incidences vary (per 1000) from 20 to 42.

In Figure \ref{case-study-South-Africa}(b), the weekly treatment cases are shown. Which suggests that as the number of infected cases increased, treatment standards climbed. Between 85 and 120 treatment instances (per 1000) are documented from weeks 42 to 80. After 80 weeks, this rate ranged over time, ranging from 55 to 75 (per 100).
\noindent
In Figure \ref{case-study-South-Africa}(c), with our suggested approach, the ongoing vaccination cases are projected. Our model predicts that the immunization rate will peak between 20 and 35 weeks. A preventive policy was applied to 130 populations (per 1000) on average. After then, the number of vaccination cases increased along with the number of infected patients. The range of reported immunization cases 170–210 weeks is crucial for managing the epidemic.\\
On the other hand, Figure \ref{case-study-South-Africa}(d) explains the current South African fatality cases. According to data analysis, from October 1, 2021, the weeks 30 to 50 have the greatest death cases. The death rate is approximately 0 to 20 (per 1000) from week 0 to week 20. After 60 weeks, the number of deaths is under control and ranges from 0 to 10 (per 1000). After doing a data analysis, we discovered that South Africa's death rates fluctuate quickly.\\
Figures   \ref{case-study-Italy}, and \ref{case-study-South-Africa} shows that, insightful conclusions were drawn from our model's analysis of the influenza cases in  Italy, and South Africa. The information paints a complete picture of the influenza scenario and offers useful details for comprehending the virus's propagation and effects. We were able to make precise forecasts and spot important patterns and trends by fitting the data to our model. These findings can help academics, politicians, and public health professionals develop efficient plans to fight influenza and lessen its effects on the populace. In order to support ongoing efforts in public health and disease surveillance, the results of our model provide a compelling and instructive account of the influenza figures in Mexico, Italy,  and South Africa.

\section{Conclusion}\label{Section-Concluding-Remarks}
We have considered epidemic control model for influenza infection by taking vaccination and treatment compartments. The model's existence, positivity, and boundlessness are verified. 
We have performed relative bias, relative influence analysis to present several scenarios and impact on the basic reproduction number. We have looked at both local and global sensitivity analyses for every model parameter. We ultimately determine that the most sensitive parameter for containing disease outbreaks is the contact rates $\beta_1$, $\beta_2$, infection rates $\alpha$, recovery and treatment rates $\gamma$, $\gamma_1$. This analysis is done after taking into account PRCC values and p-values. That means, by controlling these parameters in minimal range, the spread of influenza transmission can be mitigated. Additionally, the numerical and analytical findings showed that, after a predetermined amount of time, the number of exposed and infected people will steadily decline. Therefore, our ideal model guarantees the efficacy of the plan to lessen and eventually stop illness outbreaks. The numerical results supports the analytical results and verified that 
the illness will continue to exist in the community because the endemic equilibrium is globally asymptotically stable. On the other hand, 
the equilibrium point that is devoid of illness is globally asymptotically stable, and the disease will eventually vanish from the community. More significantly, the model is optimized in this work as we used certain suitable control measures and examined optimum control theory. To investigate the possibility of optimum control, Hamiltonian and Lagrangian formulations are utilized, and Pontryagin's maximum principle is applied to determine the essential optimality criteria. We used a variety of control techniques to try and limit the infection, but we found that implementing treatments all at once is a more effective way to lower the overall disease burden. Further, we have analyzed recent trends in Mexico, Italy, and South Africa data in which analytical, numerical and statistical studies were implemented, respectively. We successfully depicted the respective disease incidence by showing numerical results with proper fittings. We have validated our framework by comparing its predictions with simulation results. Ultimately, vaccinations, treatment plans, maintaining a preventative lifestyle, and ensuring that everyone has access to a sufficiently supported medical care system may all significantly minimize the likelihood of an illness epidemic.
We can observe from the model's simulations that if the transmission rates $\beta_1$, $\beta_2$ and rate of acquiring infection $\alpha$, increases that results into the rapid increase of the disease burden. Further simulation of the model suggests that, the treatment and vaccination to the exposed individuals can reduce the disease burden significantly. Quarantine policy, using masks, hand sanitizer, droplets in a proper way is important to increase the recovery and reduction of interaction rate.  Results of  Relative bias, According to a local and global sensitivity analysis of the threshold quantity $\mathcal{R}_0$, the contact rates $\beta_1$, $\beta_2$, infection rate $\alpha$, recovery rate $\gamma$, and treatment rate $\gamma_1$ are the most important factors. Reducing the value of these factors can help to stop the spread of the illness. According to the CDC, the vaccinations that have proven most successful in treating influenza interactions thus far are Afluria, Fluarix, FluLaval, and Fluzone. Therefore, vaccination programme should be maintained properly to decrease the levels of exposed and infected compartments.

\begin{appendices}  
	\appendix
	\section*{Appendix}
	\renewcommand\thefigure{\thesection.\arabic{figure}}  
	\renewcommand\thetable{\thesection.\arabic{table}}  
	\section{ Mathematical Model and Auxiliary Results}\label{appen}
	The appendix contains the model description and all preliminary results utilized in the main body of the paper.
	
\subsection{Model Formulation}\label{Subsection-Model Formulation}
When it comes to infectious disease modeling, the traditional SIR model determines the population's critical state for the development of the illness while taking the size of the entire population into account \cite{Marcheva Book,Stability Bound-11}. Influenza, categorized as a person-to-person transmissible disease, is of particular interest. In this study, we introduce a potential SVEIRT mathematical model comprising six compartments.
	\begin{align}\label{new_model}
		\begin{cases}
			\vspace{0.2cm}
			\displaystyle\frac{dS}{dt} = \Lambda - \left(\beta_1E+\beta_2I\right)S-(\mu+\phi) S \\
			\vspace{0.2cm}
			\displaystyle\frac{dV}{dt} = \phi S-(1-\varepsilon)\left(\beta_1E+\beta_2I\right) V-\mu V\\
			\vspace{0.2cm}
			\displaystyle\frac{dE}{dt} =\left(\beta_1E+\beta_2I\right)S-(\alpha+\mu) E\\
			\vspace{0.2cm}
			\displaystyle\frac{dI}{dt} =\alpha E+ (1-\varepsilon)\left(\beta_1E+\beta_2I\right) V-(\mu+\delta+\gamma+\gamma_1) I\\
			\vspace{0.2cm}
			\displaystyle\frac{dR}{dt} = \gamma I-\mu R\\
			\displaystyle\frac{dT}{dt} =\gamma_1 I-\mu T, 
		\end{cases}
	\end{align}
	for $ t \in (0,\infty) $ under the starting circumstances
	\begin{align}\label{ic}
		S(0) = S_{0},\;\;V(0) = V_0. \;\;  E(0) = E_0, \;\;  I(0) = I_0, \;\;  R(0) = R_0, \;\; \text{and}\;  T(0) = T_0.
	\end{align}

	\begin{figure}[H]
		\centering
		\begin{tikzpicture}[node distance=2cm]
			\node (cs) [cs] {$S$};
			\node (cv) [cv,below of=cs,yshift=-1.0cm,xshift=0cm] {$V$};
			\node (ce) [ce,right of=cs,xshift=2.0cm] {$E$};
			\node (ci) [ci,right of=ce,xshift=2.0cm] {$I$};
			\node (ct) [ct,below of=ci,yshift=-1.0cm,xshift=0cm] {$T$};
			\node (cr) [cr, right of=ci,xshift=2.0cm] {$R$};
			\node (a) [draw=none, fill=none, left of=cs,yshift=-.0cm] {};
			\node (b) [draw=none, fill=none, above of=cs,yshift=0cm] {};
			\node (c) [draw=none, fill=none, above of=ce,yshift=0cm] {};
			\node (d) [draw=none, fill=none, above of=ci,yshift=0cm] {};
			\node (e) [draw=none, fill=none, above of=cr,yshift=0cm] {};
			\node (f) [draw=none, fill=none, below of=cv,yshift=0cm] {};
			\node (g) [draw=none, fill=none, below of=ct,yshift=0cm] {};
			\node (A) [draw=none, fill=none,left of=cs,yshift=.4cm,xshift=1.0cm]{$\Lambda$};
			\node (B) [draw=none, fill=none,left of=cs,yshift=1.2cm,xshift=2.6cm]{$\mu S$};
			\node (C) [draw=none, fill=none,left of=ce,yshift=.5cm,xshift=-.1cm]{$(\beta_1 E+\beta_2 I)S$};
			\node (D) [draw=none, fill=none,right of=cs,yshift=-1.2cm,xshift=-2.6cm]{$\phi S$};
			\node (E) [draw=none, fill=none,right of=ce,yshift=.5cm,xshift=-.1cm]{$\alpha E$};
			\node (F) [draw=none, fill=none,left of=ce,yshift=1.2cm,xshift=2.6cm]{$\mu E$};
			\node (G) [draw=none, fill=none,left of=ci,yshift=1.2cm,xshift=2.9cm]{$(\mu + \delta) I$};
			\node (H) [draw=none, fill=none,left of=cr,yshift=1.2cm,xshift=2.6cm]{$\mu R$};
			\node (I) [draw=none, fill=none,right of=cv,yshift=-1.2cm,xshift=-2.6cm]{$\mu V$};
			\node (J) [draw=none, fill=none,right of=ct,yshift=-1.2cm,xshift=-2.6cm]{$\mu T$};
			\node (K) [draw=none, fill=none,right of=ci,yshift=.5cm,xshift=-.1cm]{$\gamma I$};
			\node (L) [draw=none, fill=none,right of=ci,yshift=-1.2cm,xshift=-2.6cm]{$\gamma_1 I$};
			\node (M) [draw=none, fill=none,right of=cv,yshift=-.02cm,xshift=1.1cm]{$(1-\varepsilon)(\beta_1 E+\beta_2 I)V$};
			\draw [line width=0.55mm, blue,dashed,->] (cs)--(ce);
			\draw [line width=0.55mm, blue,dashed,->] (ce)--(ci);
			\draw [line width=0.55mm, blue,dashed,->] (ci)--(cr);
			\draw [line width=0.55mm, blue,dashed,->] (cs)--(cv);
			\draw [line width=0.55mm, blue,dashed,->] (ci)--(ct);
			\draw [line width=0.55mm, blue,dashed,->] (cv)--(ci);
			
			\draw [line width=0.55mm,blue,dashed,->] (a)--(cs);
			\draw [line width=0.55mm,blue,dashed,->] (cs)--(b);
			\draw [line width=0.55mm,blue,dashed,->] (ce)--(c);
			\draw [line width=0.55mm,blue,dashed,->] (ci)--(d);
			\draw [line width=0.55mm,blue,dashed,->] (cr)--(e);
			\draw [line width=0.55mm,blue,dashed,->] (cv)--(f);
			\draw [line width=0.55mm,blue,dashed,->] (ct)--(g);
		\end{tikzpicture}
		\caption{Flow diagram of the model.}\label{fignewmodel}  
	\end{figure}
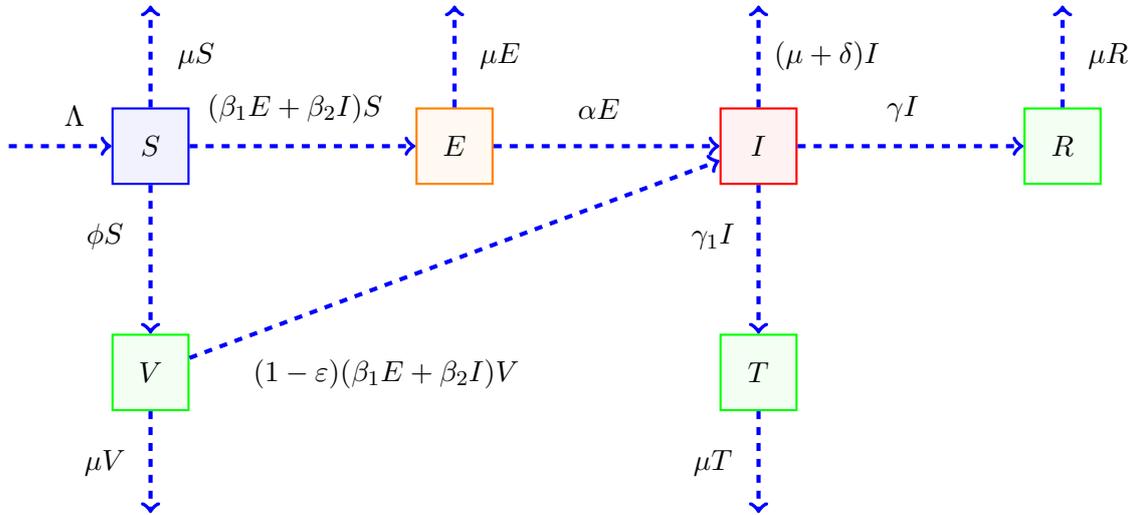
	\noindent	
	Six mutually-exclusive classes are created here from the human population, where $ S(t),\;V(t),\;E(t)$, $\;I(t),\;R(t) $, and $ T(t) $ denote the number of people in the susceptible, vaccinated, exposed, infected, recovered, and treated compartments, respectively, at time $ t $, expressed as a unit per week. Table~\ref{description-params-states} provides a concise explanation of each definition for state variable and model parameter.  In the Figure~\ref{fignewmodel} The illness transfer from one compartment to the other is shown by the black arrows.
	\begin{table}[H]
		\begin{center}
			\caption{Model state variables and parameters with their descriptions.}
			\scriptsize
			\label{description-params-states}
			\begin{tabular}{|l|l|}
				\hline\noalign{\smallskip}
				\textbf{Notation} & \textbf{Definition} \\
				\noalign{\smallskip}\hline\noalign{\smallskip}
				$N$ & Total number of human population\\
				$S$ & Total number of susceptible population \\
				$V$ & Total number of vaccinated population \\
				$E$ & Total number of exposed population \\
				$I$ & Total number of infected population \\
				$R$ & Total number of recovered population \\
				$T$ & Total number of treated population \\
				$\Lambda$ & Recruitment/birth rate in $S$ class\\
				$\beta_1$ & Transmission probability of disease from contact with $S$ and $E$ \\
				$\beta_2$ & Transmission probability of disease from contact with $S$ and $I$\\
				$\phi$ & Rate of vaccination to the $S$ population\\
				$\alpha$ & Transmission rate from $E$ to $I$ class. \\
				$\gamma$ & Recovery rate from $I$ class.\\
				$\gamma_1$ & Treatment rate of $I$ class to the treated class $T$.\\
				$\mu$ & Natural death rate (per unit of time)\\
				$\delta$ & Disease induced death rate in $I$ compartment(per unit of time)\\
				$\varepsilon$ & The effectiveness of vaccine\\
				\noalign{\smallskip}\hline		
			\end{tabular}
		\end{center}
	\end{table}

	\section*{Determination of Fixed Points and Threshold numbers}\label{Section-Determination Fixed Points}
	The model \eqref{new_model} possess two different kinds of equilibrium: endemic equilibrium (EE) and disease free equilibrium (DFE). To determine the points of equilibrium $(\widetilde{S},\widetilde{V},\widetilde{E},\widetilde{I},\widetilde{R},\widetilde{T})$ of the system \eqref{new_model}, We set the value of each derivative to zero. 
	Replacing the variables for the disease-free equilibrium (DFE) is as follows:
	$$ (\widetilde{S},\widetilde{V},\widetilde{E},\widetilde{I},\widetilde{R},\widetilde{T}) \equiv (S_{0},V_{0},E_{0},I_{0},R_{0},T_{0})$$
	such that
	\begin{align}
		\left( S_0,V_0,E_0,I_0,R_0,T_0 \right) \equiv \left( \dfrac{\mu N}{\mu + \phi},\dfrac{\phi N}{\mu + \phi}, 0, 0, 0,0 \right).
	\end{align}
	
\noindent	The variables for the endemic equilibrium (EE) are substituted as \\
	$(\widetilde{S},\widetilde{V},\widetilde{E},\widetilde{I},\widetilde{R},\widetilde{T}) \equiv (S^{*},V^{*},E^{*},I^{*},R^{*},T^{*}), $ where, $ E^{*} > 0 \; \text{and} \; I^{*} > 0$\ also $E^{*}\neq 0,\; I^{*}\neq 0$.
 Here,
	$S^{*}=\frac{\displaystyle \Lambda-(\alpha+\mu)E^{*}}{\displaystyle (\mu+\phi)}$, $V^{*}= \frac{\displaystyle \phi(\Lambda-a_1 E^{*})}{\displaystyle a_2(\mu+\lambda\lambda_1)},$ $R^* = \frac{\displaystyle \lambda I^*}{\displaystyle \mu}, \;
	T^* = \frac{\displaystyle \lambda_1 I^*}{\displaystyle \mu},$\\
	$E^*=\frac{\displaystyle (\Lambda \beta_1-a_1a_2-a_1\beta_2 I^{*})\pm \sqrt{(\Lambda \beta_1-a_1a_2-a_1\beta_2 I^{*})^2 + 4\Lambda\beta_2I^* a_1\beta_1}}{\displaystyle 2a_1\beta_1}.$\\ Where, $a_1=\alpha+\mu, \; a_2=\mu+\phi, \; \lambda_1=\beta_1E^*+\beta_2I^*, \;\text{and}\; \lambda=1-\epsilon.$ 	The details analysis of EE can be found in \cite{MK24}.
	
	\subsection*{Basic Reproduction Number with Control}\label{Subsection-R0-with-Control}
	The threshold value for the system (\ref{new_model}) associated with controlling strategies can be presented as follows,
	\begin{align} \label{repro1}
		\mathcal{R}_{0V} 
		&=\frac{\mu N\left[\alpha\beta_2+\beta_1(\mu+\delta+\gamma+\gamma_1)\right]}{(\mu+\phi)(\alpha+\mu)(\mu+\delta+\gamma+\gamma_1)}+\frac{N\phi\beta_2\lambda}{(\mu+\phi)(\mu+\delta+\gamma+\gamma_1)}.
	\end{align}
	
	\subsection*{Basic Reproduction Number without Control}\label{Subsection-R0-without control}
	The basic reproduction number of the apparatus (\ref{new_model}) in the absence of control strategies can be expressed as follows,
	\begin{align} \label{repro2}
		\mathcal{R}_{0} &=\frac{S_0\left[\alpha\beta_2+\beta_1(\gamma+\gamma_1+\mu+\delta)\right]}{(\alpha+\mu)(\gamma+\gamma_1+\delta+\mu)}.
	\end{align}

	\subsection{Qualitative Study of the Model}\label{qaulity}
	In order to maintain the biological validity of the model, it is now necessary to demonstrate that solutions to the system of differential equations \eqref{op_syseqn} exist, are limited, and are positive for all values of time.
	\begin{Th}
		(Positivity). Let $t_0>0$, Within the model \eqref{op_syseqn}, if the preconditions are met $S_0>0,V_0>0,\;E_0>0,\;I_0>0,\;I_0>0,\;R_0>0$ and $T_0>0$, then $\forall\;t\in[0,t_0]$ the functions $S(t),\;V(t),\;E(t),\;I(t),\;R(t)$ and $T(t)$will stay in $\mathbb{R}_+^6.$
	\end{Th}
	\begin{proof}
		Lower bounds may be set for each of the equations given in the model \eqref{op_syseqn} because all of the parameters used in the system are positive. Consequently, we have
		\begin{align*}
			\frac{dS}{dt} =&\Lambda-(1-w_1(t))(\beta_1E+\beta_2I)S-(\mu+\phi)S
			\geq  -(1-w_1(t))(\beta_1E+\beta_2I)S-(\mu+\phi)S\\
			\frac{dV}{dt}=&\phi S-(1-\varepsilon)(\beta_1E+\beta_2I)V-\mu V 
			\geq -(1-\varepsilon)(\beta_1E+\beta_2I)V-\mu V\\
			\frac{dE}{dt}=&(1-w_1(t))(\beta_1E+\beta_2I)S-(\alpha+\mu)E \geq -(\alpha+\mu)E\\
			\frac{dI}{dt}=&\alpha E+(1-\varepsilon)(\beta_1E+\beta_2I)V-(\mu+\delta)I-(1+w_2(t))\gamma_1 I-(1+w_3(t))\gamma I\\
			\geq& -(1+w_2(t))\gamma_1 I-(1+w_3(t))\gamma I\\
			\frac{dR}{dt}=&\gamma(1+w_3(t))I-\mu R \geq -\mu R\\
			\frac{dT}{dt}=&\gamma_1(1+w_2(t))I-\mu T \geq -\mu T
		\end{align*}
	\end{proof}
	Applying basic differential equations methods, the inequalities can be resolved and the procedure follows:
	\begin{align*}
		\frac{dS}{dt}\geq &S(0)e^{-(\mu+\phi)t-\int_{0}^{T_f}(1-w_1(t))(\beta_1E+\beta_2I)dt} >0\\
		\frac{dV}{dt}\geq &V(0)e^{-\mu t-(1-\varepsilon)\int_{0}^{T_f}(\beta_1E+\beta_2I)dt} >0\\
		\frac{dE}{dt}\geq &E(0)e^{-(\mu+\phi)t}>0\\
		\frac{dI}{dt}\geq&I(0)e^{-\{(1+w_2(t))\gamma_1+(1+w_3(t))\gamma+(\mu+\delta)\}t}>0\\
		\frac{dR}{dt}\geq&R(0)e^{-\mu t}>0,\;\;\; \frac{dT}{dt}\geq T(0)e^{-\mu t}>0
	\end{align*}
	That indicates $\forall \;t\in[0,t_0]$, the roles $S(t),\;V(t),\;E(t),\;I(t),\;R(t)$ and $T(t)$ will continue to be positive and stay in $\mathbb{R}_+^6.$\\
	To investigate the existence of an optimal control of our proposed model, it is necessary to establish the boundedness of solutions to system \eqref{op_syseqn} for finite time span. We are now investigating the state solutions' a priori boundedness \cite{Optimal Control-4, Optimal Control-6}.
	\begin{Th}
		(Boundedness).Considering, $(w_1,w_2,w_3)\in W$, there exists bounded solutions for the problem \eqref{op_syseqn}.
	\end{Th}
	\begin{proof}
		Here, we consider state variables that are super-solutions to the problems \eqref{op_syseqn}. Using the provided equations as a guide, we have derived,
		\begin{align*}
			(S+V+E+I+T+R)'(t)=\Lambda-\mu(S+V+E+I+R+T)-\delta I \leq \Lambda-\mu N
		\end{align*}
		Now, using $N(t)=S(t)+V(t)+E(t)+I(t)+R(t)+T(t)$ and $\mu>0$, we get, $N(t)'\leq \Lambda-\mu N,$\\
		which suggests that,
		\begin{align*}
			\limsup\limits_{t\rightarrow\infty} N(t)<\frac{\Lambda}{\mu}.
		\end{align*}
		This illustrates that the upper bound for $N$ acts in the same way as the upper bound for $S,V,E,I,R$, and $T$. At last,
		\begin{align*}
			T'(t)=\gamma_1(1+w_2(t))I-\mu T\leq \gamma_1(1+w_2(t))I\leq \frac{\gamma_1(1+w_2(t))\Lambda}{\mu}
		\end{align*}
		which leads to
		\begin{align*}
			T(t)\leq \frac{\gamma_1(1+w_2(t))\Lambda T_f}{\mu} \in \mathbb{R}_+,\;\;\text{for all }t\in[0,T_f].
		\end{align*}
		Since, $(w_1(t),w_2(t),w_3(t),w_4(t),w_5(t),w_6(t))\in W$, then along with $S(t),\;V(t),\;E(t),\;I(t),\;R(t)$ and $T(t)$ are bounded above. Solutions to the problems \eqref{op_syseqn} have been obtained through the use of a maximal principle theory for first order non-linear differential equations, that are bounded $\forall \;t\in[0,t_0]$ and stays in the compact set
		\begin{align*}
			\mathbb{D}=\left\{(S,V,E,I,R,T)\in\mathbb{R}_{+}^{6}: S,V,E,I,R,T\leq \frac{\Lambda}{\mu}, \;\;T\leq \frac{\gamma_1(1+w_2(t))\Lambda T_f}{\mu}\right\}
		\end{align*}
		where, $\displaystyle \mathbb{R}_{+}^6=\left\{(S,V,E,I,R,T): S\geq 0, V\geq 0,E\geq 0,I\geq 0,R\geq 0,T\geq 0\right\}.$
		Thus completes the proof.
	\end{proof}
	
	\begin{Th}
		(A solution's existence) Let, $t_0>0$, In the model \eqref{op_syseqn}, if the basic requirements are met $S_0>0,\;V_0>0,\;E_0>0,\;I_0>0,\;R_0>0$ and $T_0>0$, then $\forall \;t\in\mathbb{R}$ the functions $S(t),V(t),E(t),I(t),R(t)$ and $T(t)$ will exist in $\mathbb{R}_{+}^6.$
	\end{Th}
	\begin{proof}
		In the instance of our model \eqref{op_syseqn}, the following function defines the ODE system\\
		$f:\mathbb{R}^6\rightarrow\mathbb{R}^6$ as follows:
		$\mathbf{f(y)}=(f_1,f_2,f_3,f_4,f_5,f_6)^T$ where,
		\begin{align*}
			&f_1=\Lambda-(1-w_1(t))(\beta_1E+\beta_2I)S-(\mu+\phi)S\\
			&f_2=\phi S-(1-\varepsilon)(\beta_1E+\beta_2I)V-\mu V\\
			&f_3=(\beta_1E+\beta_2I)(1-w_1(t))S-(\alpha+\mu)E\\
			&f_4=\alpha E+(1-\varepsilon)(\beta_1E+\beta_2I)V-(\mu+\delta)I-\gamma_1(1+w_2(t))I-\gamma(1+w_3(t))I\\
			&f_5=\gamma(1+w_3(t))I-\mu R\\
			&f_6=\gamma_1(1+w_2(t))I-\mu T
		\end{align*}
		It is observed that on $\mathbb{R}^6$, $\mathbf{f}$ has a continuous derivative. Thus, in $\mathbb{R}^6$, $\mathbf{f}$ is locally Lipschitz. Therefore, we have examined if there is a unique, positive, and bounded solution to the ordinary differential equations presented in \eqref{op_syseqn}. This is because of the Fundamental Existence and Uniqueness Theorem as well as the Theorems demonstrated on positivity and boundedness of solutions \cite{Optimal Control-5}. 
	\end{proof}

\end{appendices}

\section*{Acknowledgments}
The research by M. Kamrujjaman was partially supported by the University Grants Commission (UGC), 
and the  University of Dhaka, Bangladesh.

\section*{Author contributions}
Conceptualization,  MK and KMM; methodology, KMM; software, KMM,  and MK; validation, MK; formal analysis, KMM and MK; investigation, MK; resources, MK;  data curation, KMM; original draft preparation, KMM and MK; review and editing,   MK; supervision, MK. All authors have read and agreed to the published version of the manuscript.

\section*{Competing interests/Conflict of interest}
The authors declare no conflict of interest. 

\section*{Data availability/Data sharing}
The  data was collected form the following links:\\
\url{https://gis.cdc.gov/grasp/fluview/fluportaldashboard.html}\\
\url{https://www.who.int/publications/m/item/influenza-update-n-389}\\
The used data will be published as csv or excel file upon acceptance of the manuscript.
We also converted the weekly data into units of thousands. 		

\section*{Ethical approval}
No consent is required to publish this manuscript.

\section*{Ethics Statement} 
None.


\end{document}